\documentclass[10pt,a4paper,reqno]{amsart}
\usepackage{amsfonts}
\usepackage{amsthm}
\usepackage{amsmath}
\usepackage{mathtools}
\usepackage{amscd}
\usepackage[utf8]{inputenc}
\usepackage{t1enc}
\usepackage[mathscr]{eucal}
\usepackage{indentfirst}
\usepackage{graphicx}
\usepackage{graphics}
\usepackage{esint}
\usepackage{pict2e}
\usepackage{epic}
\usepackage{float}
\usepackage{MnSymbol}
\usepackage{multirow}
\usepackage{shuffle} 
\usepackage[table,xcdraw]{xcolor} 
\usepackage[nocompress]{cite}
\usepackage{hyperref}
\usepackage{tcolorbox}
\numberwithin{equation}{section}
\usepackage[margin=2.9cm]{geometry}
\usepackage{enumitem}
\usepackage{epstopdf} 
\usepackage{fancyvrb} 
\usepackage{tikz}
\usetikzlibrary{decorations.markings,calc}

\usepackage[dvipsnames]{xcolor}

\def\centerarc[#1](#2)(#3:#4:#5)
{ \draw[#1] ($(#2)+({#5*cos(#3)},{#5*sin(#3)})$) arc (#3:#4:#5); }

\def\centerarcpath(#1)(#2:#3:#4)
{ ($(#1)+({#4*cos(#2)},{#4*sin(#2)})$) arc (#2:#3:#4) }
 
\allowdisplaybreaks

\VerbatimFootnotes

\newmuskip\pFqmuskip
\newcommand*\pFq[6][8]{%
	\begingroup 
	\pFqmuskip=#1mu\relax
	\mathchardef\normalcomma=\mathcode`,
	\mathcode`\,=\string"8000
	\begingroup\lccode`\~=`\,
	\lowercase{\endgroup\let~}\pFqcomma
	{}_{#2}F_{#3}{\left[\genfrac..{0pt}{}{#4}{#5};#6\right]}%
	\endgroup
}
\newcommand{\pFqcomma}{{\normalcomma}\mskip\pFqmuskip}

\theoremstyle{plain}
\newtheorem{theorem}{Theorem}[section]
\newtheorem{lemma}[theorem]{Lemma}
\newtheorem{corollary}[theorem]{Corollary}
\newtheorem{proposition}[theorem]{Proposition}

\theoremstyle{definition}

\newtheorem{conjecture}[theorem]{Conjecture}
\newtheorem{remark}[theorem]{Remark}

\newtheorem{example}[theorem]{Example}

\newcommand{\colorn}[1]{ {\color{blue} n^{#1}} }
\DeclareMathOperator{\Li}{Li}
\newcommand{\Res}[1]{\underset{#1}{\operatorname{\: Res \:}}}
\newcommand{\mycomment}[1]{}

\renewcommand{\Re}{\operatorname{Re}}
\renewcommand{\Im}{\operatorname{Im}}

\begin{document}
	
	\title{On single-variable Witten zeta functions of rank two and three}
	
	\author[Kam Cheong Au]{Kam Cheong Au}
	
	\address{University of Cologne \\ Department of Mathematics and Computer Science \\ Weyertal 86-90, 50931 Cologne, Germany} 
	
	\email{kau@uni-koeln.de}
	

	\begin{abstract} By introducing a novel integration kernel for the Mellin transform, we uncover many previously unknown and intriguing properties of the Witten zeta functions of rank two and three. Detailed results concerning their pole locations, residues, and special values are obtained. We propose a non-trivial conjecture regarding their derivatives at the origin, which seems to encode deep information about the root system. We also discuss their behavior at negative integers, highlighting a connection with Eisenstein series and a $p$-adic observation. 
	\end{abstract}

	\maketitle
	
	
	\section{Introduction}
	For a finite-dimensional simple Lie algebra $\mathfrak{g}$ over $\mathbb{C}$, its \emph{Witten zeta function} is the Dirichlet series $$\zeta_\mathfrak{g}(s) := \sum_{\rho} \frac{1}{(\text{dim } \rho)^s},\quad \rho\in \{\text{finite-dimensional irreducible representations of $\mathfrak{g}$}\}.$$
	If $\Phi$ is the root system associated with $\mathfrak{g}$, then we abuse notation by writing $\zeta_\mathfrak{g}(s)$ as $\zeta_\Phi(s)$. 
	
	There is a vast literature on Witten zeta functions \cite{matsumoto2006witten, komori2023theory, onodera2014functional, komori2012witten, komori2015witten, komori2010witten, okamoto2012multiple}, especially on the meromorphic continuation and functional relations of their multi-variable generalizations. We focus on Witten zeta functions for rank $2$ and $3$ irreducible root systems, namely
	$$\zeta_\Phi(s) \text{ with }\Phi \in \{A_2,B_2,G_2,A_3,B_3,C_3\}.$$
	For example,
	$$\begin{aligned}
		\zeta_{A_2}(s) &= 2^s \sum_{n,m\geq 1} \frac{1}{m^sn^s(m+n)^s},\\
		\zeta_{B_2}(s) &= 6^s \sum_{n,m\geq 1} \frac{1}{m^sn^s(m+n)^s (m+2n)^s}, \\
		\zeta_{G_2}(s) &= 120^s \sum_{n,m\geq 1} \frac{1}{m^sn^s(m+n)^s (m+2n)^s (m+3n)^s (2m+3n)^s},\\
		\zeta_{A_3}(s) &= 12^s \sum_{m_1,m_2,m_3\geq 1} \frac{1}{m_1^sm_2^sm_3^s (m_1+m_2)^s(m_2+m_3)^s(m_1+m_2+m_3)^s}.
	\end{aligned}$$
	It is convenient to normalize them as follows:
	\begin{alignat*}{2}
		\xi_{A_2}(s) := 2^{-s}\zeta_{A_2}(s),\qquad & \xi_{B_2}(s) := 6^{-s}\zeta_{B_2}(s),\qquad && \xi_{G_2}(s) := 120^{-s}\zeta_{G_2}(s),\\
		\xi_{A_3}(s) := 12^{-s}\zeta_{A_3}(s),\qquad & \xi_{B_3}(s):= 720^{-s}\zeta_{B_3}(s),\qquad && \xi_{C_3}(s) := 720^{-s}\zeta_{C_3}(s).
	\end{alignat*}
	
	Principal properties of these six functions are tabulated below. Many of the entries appear to be new and, to the best of our knowledge, could not be derived using previously known methods. 
	
	\begin{table}[h]
		\centering
		\begin{tabular}{|l|c|cc|}
			\hline
			$\xi_\Phi(s)$ & $A_2$ & \multicolumn{1}{c|}{$B_2$} & $G_2$ \\ \hline
			Converges when & $\Re(s)>\frac{2}{3}$ & \multicolumn{1}{c|}{$\Re(s)>\frac{1}{2}$} & $\Re(s)>\frac{1}{3}$ \\ \hline
			\begin{tabular}[c]{@{}l@{}}Residue at abscissa \\ of convergence\end{tabular} & $\dfrac{\Gamma(\frac{1}{3})^3}{2 \sqrt{3} \pi }$ & \multicolumn{1}{c|}{$\dfrac{\Gamma (\frac{1}{4})^2}{8 \sqrt{2 \pi }}$} & $\dfrac{\Gamma (\frac{1}{3})^3}{2^{\frac83} 3^{\frac32} \pi }$ \\ \hline
			$\xi_\Phi(0)$ & $\frac{1}{3}$ & \multicolumn{1}{c|}{$\frac{3}{8}$} & $\frac{5}{12}$ \\ \hline
			$\xi_\Phi(-n)$ & Vanishes for $n\in \mathbb{N}$ & \multicolumn{2}{c|}{Vanishes for $n\in 2\mathbb{N}$} \\ \hline
			\begin{tabular}[c]{@{}l@{}}Location of \\ other poles\end{tabular} & \begin{tabular}[c]{@{}c@{}}$s= \frac{k}{2}$\\ $k\leq 1, k\equiv 1 \pmod{2}$\end{tabular} & \multicolumn{1}{c|}{\begin{tabular}[c]{@{}c@{}}$s= \frac{k}{3}$\\ $k\leq 1, k\equiv 1,5 \pmod{6}$\end{tabular}} & \begin{tabular}[c]{@{}c@{}}$s = \frac{k}{5}$\\ $k\leq 1, k\equiv 1,3,7,9 \pmod{10}$\end{tabular} \\ \hline
			$\xi_\Phi'(0)$ & $\log (2\pi)$ & \multicolumn{1}{c|}{$\frac{3}{2}\log (2\pi) - \frac{1}{4}\log(2) $} & $\frac{5}{2}\log(2\pi) - \frac{1}{2}\log(2) - \frac{1}{2}\log(3)$ \\ \hline
		\end{tabular}
		\caption{\small{Major properties of rank 2 Witten zeta functions. Facts concerning $\xi_{A_2}(s)$ in the table are all known \cite{onodera2014functional, borwein2018derivatives, romik2017number, nakamura2006functional}; \cite{bridges2024asymptotic} proved first three rows for $\xi_{B_2}(s)$; Rutard in \cite{rutard2026values} establishes the rows $\xi_\Phi(0)$ and $\xi'_\Phi(0)$.}}
		\label{rank2_table}
	\end{table}  
	\begin{table}[h]
		\centering
		\begin{tabular}{|l|c|cc|}
			\hline
			$\xi_\Phi(s)$ & $A_3$ & \multicolumn{1}{c|}{$B_3$} & $C_3$ \\ \hline
			Converges when & $\Re(s)>\frac{1}{2}$ & \multicolumn{2}{c|}{$\Re(s)>\frac13$} \\ \hline
			\begin{tabular}[c]{@{}l@{}}Residue at abscissa \\ of convergence\end{tabular} & $\dfrac{\Gamma(\frac14)^4}{24\pi}$ & \multicolumn{2}{c|}{$\dfrac{\Gamma(\frac13)^6}{96\pi^2} $}\\ \hline
			$\xi_\Phi(0)$ & $-\frac{1}{4}$ & \multicolumn{2}{c|}{$-\frac{5}{16}$} \\ \hline
			$\xi_\Phi(-n)$ & Vanishes for $n\in \mathbb{N}$ & \multicolumn{2}{c|}{Vanishes for $n\in 2\mathbb{N}$} \\ \hline
			\begin{tabular}[c]{@{}c@{}}Location of poles\\ with $\Re(s)>0$\end{tabular} & $s\in \{\frac{1}{2}, \frac{2}{5}, \frac{1}{3}, \frac{1}{4}\}$ & \multicolumn{2}{c|}{$s\in \{\frac13,\frac14,\frac15,\frac16,\frac17\}$} \\ \hline
			$\xi_\Phi'(0)$ & $-\frac{3}{2}\log (2\pi)$ & \multicolumn{1}{c|}{$-\frac{45}{16}\log (2\pi) + \frac{9}{16}\log(2) $} & $-\frac{45}{16}\log(2\pi) + \frac{5}{8}\log(2)$ \\ \hline
		\end{tabular}
		\caption{\small{Major properties of rank 3 Witten zeta functions. All information beyond the first row seems new.}}
		\label{rank3_table}
	\end{table}
	\vspace*{-8mm}
	We briefly summarize our methodology and how it is different from conventional techniques in the literature. We shall analyze $\xi_{A_2}(s), \xi_{B_2}(s)$ and $\xi_{G_2}(s)$ uniformly via a more general series
	$$\xi_f(s) := \sum_{n,m\geq 1} \left(nm^{d+1} f\left(\frac{n}{m}\right)\right)^{-s},$$
	where $f(x) \in \mathbb{R}[x]$ is a degree $d$ polynomial with all roots negative. Rank two Witten zeta functions are recovered by letting $f(x)$ be \begin{equation}\label{root_system_f}f_A(x) := 1+x,\quad f_B(x):=(1+x)(1+2x),\quad f_G(x):=(1+x)(1+2x)(1+3x)(2+3x).\end{equation}
	
	Our starting point is the Mellin transform
	$$F_f(s;z) := \int_0^\infty f(x)^{-s} x^{z-1} dx.$$
	Taking the inverse Mellin transform, we have
	\begin{equation}\label{aux_15}f(x)^{-s} = \frac{1}{2\pi i} \int_{c-i\infty}^{c+i\infty} F_f(s;z) x^{-z} dz, \qquad \text{ for } x>0 \text{ and some real }c.\end{equation}
	Then, we can rewrite $\xi_f(s)$ as
	\begin{align}\label{aux_16}\xi_f(s) &= \sum_{n,m\geq 1} \frac{1}{2\pi i} \int_{c-i\infty}^{c+i\infty} n^{-s} m^{-(d+1)s} F_f(s;z) \left(\frac{n}{m}\right)^{-s} dz \nonumber \\
		&= \frac{1}{2\pi i} \int_{c-i\infty}^{c+i\infty} F_f(s;z) \zeta(s+z)\zeta((d+1)s-z) dz, 
	\end{align}
	we shall see that the exchange of summation and integral is justified. For $f(x) = 1+x$, $F_{1+x}(s;z)$ is simply the beta function and we have
	\begin{equation}\label{aux_17}(1+x)^{-s} = \frac{1}{2\pi i} \int_{c-i\infty}^{c+i\infty} \frac{\Gamma (z) \Gamma (s-z)}{\Gamma (s)} x^{-z} dz,\end{equation}
	$$\xi_{A_2}(s) = \frac{1}{2\pi i} \int_{c-i\infty}^{c+i\infty} \frac{\Gamma (z) \Gamma (s-z)}{\Gamma (s)} \zeta(s+z) \zeta(2s-z) dz.$$
	
	A well-established approach to studying Witten zeta functions\footnote{or, more generally, zeta functions with linear forms in denominator} involves repeatedly applying equation~(\ref{aux_17}), see \cite{matsumoto2003analytic, matsumoto2003asymptotic}. This method expresses $\xi_{B_2}(s)$ (respectively $\xi_{G_2}(s)$) as a 2-dimensional (respectively 4-dimensional) integral whose integrand can be written in terms of gamma and Riemann zeta functions. In general, for $\xi_\Phi(s)$ one obtains an $N$-dimensional integral, where $N$ is the number of positive roots in $\Phi$ minus the rank of $\Phi$.
	
	The advantage of this simple technique is that the integrand involves only familiar functions, namely the gamma and zeta functions. However, a drawback is the large dimension of the resulting integral, making it extremely difficult to extract finer analytic information, such as the entries presented in the tables above.
	
	Our approach departs from this traditional framework. Instead of repeatedly applying equation~(\ref{aux_17}), we treat $\xi_f(s)$ as a single integral as in equation~(\ref{aux_16}), where the integrand involves $F_f(s;z)$. The fact that $F_f(s;z)$ does not reduce in general to named special functions might be intimidating. However, this actually turns out to be advantageous: not only does it simplify the notation, but it also lets us concentrate on the essential details.
	
	In Section~2, we study the properties of $F_f(s;z)$ in detail, in preparation for our analysis of $\xi_f(s)$. Our investigation culminates in the following results:
	
	\begin{enumerate}[leftmargin=*]
				\item $\xi_f(s)$ has abscissa of convergence $s = \frac{2}{d+2}$ and has a simple pole at this point with a residue $$\frac{1}{d+2}\int_0^\infty (xf(x))^{-\frac{2}{d+2}} dx, \qquad \text{(Theorem \ref{res_theorem_rank2})}$$
		\item $\xi_f(s)$ is analytic at $s=0,-1,-2,\ldots$, with an explicit formula for its values at these points (Theorem \ref{zeta_f_negative_int}),
				\item apart from the pole at $s=\frac{2}{d+2}$, all other poles of $\xi_f(s)$ are at $s\in \frac{1}{1+d}\mathbb{Z}_{\leq 1}$ and are at most simple, with an explicit formula to calculate the residues (Theorem \ref{non_int_poles}),
		\item if roots of $f(x)$ are rational numbers, then $\xi'_f(0)$ is a $\mathbb{Q}$-linear combination of $\zeta'(-1)$ and $\log \Gamma(r)$ for various $r\in \mathbb{Q}$ (Section 6). 
	\end{enumerate}
	
	At each of the points mentioned above, the functions $\xi_{A_2}(s)$, $\xi_{B_2}(s)$ and $\xi_{G_2}(s)$ display properties that differ markedly from those of a generic $\xi_f(s)$. To the best of our knowledge, none of these phenomena have been explored in the existing literature. We briefly clarify these observations below.
	
	\subsection{Derivative at origin}
	Point (4) above allows us to calculate $\xi_f'(0)$. For illustration, one has
	$$ \xi'_f(0) = -\frac{4}{9}\zeta'(-1) - \frac{25}{108}\log(3) + \frac{4}{3}\log(2\pi) + \frac{1}{3}\log\left(\Gamma\left(\frac{1}{3}\right)\right), \qquad \text{ for }f(x)=(1+x)(1+3x).$$
	This expression for $\xi'_f(0)$ seems quite complicated in comparison to those of $\xi'_{A_2}(0), \xi'_{B_2}(0), \xi'_{G_2}(0)$ given in Table \ref{rank2_table}, namely, it contains extra terms $\zeta'(-1)$ and $\log(\Gamma(\frac{1}{3}))$. So the Witten zeta function seems to have a nice derivative at $s=0$. This leads us to the following conjecture. We first normalize the Witten zeta function $$\xi_{\Phi}(s) := K_\Phi^{-s} \zeta_\Phi(s),$$
	where $K_\Phi$ is a certain positive integer whose values are given in Section 10, Table \ref{K_Phi_value_table}. We conjecture the following shape of $\xi'_\Phi(0)$.
	
	\begin{conjecture}\label{conj_der_0}
		We have
		$$\xi_\Phi'(0) \in   \begin{cases}
			\mathbb{Q} \log(2\pi)  \qquad &\text{ for }\Phi = A_n, \\
			\mathbb{Q} \log(2\pi) + \mathbb{Q} \log (2)  \qquad &\text{ for }\Phi = B_n, C_n, D_n,\\
			\mathbb{Q} \log(2\pi) + \mathbb{Q} \log (2) + \mathbb{Q} \log (3) \qquad &\text{ for }\Phi = G_2, F_4, E_6, E_7,\\
			\mathbb{Q} \log(2\pi) + \mathbb{Q} \log (2) + \mathbb{Q} \log (3) + \mathbb{Q} \log (5) \qquad &\text{ for }\Phi = E_8.
		\end{cases}$$
	\end{conjecture}
	
	Looking at Table \ref{rank3_table}, we see the conjecture also holds for rank $3$ root systems, which we will prove in Section 10.
	
	\subsection{Poles of $\xi_f(s)$}
	We illustrate the point (3) using $\xi_{G_2}(s)$ and a generic $\xi_f(s)$ with $\deg (f) = 4$. This $\xi_f(s)$ has poles at $s=k/5$ where $k\leq 1$ is an integer and not a multiple of $5$. Using the formula of residues given in Theorem \ref{non_int_poles}, one sees that they are in general true poles\footnote{that is, the residues are non-zero}. In Example \ref{G2_skip_pole_ex}, we will be surprised to learn (and it is also non-trivial to prove) that, for $\xi_{G_2}(s)$, the points corresponding to $k\equiv 2,4,6,8\pmod{10}$ are analytic; hence only those $k\equiv 1,3,7,9 \pmod{10}$ could give poles. \par
	
	From another perspective, consider the following multi-variable generalization of $\xi_{G_2}(s)$:
	$$\xi_{G_2}(s_1,\ldots,s_6) := \sum_{m,n\geq 1}\frac{1}{m^{s_1}n^{s_2} (m+n)^{s_3} (m+2n)^{s_4} (m+3n)^{s_5} (2m+3n)^{s_6}}.$$
	In \cite[Theorem ~3.1]{komori2010witten}, it was shown that $\xi_{G_2}(s_1,\ldots,s_6)$ has meromorphic continuation to $\mathbb{C}^6$ and its possible singularities are hyperplanes given by
	\begin{equation}\begin{cases}
			s_1+s_3+s_4+s_5+s_6 = k,\quad k\in \mathbb{Z}_{\leq 1} \\
			s_2+s_3+s_4+s_5+s_6 = k,\quad k\in \mathbb{Z}_{\leq 1} \\
			s_1+s_2+s_3+s_4+s_5+s_6 = 2.
	\end{cases}\end{equation}
	Setting all $s_1,\ldots,s_6$ to be equal, we see the possible singularities of $\xi_{G_2}(s)$ are $s=1/3$ and $s=k/5, k\leq 1$. So $\xi_{G_2}(s)$ being analytic for $k\equiv 2,4,6,8\pmod{10}$ is also unexpected. This peculiar property also occurs for root systems $B_2, A_3, B_3$ and $C_3$ (the remark after Theorem \ref{A3_residues}).
	
	\subsection{Values of $\xi_f(-n), n\in \mathbb{N}$}
	The Witten zeta function displays special properties regarding the point (2) above. We introduce the notation: for a function $h(x)$ meromorphic at $x=0$, we denote the $n$-th coefficient of its Laurent series by $[h(x)][x^n]$. Then, with $g(x) = x^d f(1/x)$ another degree-$d$ polynomial, Theorem \ref{zeta_f_negative_int} gives the following formula for $\xi_f(-n)$,
	\begin{multline*}\xi_f(-n) = \frac{\zeta(-(d+2)n-1)}{d+1}\left([-f(x)^n \log f(x)][x^{1+(1+d)n}] + [-g(x)^n \log g(x)][x^{1+(1+d)n}]\right) \\ +\sum_{k=0}^{nd} [f(x)^n][x^k] \zeta(-n-k)\zeta(-(d+1)n+k).\end{multline*}
	
	A recent and deep result of the author \cite{au2024vanishing} says $\xi_\Phi(-2\mathbb{N}) = 0$. Combining this with the above formula gives three equalities, with $f_B(x), f_G(x)$ as in equation (\ref{root_system_f}),
	$$\xi_{A_2}(-2n) = 0 \iff \frac{(2n)!}{(4n+1)!}\zeta(-6n-1) = \sum_{k=0}^{2n} \frac{1}{k! (2 n-k)!} \zeta (-k-2 n) \zeta (k-4 n),$$
	$$\xi_{B_2}(-2n) = 0 \iff \frac{\zeta(-8n-1)}{3}(1+2^{-1-4n})[f_B(x)^{2n}\log f_B(x)][x^{1+6n}] = \sum_{k=0}^{4n} [f_B(x)^{2n}][x^k] \zeta (-k-2 n) \zeta (k-6 n),$$
	$$\xi_{G_2}(-2n) = 0 \iff \frac{\zeta(-12n-1)}{5}(1+3^{-1-6n})[f_G(x)^{2n}\log f_G(x)][x^{1+10n}] = \sum_{k=0}^{8n} [f_G(x)^{2n}][x^k] \zeta (-k-2 n) \zeta (k-10 n).$$
	
	Such identities for Bernoulli numbers\footnote{recall Riemann zeta at negative integers are equivalent to Bernoulli numbers}, called lacunary recurrences, are rare and highly non-trivial \cite{agoh2007convolution}. Among the three identities above, a direct proof is known only for the first one \cite{romik2017number}. More surprisingly, we conjecture that they can be lifted to the Eisenstein series by a formal replacement (Conjecture \ref{eis_conj}). \par
	
	It is well-known that the Riemann zeta function at negative odd integers carries deep $p$-adic properties \cite{koblitz2012p}. To motivate ourselves, recall the Kummer congruence for the Riemann zeta function in the $3$-adic numbers \cite[Corollary~5.14]{washington2012introduction}. Define $K(n) := 3n (1-3^{2n-1})\zeta(1-2n)$. Then, for any $n\in \mathbb{N}$, 
	\begin{itemize}[leftmargin=*]
		\item $K(n) \in \mathbb{Z}_3^\times$, i.e. a $3$-adic unit;
		\item the congruence $K(n+3^{m-1}) \equiv K(n) \pmod{3^m}$ holds.
	\end{itemize}
	
	Based on the above formula for computing $\xi_f(-n)$ and numerical experimentation, we discover the following surprising generalization of such congruences to Witten zeta function, potentially hinting at interesting underlying arithmetic phenomenon. For example, we have the $3$-adic analogue of Kummer's congruence for $\xi_{B_2}(1-2n)$. 
	\begin{conjecture}\label{p_adic_conj_intro}
				Let $n\in \mathbb{N}$ and define $$U(n) := 3^{2-2n}\xi_{B_2}(1-2n).$$ Then $U(n) \in \mathbb{Z}_3^\times$ and we have
		$$U(n+3^{m-2})\equiv U(n)\pmod{3^m}, \qquad m\geq 2, n>\lfloor \frac{m}{2}\rfloor.$$
	\end{conjecture}
	
	See Conjectures \ref{B2G2p-adic_B} and \ref{B2G2p-adic_C} for two further conjectural examples of these Kummer-type congruences. 
	
	Our formula for $\xi_f(-n)$ above is also valuable from a computational point of view. Other formulas for $\xi_f(-n)$ in the literature, such as those in~\cite{rutard2026values, bruna2025polynomials}, exhibit exponential complexity in $n$, requiring enormous calculation even for moderate $n \approx 20$. In contrast, our formula scales polynomially\footnote{assuming Bernoulli numbers can be computed in $O(1)$ time, which is a reasonable assumption for $n$ up to about $10^5$.} with $n$, and remains practical even for $n$ in the thousands. We hope our formula will facilitate further exploration of the $p$-adic aspects of Witten zeta functions.
	
	\subsection{Residue of $\xi_\Phi(s)$ at abscissa of convergence}
	We consider the point (1) above. It states that the residue of $\xi_f(s)$ at the abscissa of convergence is
	$$\frac{1}{d+2}\int_0^\infty (xf(x))^{-\frac{2}{d+2}} dx.$$
	For generic $f(x)$, this is the simplest form one could expect. It is thus surprising that the corresponding integral for $\xi_{A_2}(s), \xi_{B_2}(s), \xi_{G_2}(s)$ evaluates nicely in terms of gamma functions\footnote{for $A_2$ and $B_2$, this is less surprising, since the integrand is simple enough}. For $G_2$, the residue entry of Table \ref{rank2_table} is equivalent to the evaluation (Corollary \ref{FG_eval}) $$\int_0^\infty \frac{dx}{\sqrt[3]{x(1+x)(1+2x)(1+3x)(2+3x)}} = \frac{1}{2^{\frac53} 3^{\frac 12}\pi} \Gamma\left(\frac{1}{3}\right)^3. $$
	This is quite an unexpected evaluation, considering the fact that the algebraic curve $y^3 = x(1+x)(1+2x)(1+3x)(2+3x)$ has high genus $4$.\par
	
	It is known that the abscissa of convergence for $\zeta_\Phi(s)$ is $\frac{\text{rank}}{\#\text{positive roots}}$ \cite{hasa2019representation}. Looking at Table \ref{rank3_table}, we see that the residues for rank-three Witten zeta functions are also simple expressions involving $\Gamma(\frac 13)$ or $\Gamma(\frac14)$. Proving them requires some wizardry on hypergeometric series (Corollary \ref{FB_diag_eval}). 
	
	\begin{conjecture}\label{residue_gamma}
				For any root system $\Phi$, the residues of $\xi_\Phi(s)$ at the abscissa of convergence are always of the form $\overline{\mathbb{Q}} \times \text{(a product of gamma functions at rational numbers)}.$
	\end{conjecture}
	Numerical evidence for the conjecture exists for $\Phi = A_4$,
	$$\Res{s=2/5} \xi_{A_4}(s) \stackrel{?}{=} \frac{\sqrt{50-22 \sqrt{5}} }{240 \pi }\Gamma \left(\frac{1}{5}\right)^5.$$
	
	\vspace{3mm}

	From Section 7 onward, we shift our attention from rank two to rank three. Let $f(x_1,x_2)$ be a product of linear polynomials with non-negative coefficients, having total degree $d$. Define
	$$\xi_f(s) := \sum_{m_1,m_2,n\geq 1} \left[ n^{d+1} m_1 m_2 f\left(\frac{m_1}{n},\frac{m_2}{n}\right) \right]^{-s}.$$
	We recover $\xi_{A_3}(s), \xi_{B_3}(s), \xi_{C_3}(s)$ by specializing to certain choices of $f(x_1,x_2)$. We extend the previous methodology for rank two, starting with the Mellin transform 
	$$F_f(s;z_1,z_2) := \int_0^\infty \int_0^\infty f(x_1,x_2)^{-s} x_1^{z_1-1} x_2^{z_2-1} dx_1 dx_2,$$
	then taking its inverse transform, we have
	$$f(x_1,x_2)^{-s} = \frac{1}{(2\pi i)^2} \int_{c-i\infty}^{c+i\infty} \int_{c-i\infty}^{c+i\infty} F_f(s;z_1,z_2) x_1^{-z_1} x_2^{-z_2} dz_1 dz_2, \quad x_1,x_2\geq 0,\quad  c_i \text{ sufficiently large.}$$
	Therefore we have the following integral representation
	$$\xi_f(s) = \frac{1}{(2\pi i)^2} \int_{c-i\infty}^{c+i\infty} \int_{c-i\infty}^{c+i\infty} \zeta((d+1)s-z_1-z_2) \zeta(s+z_1)\zeta(s+z_2)F_f(s;z_1,z_2) dz_1 dz_2.$$
	When compared to the rank two case, the increase in difficulty is only computational rather than conceptual. \par
	
	We obtain in Theorems \ref{A3_residues}, \ref{B3_residues}, \ref{C3_residues} the exact values of the residues for rank-three Witten zeta functions at positive poles. For example, Theorem \ref{A3_residues} states
	$$\begin{aligned}
		\Res{s=1/2} \xi_{A_3}(s) &= \frac{\Gamma(\frac14)^4}{24\pi}, \\  \Res{s=2/5} \xi_{A_3}(s) &= \frac{\left(\sqrt{5}+5\right) \Gamma \left(\frac{1}{5}\right) \Gamma \left(\frac{3}{5}\right)}{10 \Gamma \left(\frac{4}{5}\right)} \zeta \left(\frac{2}{5} \right), \\
		\Res{s=1/3} \xi_{A_3}(s) &= \frac{2}{3}\xi_{A_2}\left(\frac{1}{3}\right), \\
		\Res{s=1/4} \xi_{A_3}(s) &= \frac{1}{4}\zeta \left(\frac14\right)^2.
	\end{aligned}$$
	and these are all poles of $\xi_{A_3}(s)$ with $\Re(s)>0$. Such nice expressions for the residues seem almost impossible to obtain by previous methods in the literature. Moreover, we will see that the poles of the $A_2$, $B_2$, $G_2$, and $A_3$ zeta functions on $\mathbb{C}$ are simple, whereas the $B_3$ and $C_3$ zeta functions have a family of double poles (Theorems \ref{B3_residues}, \ref{C3_residues}) on the negative real axis. Hence, for Witten zeta functions of rank $\leq 3$, all poles on the positive real axis are simple\footnote{However, there is no reason to believe this fact generalizes to higher-rank root systems}.
 	
	Given a root system $\Phi$, a natural counting problem is the number of $n$-dimensional representations of the Lie algebra corresponding to $\Phi$. We denote this number by $r_\Phi(n)$. When $\Phi = A_1$, this is the famous partition function. A celebrated result of Hardy and Ramanujan \cite{hardy1918asymptotic} states that
	$$r_{A_1}(n) \sim \frac{1}{4\sqrt{3}n} \exp\left( \pi \sqrt{\frac{2n}{3}}\right).$$ 
	A combination of the circle method and the saddle point method, adapted in~\cite{bridges2024asymptotic, debruyne2020saddle}, can be employed to determine the leading asymptotic behavior of $r_\mathfrak{g}(n)$, provided sufficient information is available regarding the locations of the poles and the residues of $\zeta_\mathfrak{g}(s)$. As an application of our new results on Witten zeta functions, we derive the leading asymptotic for $r_\Phi(n)$.
	\begin{theorem}
				For rank-$2$ and rank-$3$ root systems, $r_\Phi(n)$ has the following form
		$$r_\Phi(n) \sim \frac{C}{n^b} \exp\left(\sum_{j\in J} a_j n^j\right),\qquad n\to \infty,$$
				Here, $C$ and $a_j$ are explicit constants that depend on $\Phi$ (see Section 11 for their explicit values). Also\footnote{The astute readers might wonder why $0$ is included in some $J$, since it could be absorbed into the constant $C$. The reason for this will become clear later.}, 
		\begin{itemize}[leftmargin=*]
			\item for $\Phi = A_2$, $b=\frac{3}{5},\ J = \{\frac{2}{5},\frac{3}{10},\frac{1}{5},\frac{1}{10},0 \}$;
			\item for $\Phi = B_2$, $b=\frac{7}{12},\ J = \{\frac{1}{3},\frac{2}{9},\frac{1}{9},0\}$;
			\item for $\Phi = G_2$, $b=\frac{9}{16},\ J = \{\frac{1}{4},\frac{3}{20},\frac{1}{20}\}$;
			\item for $\Phi = A_3$, $b=1,\ J = \{\frac{1}{3},\frac{4}{15},\frac{2}{9},\frac{1}{5},\frac{1}{6},\frac{7}{45},\frac{2}{15},\frac{1}{9},\frac{1}{10},\frac{4}{45},\frac{1}{15},\frac{1}{18},\frac{2}{45},\frac{1}{30},\frac{1}{45},0\}$;
			\item for $\Phi = B_3$ or $C_3$, $b=\frac{71}{64},\ J = \{\frac{1}{4},\frac{3}{16},\frac{3}{20},\frac{1}{8},\frac{3}{28},\frac{7}{80},\frac{1}{16},\frac{1}{20},\frac{5}{112},\frac{1}{40},\frac{1}{140},0\}$.
		\end{itemize}
	\end{theorem}
	\par
	
	The results above for $A_2, B_2$ were already derived in \cite{romik2017number} and \cite{bridges2024asymptotic}, respectively. A partial result for $G_2$ was given in \cite{rutard2026values}. \\
	
	The structure of this article is as follows: Sections 2 to 6 focus on the rank-2 case. In Section 2, we introduce the necessary background on our Mellin transform framework for rank-2 Witten zeta functions, which is essential for establishing the meromorphic continuation of $\xi_f(s)$ in Section 3. Section 4 examines the values of $\xi_f$ at non-positive integers, and Section 5 classifies all its poles. In Section 6, we compute a certain Mellin convolution required to evaluate $\xi_f'(0)$. Sections 7 to 10 address rank-3 Witten zeta functions. Section 7 extends the framework of the Mellin transform in the previous sections to accommodate rank-3 cases; Section 8 serves as an interlude discussing certain definite integrals. In Section 9, we analyze residues with $\Re(s) > 0$, and Section 10 calculates $\xi_\Phi(0)$ and $\xi_\Phi'(0)$ for rank-three $\Phi$. Finally, in Section 11, we apply our findings to the asymptotic behavior of the number of representations of Lie algebras.  
	
	\section*{Acknowledgment}
	The author thanks Prof. Kathrin Bringmann, Dr Johann Franke and Prof. Don Zagier for valuable discussions and feedback. The author has received funding from the European Research Council (ERC) under the European Union’s Horizon 2020 research and innovation programme (grant agreement No. 101001179).

	\section{Integration kernel (rank 2)}
	In this section, we assemble materials that we will need for our investigation of $\xi_f(s)$. 
	
	\begin{lemma}\label{integral_over_gamma_entire}
		Let $f(x_1,\ldots,x_n)$ be analytic on a neighbourhood of $[0,1]^n \subset \mathbb{C}^n$ and positive on $[0,1]^n$. Then
		$$\frac{1}{\Gamma(z_1)\ldots \Gamma(z_n)}\int_{[0,1]^n} f(x_1,\ldots,x_n)^{-s} x_1^{z_1}\ldots x_n^{z_n} \frac{dx_1}{x_1} \ldots \frac{dx_n}{x_n}, \qquad s\in \mathbb{C}, \Re(z_1) > 0,\ldots \Re(z_n)>0,$$
		extends to an entire function on $(s,z_1,\ldots,z_n)\in \mathbb{C}^{n+1}$. 
	\end{lemma}
	\begin{proof}
		First we consider the case $n=1$. We first claim that \begin{equation}\label{aux_1}\int_0^1 f(x)^{-s} x^{z} \frac{dx}{x} = (e^{2\pi i z}-1)^{-1} \int_{C(1)} f(x)^{-s} x^z \frac{dx}{x}, \quad \Re (z) > 0,\end{equation}
		with the contour $C(1)$ shown: it starts at $1$, goes counterclockwise around $0$ and ends at $1$. Here we also agree $0\leq \arg x < 2\pi$. To see this, let $\varepsilon>0$ and let $C_\varepsilon$ be the small circle of radius $\varepsilon$ around the origin. Then
		$$\int_{C(1)} f(x)^{-s} x^z \frac{dx}{x} = \int_1^{\varepsilon} f(x)^{-s} x^{z-1} dx + \int_{C_\varepsilon} f(x)^{-s} x^{z-1} dx + \int_{\varepsilon}^1 f(x)^{-s} (e^{2\pi i z} x^{z-1}) dx.$$
		When $\Re(z) > 0$, the integral along $C_\varepsilon$ tends to $0$ as $\varepsilon \to 0$, so we have equality (\ref{aux_1}). 
		
		\begin{figure}[h]
			\centering
			\begin{tikzpicture}[decoration={markings,
					mark=at position 0.3cm with {\arrow[line width=1pt]{>}},
					mark=at position 5cm with {\arrow[line width=1pt]{>}},
					mark=at position 7.85cm with {\arrow[line width=1pt]{>}},
					mark=at position 13cm with {\arrow[line width=1pt]{>}}
				}
				]
				\draw[help lines,->] (-1,0) -- (8,0) coordinate (xaxis);
				\draw[help lines,->] (0,-1) -- (0,1) coordinate (yaxis);

				\path[draw,line width=0.8pt,postaction=decorate] (7.5,0) -- (7.5,0.3) -- (0.0,0.3) \centerarcpath(0,0)(90:270:0.3) -- (7.5,-0.3) -- (7.5,0);
				
				\filldraw[black] (7.5,0) circle (1pt) node[anchor=west]{$1$};
			\end{tikzpicture}\caption{The contour $C(1)$.}
		\end{figure}
		
		The contour $C(1)$ does not pass through $x=0$, thus $\int_{C(1)} f(x)^{-s} x^z \frac{dx}{x}$ defines an entire function on $(s,z)\in \mathbb{C}^2$, so the RHS of equation (\ref{aux_1}) provides an analytic continuation of the LHS. For general $n$, equation (\ref{aux_1}) can be evidently generalized to $$\int_{[0,1]^n} f(x_1,\ldots,x_n)^{-s} x_1^{z_1}\ldots x_n^{z_n} \frac{dx_1}{x_1} \ldots \frac{dx_n}{x_n} = \prod_{k=1}^n (e^{2\pi i z_k}-1)^{-1} \int_{C(1)^n} f(x_1,\ldots,x_n)^{-s} x_1^{z_1}\ldots x_n^{z_n} \frac{dx_1}{x_1} \ldots \frac{dx_n}{x_n} .$$
		We need to show $$\prod_{k=1}^n \frac{(e^{2\pi i z_k}-1)^{-1}}{\Gamma(z_k)} \int_{C(1)^n} f(x_1,\ldots,x_n)^{-s} x_1^{z_1}\ldots x_n^{z_n}  \frac{dx_1}{x_1} \ldots \frac{dx_n}{x_n} $$ is an entire function in $(s,z_1,\ldots,z_n)\in \mathbb{C}^{n+1}$. 
		As $(e^{2\pi i z}-1)^{-1}$ has simple poles at $z\in \mathbb{Z}$ and $1/\Gamma(z)$ has simple zeros at $z\in \mathbb{Z}_{\leq 0}$, it remains to show the integral on the RHS vanishes when one of $z_k$ is a positive integer. When this is the case, the $x_k$-integrand will be single-valued and analytic inside the contour $C(1)$, which is contractible to a point, so the integral is $0$.
	\end{proof}
	
	For the rest of this section, we let $f(x)$ be a degree $d$ polynomial: 
	\begin{equation}\label{cond_poly_f_1var}f(x) := c(1+\alpha_1 x)\ldots (1+\alpha_d x), \qquad \text{ such that }c, \alpha_1,\dots,\alpha_d>0.\end{equation}
	Define $$F_f(s;z) := \int_0^\infty f(x)^{-s} x^z \frac{dx}{x}.$$
	It is absolutely convergent when $d \Re(s) > \Re(z) > 0$, so $F_f(s;z)$ is an analytic function on $$\{(s,z)\in \mathbb{C}^2 : d\Re(s) > \Re(z) > 0\}.$$ We need very often another representation of $F_f(s;z)$:
	\begin{equation}\label{contourint_rep_F}F_f(s;z) = (e^{2\pi i z}-1)^{-1} \int_{C(\infty)} f(x)^{-s} x^z \frac{dx}{x}, \qquad d\Re(s)>\Re(z),\end{equation}
	here we agree $0\leq \arg x < 2\pi$ and $C(\infty)$ is the contour shown. It starts at $\infty$, goes counterclockwise around the origin and ends at $\infty$. 
	\begin{figure}[h]
		\centering
		\begin{tikzpicture}[decoration={markings,
				mark=at position 5cm with {\arrow[line width=1pt]{>}},
				mark=at position 7.85cm with {\arrow[line width=1pt]{>}},
				mark=at position 13cm with {\arrow[line width=1pt]{>}}
			}
			]
			\draw[help lines,->] (-1,0) -- (8,0) coordinate (xaxis);
			\draw[help lines,->] (0,-1) -- (0,1) coordinate (yaxis);
			
			\path[draw,line width=0.8pt,postaction=decorate] (7.5,0.3) -- (0.0,0.3) \centerarcpath(0,0)(90:270:0.3) -- (7.5,-0.3);
			
			\filldraw[black] (1,0) circle (1pt) node[anchor=west]{$1$};
		\end{tikzpicture}\caption{\small The contour $C(\infty)$.}
	\end{figure}
	
	To see this, let $\varepsilon>0$ and let $C_\varepsilon$ denote the small circle of radius $\varepsilon$ around the origin. Then
	$$\int_{C(\infty)} f(x)^{-s} x^z \frac{dx}{x} = \int_\infty^{\varepsilon} f(x)^{-s} x^{z-1} dx + \int_{C_\varepsilon} f(x)^{-s} x^{z-1} dx + \int_{\varepsilon}^\infty f(x)^{-s} (e^{2\pi i z} x^{z-1}) dx,$$
	This expression is independent of $\varepsilon$. When $\Re(z)>0$, the integral along $C_\varepsilon$ tends to $0$ as $\varepsilon \to 0$, so we have equality (\ref{contourint_rep_F}). Since the contour $C(\infty)$ does not pass through $0$, $\int_{C(\infty)} f(x)^{-s} x^z \frac{dx}{x}$ defines an analytic function on the region $\{(s,z)\in \mathbb{C}^2 : d\Re(s) > \Re(z)\}$, so equation (\ref{contourint_rep_F}) provides a meromorphic continuation of $F_f(s;z)$ to this larger region.
	
	\begin{theorem}\label{K_is_entire_function}
		The function $F_f(s;z)$ has meromorphic continuation to $\mathbb{C}^2$, and $$K_f(s;z) := \frac{\Gamma(s)}{\Gamma(z)\Gamma(ds-z)} F_f(s;z)$$ is entire.  
	\end{theorem}
	\begin{proof}
		Split the domain of integration $[0,\infty)$ into the two intervals $[0,1]$ and $[1,\infty)$. In the latter, substitute $x\mapsto 1/x$ and denote $g(x) = x^d f(1/x)$. Then
		$$F_f(s;z) = \int_0^1 f(x)^{-s} x^z \frac{dx}{x} + \int_0^1 g(x)^{-s} x^{ds-z} \frac{dx}{x}.$$
		Therefore
		$$\frac{F_f(s;z)}{\Gamma(z)\Gamma(ds-z)} = \frac{1}{\Gamma(ds-z)} \left(\frac{1}{\Gamma(z)}\int_0^1 f(x)^{-s} x^z \frac{dx}{x} \right) + \frac{1}{\Gamma(z)} \left( \frac{1}{\Gamma(ds-z)}\int_0^1 g(x)^{-s} x^{ds-z} \frac{dx}{x}\right).$$
		Both $f(x)$ and $g(x)$ are positive on $[0,1]$; hence Lemma \ref{integral_over_gamma_entire} implies that the two expressions in parentheses extend to entire functions. Since the reciprocal of the gamma function is also entire, we see that $\frac{F_f(s;z)}{\Gamma(z)\Gamma(ds-z)}$ extends to an entire function. 
		
		\par To prove that $K_f(s;z) = \Gamma(s) \left( \frac{F_f(s;z)}{\Gamma(z)\Gamma(ds-z)}\right)$ is still entire, we need to show that $F_f(-n;z) = 0$ for generic $z$ and $n\in \mathbb{Z}_{\geq 0}$. Fix $n$. By the paragraph above, $F_f(-n;z)$ is a meromorphic function of $z$ on $\mathbb{C}$; thus, by the principle of analytic continuation, it suffices to prove that it vanishes when $z$ is real and sufficiently negative. From representation (\ref{contourint_rep_F}), we need to prove \begin{equation}\label{aux_3}\int_{C(\infty)} f(x)^n x^z \frac{dx}{x} = 0 \qquad \text{ when }z\text{ is real and sufficiently negative}.\end{equation}
		Let $C(R)$ be the truncation of $C(\infty)$ at a large real number $R$, and let $C'(R)$ be the circle of radius $R$ as in the figure. \begin{figure}[h]
			\centering
			\begin{tikzpicture}[decoration={markings,
					mark=at position 0.3cm with {\arrow[line width=1pt]{>}},
					mark=at position 5cm with {\arrow[line width=1pt]{>}},
					mark=at position 7.85cm with {\arrow[line width=1pt]{>}},
					mark=at position 13cm with {\arrow[line width=1pt]{>}}
				}
				]
				\draw[help lines,->] (-4,0) -- (4,0) coordinate (xaxis);
				\draw[help lines,->] (0,-4) -- (0,4) coordinate (yaxis);
				
				\path[draw,line width=0.8pt,postaction=decorate] (3.5,0.3) -- (0.0,0.3) -- \centerarcpath(0,0)(90:270:0.3) -- (3.5,-0.3) -- \centerarcpath(0,0)(355.1:4.9:3.51283);
				
				\node at (-2.8,1) {$C'(R)$};
				\node at (1.5,0.5) {$C(R)$};
				\filldraw[black] (3.5,0) circle (1pt) node[anchor=west]{$R$};
			\end{tikzpicture}
		\end{figure}
		Because $f(x)^n x^z \frac{1}{x}$ is analytic inside the contour, $$\int_{C(R)} f(x)^n x^z \frac{dx}{x} + \int_{C'(R)} f(x)^n x^z \frac{dx}{x} = 0.$$
		For $z$ real and sufficiently negative, the second term above tends to $0$ as $R\to \infty$ because $$\left| \int_{C'(R)} f(x)^n x^z \frac{dx}{x}\right| = O(R^{dn+z}),\quad R\to \infty;$$ the first term tends to $\int_{C(\infty)} f(x)^n x^z \frac{dx}{x}$. Hence we have (\ref{aux_3}).
	\end{proof}
	
	\begin{remark}\label{remark_d_12}
		For $d=1,2$, $K_f(s;z)$ can be expressed in terms of well-known special functions:
		\begin{alignat}{2}
			K_f(s;z) &= \alpha_1^{-z}   && \text{ when }f(x)= 1+\alpha_1 x, \\
			K_f(s;z) &= \frac{\Gamma(s)}{\Gamma(2s)} \alpha_1^{-z} \: \pFq{2}{1}{s,z}{2s}{1-\frac{\alpha_2}{\alpha_1}} \qquad \qquad && \text{ when }f(x) =(1+\alpha_1 x)(1+\alpha_2 x),
		\end{alignat}
		where $_2F_1$ is the Gauss hypergeometric function. These identities follow from formulas in \cite{gradshteyn2014table}: 
		\begin{align}\label{2F1integral}
			\int_0^\infty \frac{x^{z-1}}{(1+\alpha x)^a} dx &= \alpha^{-z} \frac{\Gamma(z)\Gamma(a-z)}{\Gamma(a)},\nonumber \\
			\int_0^\infty x^{z-1} (1+x)^{-a} (1+\lambda x)^{-b} dx &= \frac{\Gamma(z)\Gamma(a+b-z)}{\Gamma(a+b)} \times \pFq{2}{1}{b,z}{a+b}{1-\lambda}
		\end{align}
		
		For higher $d$, an explicit formula is in general not possible. Nonetheless, most properties of $K_f(s;z)$ that we need can be derived without needing an explicit formula.
	\end{remark}
	
	\begin{remark}
	Theorem \ref{K_is_entire_function} and its two-variable analogue, Proposition \ref{K_is_entire_function_2}, are special cases of the general theory of Euler--Mellin integrals developed in \cite[Theorem~2.4 and Remark~2.6]{berkesch2014euler}. We nevertheless include self-contained proofs, since our later arguments require more explicit information about $K_f$ than these general results provide.
	\end{remark}
	
	If $h(x)$ is a function analytic at the origin, we use $[h(x)][x^n]$ to denote its $n$-th Taylor coefficient around $x=0$. For a polynomial $f(x)$ considered above and any $s\in \mathbb{C}$, $f(x)^{-s}$ is analytic at the origin. It is easy to see that $[f(x)^{-s}][x^n]$ is a polynomial in $s$ for any $n\geq 0$.
	
	\begin{proposition}\label{Kf(s,-n)}
		When $n$ is a non-negative integer, 
		$$K_f(s;-n) = (-1)^n n!\frac{\Gamma(s)}{\Gamma(ds+n)} [f(x)^{-s}][x^n],$$
		as an equality of meromorphic functions in $s$. 
	\end{proposition}
	\begin{proof}
		By analytic continuation, we only need to prove this when $s$ is real and sufficiently large. In this case, we can use representation (\ref{contourint_rep_F}):
		$$K_f(s;z) = \frac{1}{e^{2\pi i z}-1} \frac{\Gamma(s)}{\Gamma(z)\Gamma(ds-z)} \int_{C(\infty)} f(x)^{-s} x^z \frac{dx}{x}.$$
		Letting $z$ tend to $-n$, we have
		$$K_f(s;-n) = \frac{(-1)^n n! \Gamma(s)}{\Gamma(ds+n)} (\frac{1}{2\pi i})\int_{C(\infty)} f(x)^{-s} x^{-n} \frac{dx}{x}.$$
		The exponent of $x$ is now an integer, so the integrand is single-valued across the real axis. The only pole inside is $x=0$, and the residue theorem gives
		$$\frac{1}{2\pi i}\int_{C(\infty)} f(x)^{-s} x^{-n} \frac{dx}{x} = [f(x)^{-s}][x^n],$$
		as desired. 
	\end{proof}
	
	\begin{proposition}\label{K_f_atnon-negative-integers}
		When $m,n$ are non-negative integers,
		$$K_f(-m;-n) = \begin{cases} (-1)^{(d+1)m} d \frac{n! (dm-n)!}{m!} [f(x)^m][x^n], &\quad n\leq dm, \\ (-1)^{n+m} \frac{n!}{m! (n-dm-1)!} [-f(x)^m \log f(x)][x^n], &\quad n> dm. \end{cases}$$
	\end{proposition}
	\begin{proof}
		When $n\leq dm$, we simply let $s\to -m$ in the above proposition. Here we note that both $\Gamma(s)$ and $\Gamma(ds+n)$ have simple poles at $s=-m$; this gives the formula in the case $n\leq dm$. \par When $n>dm$, $\Gamma(s)$ has a pole at $s=-m$, $\Gamma(ds+n)$ is regular with value $(n-dm-1)!$, and $[f(x)^{-s}][x^n]$ has a zero at $s=-m$.\footnote{Because $f(x)^m$ is a polynomial of degree $dm$ and we are extracting a coefficient higher than its degree.} By differentiating with respect to $s$, we obtain its first-order term\footnote{Our hypothesis on $f(x)$ says that it maps $\mathbb{R}_{\geq 0}$ to $\mathbb{R}_{>0}$, so with the principal branch of the logarithm, $\log f(x)$ is well-defined on a neighbourhood of $\mathbb{R}_{\geq 0} \subset \mathbb{C}$.}:
		$$[f(x)^{-s}][x^n] = [-f(x)^m \log f(x)][x^n] (s+m) + O((s+m)^2), \qquad s\to -m,$$
		Therefore, we obtain the formula in the case $n>dm$. 
	\end{proof}
	
	\begin{proposition}\label{f_g_relation}
		Let $g(x) = x^d f(1/x)$, which is another degree-$d$ polynomial. Then
		$$F_f(s;z) = F_g(s;ds-z), \qquad K_f(s;z) = K_g(s;ds-z).$$
	\end{proposition}
	\begin{proof}
		The first equality follows by substituting $x\mapsto 1/x$ in the integral representation
		$$F_f(s;z) = \int_0^\infty f(x)^{-s} x^z\frac{dx}{x} = \int_0^\infty g(x)^{-s} x^{ds-z} \frac{dx}{x} = F_g(s;ds-z).$$
		The second equality follows since $\frac{\Gamma(s)}{\Gamma(z)\Gamma(ds-z)}$ is unchanged after replacing $z$ by $ds-z$.
	\end{proof}

	\begin{proposition}\label{Kf(0,z)}
		Let $f(x) = c(1+\alpha_1 x)\ldots (1+\alpha_d x)$, then for any $z\in \mathbb{C}$,
		$$K_f(0;z) = \alpha_1^{-z} + \ldots + \alpha_d^{-z}.$$
	\end{proposition}
	\begin{proof}By analytic continuation, it suffices to prove this when $\Re(z)$ is sufficiently negative, so the following representation is valid:
		\begin{equation}\label{aux_7}K_f(s;z) = \frac{1}{e^{2\pi i z}-1} \frac{\Gamma(s)}{\Gamma(z)\Gamma(ds-z)} \int_{C(\infty)} f(x)^{-s} x^z \frac{dx}{x}.\end{equation}
		Writing $h(z) := -(e^{2\pi i z}-1)^{-1} \frac{1}{\Gamma(z)\Gamma(-z)}$, we first claim 
		\begin{equation}\label{aux_4}K_f(0;z) =  h(z) \int_{C(\infty)}  x^z \log f(x) \frac{dx}{x}.\end{equation}
		Indeed, fix a generic $z$. Since $\Gamma(s)$ has a pole at $s=0$ and $K_f(s;z)$ is entire, the integral in equation (\ref{aux_7}) must vanish when $s=0$. Its first-order term is $$\int_{C(\infty)} f(x)^{-s} x^z \frac{dx}{x} = \left(- \int_{C(\infty)} x^z \log f(x) \frac{dx}{x} \right) s+ O(s^2), \qquad s\to 0.$$
		As $\Gamma(s) = s^{-1} + O(1)$, we see $\Gamma(s)\int_{C(\infty)} f(x)^{-s} x^z \frac{dx}{x}$ tends to $-\int_{C(\infty)} x^z \log f(x) \frac{dx}{x}$ as $s$ tends to $0$, which is what we claimed. \par
		
		Next, with $f(x) =  c(1+\alpha_1 x)\ldots (1+\alpha_d x)$, 
		$$\begin{aligned}K_f(0;z) &= h(z)\int_{C(\infty)} x^z \log[c(1+\alpha_1x)\ldots(1+\alpha_d x)] \frac{dx}{x}\\ &= h(z) (\log c) \underbrace{\int_{C(\infty)} x^z \frac{dx}{x}}_{=0} + \sum_{k=1}^d h(z)\int_{C(\infty)} x^z \log(1+\alpha_k x) \frac{dx}{x} \\ &= h(z) \left( \int_{C(\infty)} x^z \log(1+ x) \frac{dx}{x}\right)\times \sum_{k=1}^d \alpha_k^{-z},
		\end{aligned}$$
		where we used the substitution $x\mapsto x/\alpha_k$ in the last step. It remains to show that the factor in front of the expression above is $1$. This follows from equation (\ref{aux_4}) and the fact that when $f(x) = 1+x$, $K_f(s;z)$ equals $1$ identically (Remark \ref{remark_d_12}).
	\end{proof}
	
	\begin{proposition}\label{Kf(-n,z)}
		Let $f(x) = c(1+\alpha_1 x)\ldots (1+\alpha_d x)$, then
		$$K_f(-n;z) = \frac{(-1)^n}{n!} \sum_{j=0}^{nd} \sum_{k=1}^d \alpha_k^{-z-j} [f(x)^n][x^j] (z)_j (-z-nd)_{nd-j},$$
		where $(z)_j= z(z+1)\ldots (z+j-1)$ is the rising factorial.
	\end{proposition}
	\begin{proof}
		The proof uses essentially the same idea as above. Letting $s\to -n$ in equation (\ref{aux_7}) gives
		$$K_f(-n;z) = \frac{1}{e^{2\pi i z}-1} \frac{-(-1)^n}{n! \Gamma(z)\Gamma(ds-z)} \int_{C(\infty)} \log(f(x)) f(x)^n x^z \frac{dx}{x}.$$
		Rewrite the last integral as
		$$\sum_{j=0}^{nd} \sum_{k=1}^d [f(x)^n][x^j]  \int_{C(\infty)}\log(1+\alpha_k x) x^{j+z} \frac{dx}{x} = \sum_{j=0}^{nd} \sum_{k=1}^d \alpha_k^{-z-j} [f(x)^n][x^j]  \int_{C(\infty)}\log(1+x) x^{j+z} \frac{dx}{x}.$$
		From the proof above, we know that
		$$\int_{C(\infty)} \log(1+x) x^z \frac{dx}{x} = -\Gamma(z)\Gamma(-z) (e^{2\pi i z}-1).$$
		Substituting this back and some simplifications give the result. 
	\end{proof}
	
	The next lemma says that $F_f(s;z)$ has exponential decay for $z$ in the vertical direction, provided that we fix $s$.
	
	\begin{lemma}\label{rapid_decrease_Ff}
		Write $z=\sigma+it$, where $\sigma$ and $t$ are its real and imaginary parts, respectively. Let $K$ be a compact subset of $\mathbb{R}$. Then, for any $\varepsilon > 0$ and $s \in \mathbb{C}$, we have
		$$F_f(s;\sigma+ it) = O(e^{-(\pi - \varepsilon)|t|}), \qquad |t|>1, \sigma \in K,$$
		where the implied constant depends only on $f, K, s$ and $\varepsilon$. 
	\end{lemma}
	\begin{proof}
		Since $f\in \mathbb{R}[x]$, $\overline{F_f(s;z)} = F_f(\overline{s};\overline{z})$, we only need to handle the case $t>0$. We first prove the statement under the additional assumption $\sup K < d\Re(s) $. If this holds, equation (\ref{contourint_rep_F}) gives
		$$F_f(s;z) = \frac{1}{e^{2\pi iz}-1} \int_{C(\infty)} f(x)^{-s} x^z \frac{dx}{x}, \quad d\Re(s)>\Re(z).$$
		Fix $r>0$ such that $f(x)$ has no roots with absolute value $<r$. For a given small $\varepsilon>0$, write $\delta = \pi - \varepsilon$. Since all roots of $f(x)$ lie on the negative real axis, $f(x)^{-s}$ extends to an analytic function on the simply-connected region $x\in \mathbb{C}-(-\infty,0]$. By rotating the two arms of the contour $C(\infty)$ towards the negative real axis, we can deform it into another contour that consists of three parts: \begin{itemize}[leftmargin=*]
			\item $\gamma_1$: ray from $e^{\delta i}\infty$ to $e^{\delta i}r$,
			\item $\gamma_2$: circular arc from $re^{\delta i}$ to $re^{-\delta i}$,
			\item $\gamma_3$: ray from $e^{(2\pi-\delta)i}r$ to $e^{(2\pi-\delta) i}\infty$. 
		\end{itemize}
		
		\begin{figure}[h]
			\centering
			\scalebox{1.2}{\begin{tikzpicture}[decoration={markings,
						mark=at position 0.3cm with {\arrow[line width=1pt]{>}},
						mark=at position 5.5cm with {\arrow[line width=1pt]{>}},
						mark=at position 10cm with {\arrow[line width=1pt]{>}}
					}
					]
					\draw[help lines,->] (-6,0) -- (1.5,0) coordinate (xaxis);
					\draw[help lines,->] (0,-1.5) -- (0,1.5) coordinate (yaxis);
					
					\path[draw,line width=0.8pt,postaction=decorate] (-6,1.5) -- (-0.970143,0.242536)  \centerarcpath(0,0)(166:194:1) -- (-6,-1.5);
					
					\node at (-4,1.3) {$\gamma_1$};
					\node at (-4,-1.3) {$\gamma_3$};
					\node at (-1.3,0) {$\gamma_2$};
					\draw [dashed] (0,0) -- (-4,1);
					\draw [dashed] (0,0) -- (-4,-1);
					\filldraw[black] (-0.970143,0.242536) circle (1pt) node[anchor=south]{$re^{i\delta}$};
			\end{tikzpicture}}
		\end{figure}
		Hence $$F_f(s;z) = \frac{1}{e^{2\pi iz}-1} \left(\int_{\gamma_1} + \int_{\gamma_2}+ \int_{\gamma_3} \right)f(x)^{-s} x^z \frac{dx}{x}, \quad d\Re(s)>\Re(z).$$
		The integral along $\gamma_1$ satisfies $$\begin{aligned}\left| \int_{\gamma_1} f(x)^{-s} x^z \frac{dx}{x} \right| &= \left| - e^{\delta i(\sigma+it)}\int_r^\infty f(e^{\delta i}x)^{-s} x^{\sigma+it} \frac{dx}{x}\right| \\
			&\leq e^{-\delta t}  \int_r^\infty |f(e^{\delta i}x)^{-s}| x^\sigma \frac{dx}{x}\\
			&\leq e^{-\delta t}  O\left( \int_r^\infty (1+x)^{-d\Re(s)} x^\sigma \frac{dx}{x} \right)\\
			&= O(e^{-\delta t}),
		\end{aligned}$$
		where the implied constant depends only on $K, s , f$ and $\varepsilon$. Similarly, one shows 
		$$\left| \int_{\gamma_3} f(x)^{-s} x^z \frac{dx}{x} \right| = O(e^{-(2\pi - \delta)t}).$$
		For the integral along $\gamma_2$, its absolute value is
		\begin{align*}\left|\int_{\delta}^{2\pi-\delta} f(re^{i\theta})^{-s} (re^{i\theta})^z d\theta\right| &= O\left(\int_{\delta}^{2\pi - \delta} (re^{i\theta})^{\sigma+it} d\theta \right) \\ &= O\left(\int_{\delta}^{2\pi - \delta} e^{-t\theta} d\theta \right) = O(e^{-\delta t}) + O(e^{-(2\pi-\delta)t})
			,\end{align*}
		Combining these three estimates yields $$F_f(s;\sigma+it) = \frac{O(e^{-\delta t}) + O(e^{-(2\pi-\delta)t})}{e^{2\pi i (\sigma +it)}-1}.$$
		Remember our assumption $t>0$, we see the above is $O(e^{-\delta t})$, as desired. \par
		Next we remove the assumption $\sup K< d\Re(s)$. For any positive integer $N$, there exists a polynomial $R_N(z,x) \in \mathbb{R}[z,x]$ such that
		$$f(x)^{-s} = x^{-ds} R_N(s,\frac{1}{x}) + O(x^{-ds-N}), \qquad |x| \text{ large},$$
		this follows by looking at the expansion of $f(x)^{-s}$ around $x=\infty$. Therefore $$F_f(s;z) = \frac{1}{e^{2\pi iz}-1} \int_{C(\infty)} \left(f(x)^{-s} - x^{-ds} R_N(s,\frac{1}{x})\right) x^z\frac{dx}{x}, \quad d\Re(s)+N>\Re(z).$$ 
		For compact $K$, choose $N$ large enough that $\sup K < d\Re(s) + N$. Then we can use the above integral representation for $z\in K$; the rest proceeds exactly as above. Deforming into the contours $\gamma_1, \gamma_2, \gamma_3$, it is easy to see that the conclusion remains the same.
	\end{proof}

	\section{Meromorphic continuation of $\xi_f(s)$}
	We retain all notation from the previous section: $f(x)$ as in equation (\ref{cond_poly_f_1var}), $d = \deg(f)$, $$F_f(s;z) = \int_0^\infty f(x)^{-s} x^{z-1} dx, \qquad K_f(s;z) = \frac{\Gamma(s)}{\Gamma(z)\Gamma(ds-z)} F_f(s;z).$$
	
	Our central object of investigation is the series
	$$\xi_f(s) := \sum_{n,m\geq 1} \left[nm^{d+1} f\left(\frac{n}{m}\right)\right]^{-s}.$$
	It converges when $\Re(s)$ is sufficiently large. We will see below that its abscissa of convergence is $2/(d+2)$; that is, it converges for $\Re(s)>2/(d+2)$. 
	
	For the rest of the article, we introduce the abbreviation 
	$$\fint_{(c)} := \frac{1}{2\pi i} \int_{c-i\infty}^{c+i\infty}, \qquad \text{ for }c\in \mathbb{R}.$$
	
	By the inverse Mellin transform, we have
	$$f(x)^{-s} = \fint_{(c)} F_f(s;z) x^{-z} dz  \quad \text{ when } d\Re(s) > c > 0, x>0.$$
	Therefore
	$$\begin{aligned}\xi_f(s) &= \sum_{n,m\geq 1} \fint_{(c)} n^{-s} m^{-(d+1)s} F_f(s;z) \left(\frac{n}{m}\right)^{-z} dz \\
		&= \fint_{(c)} F_f(s;z) \zeta(s+z)\zeta((d+1)s-z) dz \\
		&= \frac{1}{\Gamma(s)} \fint_{(c)} \Gamma(z)\Gamma(sd-z) K_f(s;z) \zeta(s+z)\zeta((d+1)s-z) dz,
	\end{aligned}$$
	The exchange of summation and integration holds when $\Re(s)$ is sufficiently large. By Lemma \ref{rapid_decrease_Ff}, for fixed $s$, $F_f(s;z)$ has exponential decay in the vertical direction, while $\zeta$ has moderate growth, so the above integral is absolutely convergent.\par 
	
	We analyze the integrand using a "singularity diagram". The integrand (in variable $z$)
	$$\Gamma(z)\Gamma(ds-z) K_f(s;z) \zeta(s+z)\zeta((d+1)s-z)$$
	has poles\footnote{here we highlight the importance of Theorem \ref{K_is_entire_function}, which says $K_f(s;z)$ is an entire function.}  at $$z=-n,\quad z=d s+n,\quad z=1-s,\quad z=(d+1)s-1 \text{ with }n\in \mathbb{Z}_{\geq 0}.$$
	The diagram is constructed by plotting these lines with respect to coordinates $(\Re(s),\Re(z)) = (\Re(s),c)$. The diagram depends only on $d = \deg(f)$.
	\begin{figure}[h]
		\centering
		\begin{tikzpicture}[x=1.05cm,y=0.72cm]
			\begin{scope}
			\clip (-3.25,-4.45) rectangle (3.35,4.65);
			\foreach \n in {0,...,10}
				\draw[gray!55,thin] (-3.2,{-6.4+\n}) -- (3.2,{6.4+\n});
			\foreach \n in {1,...,4}
				\draw[gray!55,thin] (-3.2,-\n) -- (3.2,-\n);
			\draw[black,thick,dashed] (-3.2,4.2) -- (3.2,-2.2);
			\draw[black,thick,dashed] (-1.15,-4.45) -- (1.883,4.65);
			\end{scope}
			\draw[->] (-3.2,0) -- (3.3,0) node[right] {$\Re(s)$};
			\draw[->] (0,-4.4) -- (0,4.6) node[above] {$c$};
			\foreach \x in {-3,-2,-1,1,2,3}
				\draw (\x,0.08) -- (\x,-0.08) node[below=2pt] {\small $\x$};
			\foreach \y in {-4,-2,2,4}
				\draw (0.06,\y) -- (-0.06,\y) node[left=2pt] {\small $\y$};
			\fill (1,1) circle (1.6pt) node[right=2.5pt] {$A$};
			\fill (0.65,0.2) circle (1.6pt) node[left=2.5pt] {$B$};
			\fill (0.65,-0.8) circle (1.6pt) node[left=2.5pt] {$C$};
			\fill (0.65,-1.8) circle (1.6pt) node[left=2.5pt] {$D$};
		\end{tikzpicture}
		\caption{\small The singularity diagram when $d=2$. Its shape remains essentially the same for other $d$. The two dashed lines intersect at $\Re(s) = 2/(d+2)$.}
		\label{plot}
	\end{figure}
	These lines cut the plane into different connected components. The component containing point $A$ in the figure corresponds to a choice of $(\Re(s),c)$ for which our original integral representation 
	$$\xi_f(s) =\frac{1}{\Gamma(s)} \fint_{(c)} \Gamma(z)\Gamma(sd-z) K_f(s;z) \zeta(s+z)\zeta((d+1)s-z) dz =  \frac{1}{\Gamma(s)} \fint_{(c)}  H_f(s;z) dz$$
	is valid, where $$H_f(s;z) := F_f(s;z)\zeta(s+z) \zeta((d+1)s-z).$$ Lemma \ref{rapid_decrease_Ff} tells us that the integrand has exponential decay along the path of integration; thus we can shift the contour of integration from $\Re(z) = c$ to another $c' \in \mathbb{C}$, picking up residues in between:
	$$ \fint_{(c)}  H_f(s;z) dz =  \fint_{(c')}  H_f(s;z) dz \pm \left(\text{Sum of residues at poles with real part between }c\text{ and }c' \right),$$
	Here the sign is plus if $c'<c$ and minus if $c'>c$.
	
	For example, we can cross the line $\Re(z)=1-\Re(s)$ to reach point $B$, then 
	$$\begin{aligned}\xi_f(s) &= \Res{z=1-s}[H_f(s,z)] +  \fint_{(c)} H_f(s;z) dz 
		\\ &= \frac{\Gamma (1-s)\Gamma(d s+s-1)}{\Gamma(s)}\zeta ((d+2)s-1) K_f(s;1-s) +  \fint_{(c')}  H_f(s;z) dz, \end{aligned}$$
	where the point $(\Re(s),c')$ corresponding to the last integral now lies in the connected component containing $B$. As a function of $s$, this integral is analytic throughout that region and, in particular, near $s=2/(d+2)$.
	\begin{theorem}\label{res_theorem_rank2}
		The abscissa of convergence of $\xi_f(s)$ is $2/(d+2)$. The function has a simple pole at $s=2/(d+2)$, with residue 
		$$\frac{1}{d+2} F_f(\frac{2}{d+2};\frac{d}{d+2}) = \frac{1}{d+2}\int_0^\infty (xf(x))^{-2/(d+2)} dx.$$
	\end{theorem}
	\begin{proof}
		From the above discussion, the integral $\fint_{(c')}  H_f(s;z) dz$ is analytic near $s=2/(d+2)$ and the first term $$ \frac{\Gamma (1-s)\Gamma(d s+s-1)}{\Gamma(s)}\zeta ((d+2)s-1) K_f(s;1-s) = \zeta((d+2)s-1) F_f(s;1-s)$$ has a simple pole at this point, with residue $\frac{1}{d+2} F_f(\frac{2}{d+2};\frac{d}{d+2})$. Since $\xi_f(s) =\sum_{n,m\geq 1} \left[nm^{d+1} f\left(\frac{n}{m}\right)\right]^{-s}$ corresponds to a Dirichlet series with non-negative coefficients, its largest pole must be the abscissa of convergence.  
	\end{proof}
	
	For general $f(x)$, one does not expect this integral to simplify further. Quite surprisingly, for $f(x)$ arising from rank-$2$ root systems, the integral can be evaluated nicely:
	\begin{theorem}\label{rank2-residues-abiscssa}
		The residue of $\xi_{A_2}(s)$ at $s=2/3$ is $$\frac{\Gamma(\frac{1}{3})^3}{2 \sqrt{3} \pi }.$$ The residue of $\xi_{B_2}(s)$ at $s=1/2$ is $$\frac{\Gamma \left(\frac{1}{4}\right)^2}{8 \sqrt{2 \pi }}.$$
		The residue of $\xi_{G_2}(s)$ at $s=1/3$ is $$\frac{\Gamma \left(\frac{1}{3}\right)^3}{2^{8/3} 3^{3/2} \pi }.$$
	\end{theorem}
	\begin{proof}
		We skip the proofs for $\xi_{A_2}$ and $\xi_{B_2}$, not only because the evaluations of the integrals $\int_0^\infty (x f(x))^{-2/(d+2)} dx$ are trivial in these cases, but also because they are known results (\cite{romik2017number} for $\xi_{A_2}$ and \cite{bridges2024asymptotic} for $\xi_{B_2}$). For $\xi_{G_2}$, the claim is equivalent to the \emph{non-trivial} evaluation $$\int_0^\infty \frac{dx}{\sqrt[3]{x(1+x)(1+2x)(1+3x)(2+3x)}} = \frac{1}{2^{5/3} 3^{1/2}\pi} \Gamma\left(\frac{1}{3}\right)^3. $$
		In order not to interrupt our discussion, we postpone its proof to a later section (Corollary \ref{FG_eval}).
	\end{proof}
	
	\begin{remark}
	In \cite[Equation~(5.1)]{bridges2024asymptotic}, the authors use a different inverse Mellin representation of $\xi_{B_2}(s)$, namely
	$$\xi_{B_2}(s) = \frac{1}{\Gamma(s)} \fint_{(c)} \Gamma(s+z)\Gamma(-z) \xi_{A_2}(s,s-z,2s+z) dz,$$
	where $\xi_{A_2}(s_1,s_2,s_3) := \sum_{m,n\geq 1} m^{-s_1} n^{-s_2} (m+n)^{-s_3}$ is the multivariable analogue of $\xi_{A_2}$, also known as the Mordell--Tornheim zeta function. Using this formula to analyze $\xi_{B_2}(s)$ entails two nested contour-shifting arguments: one first continues the Mordell--Tornheim factor $\xi_{A_2}(s,s-z,2s+z)$ and then analyzes the outer integral. By contrast, our representation expresses $\xi_{B_2}(s)$ directly in terms of $K_f(s;z)$, so only one singularity-diagram analysis is needed. This advantage was already noted in the introduction; see the discussion following \eqref{aux_17}.
	\end{remark}
	
	Returning to the meromorphic continuation of $\xi_f(s)$, in the diagram we can continue from $B$ to $C$, then to $D$, and so on, picking up residues at $s=0,-1,-2,\ldots$ in the process and obtaining
	\begin{multline}\label{continuation_M}\xi_f(s) = \frac{\Gamma (1-s)\Gamma(d s+s-1)}{\Gamma(s)}\zeta ((d+2)s-1) K_f(s;1-s) + \sum_{k=0}^M \frac{(-1)^k}{k!} \frac{\Gamma(sd+k)}{\Gamma(s)} K_f(s;-k) \zeta(s-k) \zeta((d+1)s+k) \\ + \frac{1}{\Gamma(s)}\fint_{(-M-\varepsilon)} \Gamma(z)\Gamma(sd-z) K_f(s;z) \zeta(s+z)\zeta((d+1)s-z) dz,\end{multline}
	for any non-negative integer $M$ and small positive $\varepsilon > 0$. From the diagram, we see that as $M$ becomes sufficiently large, the latter integral is analytic at any point $s\in \mathbb{C}$. We have therefore proved the meromorphic continuation of $\xi_f(s)$. 
	
	\section{Values of $\xi_f(s)$ at non-positive integers}
	In equation (\ref{continuation_M}), let $s$ tend to $-n$, where $n$ is a non-negative integer. The last term involving the integral is $0$ because of the factor $1/\Gamma(s)$ in front. As for the sum, 
	$$\sum_{k=0}^M \frac{(-1)^k}{k!} \frac{\Gamma(sd+k)}{\Gamma(s)} K_f(s;-k) \zeta(s-k) \zeta((d+1)s+k),$$
	when $s\to -n$, the only terms that could be non-zero come from $k\in\{0,1,\ldots,nd-1,nd\} \cup \{n(d+1)+1\}$. The former arise from poles of $\Gamma(sd+k)$ and the latter from the single pole of $\zeta((d+1)s+k)$. Therefore
	\begin{multline*}
		\xi_f(-n) = \frac{(n!)^2}{d+1}\frac{(-1)^{d n+1}}{((d+1) n+1)!} \zeta(-(d+2)n-1) \left[K_f(-n;1+n) + K_f(-n;-1-(d+1)n)\right] \\ +\sum_{k=0}^{nd} \frac{(-1)^{n(d+1)} n!}{d\times k!(nd-k)!} K_f(-n;-k) \zeta(-n-k)\zeta(-(d+1)n+k).
	\end{multline*}
	Note that $K_g(-n;1+n) = K_f(-n;-1-(d+1)n)$ with $g(x) = x^d f(1/x)$ (Proposition \ref{f_g_relation}). Some simplifications using Proposition \ref{K_f_atnon-negative-integers} produce the following form.
	\begin{theorem}\label{zeta_f_negative_int}$\xi_f(s)$ is analytic at $s=0,-1,-2,\ldots$ and \begin{multline*}
			\xi_f(-n) = \frac{\zeta(-(d+2)n-1)}{d+1}\left([-f(x)^n \log f(x)][x^{1+(1+d)n}] + [-g(x)^n \log g(x)][x^{1+(1+d)n}]\right) \\ +\sum_{k=0}^{nd} [f(x)^n][x^k] \zeta(-n-k)\zeta(-(d+1)n+k),
		\end{multline*}
		where $g(x) = x^d f(1/x)$.\end{theorem}
	
	\begin{corollary}
		When $d=\deg (f)$ is odd and $n\in \mathbb{Z}_{\geq 0}$, we have
		$$\xi_f(-2n-1) = 0.$$
	\end{corollary}
	\begin{proof}
		The parity conditions on $d$ and $-2n-1$ ensure that every term in the above expression is zero because $\zeta(-2\mathbb{N}) = 0$.
	\end{proof}
	
	\begin{example}Let $f(x) = (1+x)(1+3x)$, $$\xi_f(s) = \sum_{n,m\geq 1}\frac{1}{n^s m^s (n+m)^s (n+3m)^s}.$$  $\xi_f(0),\xi_f(-1),\ldots,\xi_f(-4)$ equal respectively
		$$\frac{43}{108},\quad -\frac{23}{87480},\quad \frac{809}{72171},\quad \frac{10828231}{51963120},\quad \frac{2383543667}{157837977}.$$
	\end{example}
	
	The situation for the Witten zeta function is more interesting. 
	
	\begin{example}In Theorem \ref{zeta_f_negative_int}, let $f(x)$ be one of those listed in equation (\ref{root_system_f}). We obtain the values of the rank-$2$ Witten zeta functions at non-positive integers.
		\begin{table}[h]
			\resizebox{\textwidth}{!}{%
				\begin{tabular}{|c|c|ccccc|c|}
					\hline
					$s$& $0$           & \multicolumn{1}{c|}{$-1$}                    & \multicolumn{1}{c|}{$-2$} & \multicolumn{1}{c|}{$-3$}                                          & \multicolumn{1}{c|}{$-4$} & $-5$                                                                    & $-6$                  \\ \hline
					$\xi_{A_2}(s)$ & $\frac{1}{3}$          & \multicolumn{6}{c|}{$0$}                                                                                                                                                                                                                           \\ \hline
					$\xi_{B_2}(s)$ & $\frac{3}{8}$  & \multicolumn{1}{c|}{$\frac{-11}{4480}$}      & \multicolumn{1}{c|}{$0$}  & \multicolumn{1}{c|}{$\frac{-4581}{1576960}$}                       & \multicolumn{1}{c|}{$0$}  & $\frac{-287820799443}{256502005760}$                                    & $0$                   \\ \hline
					$\xi_{G_2}(s)$ & $\frac{5}{12}$ & \multicolumn{1}{c|}{$\frac{33205}{2612736}$} & \multicolumn{1}{c|}{$0$}  & \multicolumn{1}{c|}{$\frac{313071820474581425}{1580471680499712}$} & \multicolumn{1}{c|}{$0$}  & $\frac{34067192614177690045793856810696875}{2320107746327814532497408}$ & $0$                   \\ \hline
				\end{tabular}
			}
			\caption{\small Values of $\xi_{A_2}$, $\xi_{B_2}$, and $\xi_{G_2}$ at the first few non-positive integers. The values in the table have already been reported in \cite{rutard2026values, bruna2025polynomials}, using different methods. }
		\end{table}
		\vspace*{-8mm}
	\end{example}
	
	The vanishing of $\xi_{A_2}(-2n-1)$ from the table is explained by the corollary above. On the other hand, the vanishing of $\xi_{\Phi}(-2n)$ is a much deeper fact.
	
	\begin{theorem}\cite{au2024vanishing}
		Let $\Phi$ be a root system. The order of vanishing of $\zeta_\Phi(s)$ at negative even integers is at least the rank of $\Phi$.
	\end{theorem}
	
	Its proof relies on the root system's symmetry and employs techniques of an entirely different nature. Thus $\xi_{A_2}(s),\xi_{B_2}(s)$, and $\xi_{G_2}(s)$ not only vanish at $s\in -2\mathbb{N}$, but also have at least double zeros there. For $\xi_{A_2}(s)$, this was already shown in \cite{onodera2014functional}.
	
	\subsection{A connection to Eisenstein series}
	We draw some interesting consequences from the vanishing of $\xi_{B_2}(-2n), \xi_{G_2}(-2n)$. We write out, using Theorem \ref{zeta_f_negative_int}, an explicit expression that is equivalent to their vanishing.
	
	For $\xi_{A_2}$, it is better to directly use the fact that $K_f(s;z)\equiv 1$, arriving at:
	\begin{equation}\label{xiAvanish}\xi_{A_2}(-2n) = 0 \iff \frac{(2n)!}{(4n+1)!}\zeta(-6n-1) = \sum_{k=0}^{2n} \frac{1}{k! (2 n-k)!} \zeta (-k-2 n) \zeta (k-4 n).\end{equation}
	
	For $\xi_{B_2}$ and $\xi_{G_2}$, denote $f_B(x) = (1+x)(1+2x)$ and $f_G(x) = (1+x)(1+2x)(1+3x)(2+3x)$, 
	$$\xi_{B_2}(-2n) = 0 \iff \frac{\zeta(-8n-1)}{3}(1+2^{-1-4n})[f_B(x)^{2n}\log f_B(x)][x^{1+6n}] = \sum_{k=0}^{4n} [f_B(x)^{2n}][x^k] \zeta (-k-2 n) \zeta (k-6 n).$$
	$$\xi_{G_2}(-2n) = 0 \iff \frac{\zeta(-12n-1)}{5}(1+3^{-1-6n})[f_G(x)^{2n}\log f_G(x)][x^{1+10n}] = \sum_{k=0}^{8n} [f_G(x)^{2n}][x^k] \zeta (-k-2 n) \zeta (k-10 n).$$
	
	Converting the $\zeta(1-2\mathbb{N})$ to $\zeta(2\mathbb{N})$, we have
	$$\begin{aligned}\frac{6n+1}{2}\zeta(6n+2) \frac{(2n)!^2}{(4n+1)!} = &\sum_{k=1}^{n-1} \frac{\binom{2n}{2k+1}}{\binom{6n}{2n+2k+1}} \zeta(2k+2n+2) \zeta(4n-2k), \\
		\frac{8n+1}{6}\zeta(8 n+2)(2^{-4 n-1}+1)  \times [f_B(x)^{2n}\log f_B(x)][x^{1+6n}] = &\sum_{k=1}^{2n-1} \frac{[f_B(x)^{2n}][x^{2k+1}]}{\binom{8 n}{2 k+2 n+1}} \zeta (2 k+2 n+2) \zeta (6 n-2 k), \\
		\frac{12n+1}{10}\zeta(12n+2)(3^{-6 n-1}+1)  \times [f_G(x)^{2n}\log f_G(x)][x^{1+10n}] = &\sum_{k=1}^{4n-1} \frac{[f_G(x)^{2n}][x^{2k+1}]}{\binom{12 n}{2k+2 n+1}} \zeta (2 k+2 n+2) \zeta (10 n-2 k).\end{aligned}$$
	In terms of Bernoulli numbers $B_n$, these three identities relate $B_{6n+2}$ to $B_{2n+2}, \ldots, B_{4n}$; $B_{8n+2}$ to $B_{2n+2}, \ldots, B_{6n}$; and $B_{12n+2}$ to $B_{2n+2}, \ldots, B_{10n}$, respectively. Identities of this type (called lacunary recurrences; see \cite{agoh2007convolution}) are quite rare.
	
	\mycomment{
		\begin{remark}It can be shown that $$[f_B(x)^{2n} \log f_B(x)][x^{1+6n}] = \frac{\pi  27^{-n} \Gamma (2 n+1)}{\Gamma \left(n+\frac{7}{6}\right) \Gamma \left(n+\frac{5}{6}\right)}.$$ On the other hand, $[f_G(x)^{2n}\log f_G(x)][x^{1+10n}]$ does not admit such a nice form. \end{remark}}
	
	Recall the $\text{SL}_2(\mathbb{Z})$-Eisenstein series: for $\tau$ in the upper half-plane, $$G_{2k}(\tau) := \sum_{\substack{(m,n)\in \mathbb{Z}\times \mathbb{Z} \\(m,n)\neq (0,0)}} \frac{1}{(m+n\tau)^{2k}},\qquad k\geq 2.$$ Note that $\lim_{\tau \to i\infty} G_{2k}(\tau) = 2\zeta(2k)$. 
	
	Romik \cite{romik2017number} observed that, in the equation (\ref{xiAvanish}) above, if we replace every occurrence of $\zeta(2k)$ by $G_{2k}(\tau)/2$, we still obtain a valid equality. That is
	$$(6n+1)G_{6n+2}(\tau) \frac{(2n)!^2}{(4n+1)!} = \sum_{k=1}^{n-1} \frac{\binom{2n}{2k+1}}{\binom{6n}{2n+2k+1}} G_{2k+2n+2}(\tau) G_{4n-2k}(\tau),$$
	The original equality is recovered by letting $\tau  \to i\infty$. Romik gave a proof that applies to both (\ref{xiAvanish}) and the Eisenstein series version. Mertens and Rolen \cite{mertens2015lacunary} provide an alternative proof. \par
	
	Surprisingly, when one performs the same replacement to the identities coming from $\xi_{B_2}$ and $\xi_{G_2}$, they still seem to hold. We formulate them as conjectures.
	
	\begin{conjecture}\label{eis_conj}
		For $n\in \mathbb{N}$, the following two equalities hold:
		$$\begin{aligned}
			\frac{8n+1}{3}G_{8 n+2}(\tau)(2^{-4 n-1}+1)  \times [f_B(x)^{2n}\log f_B(x)][x^{1+6n}] \stackrel{?}{=} &\sum_{k=1}^{2n-1} \frac{[f_B(x)^{2n}][x^{2k+1}]}{\binom{8 n}{2k+2 n+1}} G_{2k+2 n+2}(\tau) G_{6 n-2k}(\tau), \\
			\frac{12n+1}{5}G_{12n+2}(\tau)(3^{-6 n-1}+1)  \times [f_G(x)^{2n}\log f_G(x)][x^{1+10n}] \stackrel{?}{=} &\sum_{k=1}^{4n-1} \frac{[f_G(x)^{2n}][x^{2k+1}]}{\binom{12 n}{2k+2 n+1}} G_{2 k+2 n+2}(\tau) G_{10n-2k}(\tau).\end{aligned}$$
	\end{conjecture}

	\subsection{$p$-adic properties of $\xi_{B_2}(1-2n)$ and $\xi_{G_2}(1-2n)$}
	We briefly look at values of $\xi_{B_2}(s), \xi_{G_2}(s)$ at negative odd integers, which can be computed efficiently by Theorem \ref{zeta_f_negative_int}. \par
	Here we report three experimentally discovered $p$-adic Kummer-type congruences for Witten zeta function. The first one was already mentioned in the introduction; see Conjecture \ref{p_adic_conj_intro}.
	
	\begingroup
	\renewcommand{\thetheorem}{\arabic{section}.\arabic{theorem}A}
	\renewcommand{\theHtheorem}{\thetheorem}
	\begin{conjecture}\label{B2G2p-adic}
		($3$-adic Kummer congruence for $B_2$) Let $$U(n) := 3^{2-2n}\xi_{B_2}(1-2n),$$ then $U(n) \in \mathbb{Z}_3^\times$ is a $3$-adic unit and we have
		$$U(n+3^{m-2})\equiv U(n)\pmod{3^m}, \qquad m\geq 2, n>\lfloor \frac{m}{2}\rfloor.$$
	\end{conjecture}
	\endgroup
	\addtocounter{theorem}{-1}
	\begingroup
	\renewcommand{\thetheorem}{\arabic{section}.\arabic{theorem}B}
	\renewcommand{\theHtheorem}{\thetheorem}
	\begin{conjecture}\label{B2G2p-adic_B}
		($2$-adic Kummer congruence for $B_2$) Let $$W(n) := \frac{2^{3 n+3}  n (3 n)!}{(n!)^2} \xi_{B_2}(1-2n),$$ then $W(n) \in \mathbb{Z}_2^\times$ is a $2$-adic unit and we have
		$$W(n+2^{m-1})\equiv W(n)\pmod{2^m}, \qquad m\geq 3, n>\lfloor \frac{m+1}{4}\rfloor.$$
	\end{conjecture}
	\endgroup
	\addtocounter{theorem}{-1}
	\begingroup
	\renewcommand{\thetheorem}{\arabic{section}.\arabic{theorem}C}
	\renewcommand{\theHtheorem}{\thetheorem}
	\begin{conjecture}\label{B2G2p-adic_C}
		($5$-adic Kummer congruence for $G_2$) Let $$V(n) := 5^{2-2n-\epsilon_n}\xi_{G_2}(1-2n), \qquad \epsilon_n :=\begin{cases}1& n\text{ is odd},\\0& n\text{ is even},\end{cases}$$ then $V(n) \in \mathbb{Z}_5^\times$ is a $5$-adic unit and we have
		$$V(n+2\cdot 5^{m-1})\equiv V(n)\pmod{5^m}, \qquad m\geq 1, n>\lfloor \frac{m}{2}\rfloor.$$
	\end{conjecture}
	\endgroup
	
	The above conjectures have been verified as long as the corresponding $\xi_\Phi(1-2n)$ has $n\leq 600$ for $\Phi = B_2$ and $n\leq 200$ for $\Phi = G_2$. \par
	In particular, the conjectures imply that $\xi_{B_2}(1-2n)$ and $\xi_{G_2}(1-2n)$ are always non-zero. It might be interesting to investigate this type of congruence at other primes $p$. 
	\par A Mathematica implementation to compute $\xi_{B_2}(1-2n), \xi_{G_2}(1-2n)$ can be found at \url{https://sites.google.com/view/kc-au/2412-17196}. Readers interested in exploring further arithmetic properties of these numbers are invited to try it.

	\section{Poles of $\xi_f(s)$}
	Apart from the pole of $\xi_f(s)$ at the abscissa of convergence $s=2/(d+2)$, it also has a family of simple poles that are easily described. 
	
	\begin{theorem}\label{non_int_poles}
		Let $s_0\neq 2/(d+2)$ be such that $\xi_f(s)$ has a pole at $s=s_0$. Then $s_0\notin \mathbb{Z}$ and $$s_0(1+d)-1 = -n$$ for some $n\in \mathbb{Z}_{\geq 0}$. The pole is simple, with residue
		\begin{multline*}\frac{(-1)^{n}}{n!} \frac{\Gamma(1-s_0)}{\Gamma(s_0)} \zeta((d+2)s_0-1)\frac{1}{1+d} \left[K_f(s_0;-n) + K_g(s_0;-n)\right] \\
			= \frac{\zeta((d+2)s_0-1)}{d+1} \left([f(x)^{-s_0}][x^{n}] + [g(x)^{-s_0}][x^{n}] \right),\end{multline*}
		where $g(x)=x^d f(1/x)$.
	\end{theorem}
	\begin{proof}We saw previously that $\xi_f(s)$ is analytic at integer $s$, so $s_0$ is not an integer. The discussion preceding equation (\ref{continuation_M}) gives \begin{multline*}\xi_f(s) = \frac{\Gamma (1-s)\Gamma((d+1)s-1)}{\Gamma(s)}\zeta ((d+2)s-1) K_f(s;1-s) + \sum_{k=0}^M \frac{(-1)^k}{k!} \frac{\Gamma(sd+k)}{\Gamma(s)} K_f(s;-k) \zeta(s-k) \zeta((d+1)s+k) \\ + \frac{1}{\Gamma(s)} \fint_{(-M-\varepsilon)} \Gamma(z)\Gamma(sd-z) K_f(s;z) \zeta(s+z)\zeta((d+1)s-z) dz.\end{multline*}
		By looking at the singularity diagram, we see that the integral becomes analytic at any point $s\in \mathbb{C}$ when $M$ gets large\footnote{More precisely, it is analytic on $M+1+\varepsilon>\Re (s)>
			\max\!\left(\frac{-M-\varepsilon}{d},
			\frac{1-M-\varepsilon}{d+1}\right)$.}; hence we can ignore the last integral. For $s\notin \mathbb{Z}$, the only possible poles are those coming from $$\zeta((d+2)s-1) \text{ or } \Gamma((d+1)s-1) \text{ or }\Gamma(sd+k).$$ In the first case, the pole is at $s=2/(d+2)$, a case we excluded. The third case cannot give a pole either: when $s\to -n_0/d$ where $n_0\in \mathbb{Z}_{\geq 0}$, we have $K_f(s;-k) \to 0$ (Proposition \ref{Kf(s,-n)}). \par
		So only the second case remains: $s\to s_0$ where $s_0(1+d)-1 = -n$ with $n\in \mathbb{Z}_{\geq 0}$. At this point, the only contribution to the residue comes from $k=-s_0(1+d)+1 = n$. A straightforward calculation shows that the residue is $$\frac{(-1)^{n}}{n!} \frac{\Gamma(1-s_0)}{\Gamma(s_0)} \zeta((d+2)s_0-1)\frac{1}{1+d} \left[K_f(s_0;1-s_0) + K_f(s_0;-n)\right].$$
		By Proposition \ref{f_g_relation} $K_f(s_0;1-s_0) = K_g(s_0;-n)$, so we have the first equality. The second equality is obtained by applying Proposition \ref{Kf(s,-n)}. 
	\end{proof}

	\begin{corollary}\label{second_pole_residue}
		Let $f(x) = c_0 +c_1 x +\ldots + c_d x^d$. Apart from $s=2/(2+d)$, the only other pole of $\xi_f(s)$ with $\Re(s)\geq 0$ is $s=1/(1+d)$. It is simple, and the residue at this point is
		$$\frac{c_0^{-1/(d+1)} + c_d^{-1/(d+1)}}{d+1}\zeta(\frac{1}{d+1}).$$
	\end{corollary}
	\begin{proof}
		This is the special case $s_0 = 1/(d+1), n = 0$ in the above theorem, for which $[f(x)^{-s_0}][x^0]=f(0)^{-s_0} = c_0^{-s_0}$ and $[g(x)^{-s_0}][x^0] = g(0) ^{-s_0} = c_d^{-s_0}$.
	\end{proof}
	
	\begin{example}
		For $\xi_{A_2}(s)$, we have $d=1$, $2s_0 -1 = -n \in \mathbb{Z}_{\leq 0}$. Therefore the poles are located at $s=\frac{1-n}{2}\in \{\frac12,-\frac12,-\frac32,-\frac52,\ldots\}$. Since $K_f \equiv K_g \equiv 1$, we have
		$$\Res{s=s_0} \xi_{A_2}(s) = \frac{(-1)^{n}}{n!} \frac{\Gamma(1-s_0)}{\Gamma(s_0)} \zeta(3s_0-1).$$
		The information in this example was already obtained by Romik \cite{romik2017number}.
	\end{example}
	
	\begin{example}
		For $\xi_{B_2}(s)$, we have $d=2$, $3s_0 - 1 = -n \in \mathbb{Z}_{\leq 0}$. In this case $g(x) = 2f(x/2)$, so $K_g(s;z) = 2^{z-s} K_f(s;z)$. Moreover, we have the evaluation $K_f(s;1-s) = \frac{\Gamma(s/2)}{2\Gamma(3s/2)}$. Therefore 
		$$\begin{aligned}\Res{s=s_0} \xi_{B_2}(s) &= \frac{(-1)^{n}}{n!} \frac{\Gamma(1-s_0)}{\Gamma(s_0)} \frac{\zeta(4s_0-1)}{3} [K_f(s_0;1-s_0) + K_g(s_0;1-s_0)] \\ &=  \frac{(-1)^{n}}{n!} \frac{\Gamma(1-s_0)}{\Gamma(s_0)} \frac{\zeta(4s_0-1)}{3} (1+2^{1-2s_0})K_f(s_0;1-s_0)  \\ &= \frac{(-1)^{n}}{n!} \frac{\Gamma(1-s_0)}{\Gamma(s_0)} \frac{\zeta(4s_0-1)}{3} (1+2^{1-2s_0}) \frac{\Gamma(s_0/2)}{2\Gamma(3s_0/2)}.\end{aligned}$$
		For an integer $k\leq 1$, the factor $\Gamma(3s_0/2)$ in the denominator shows that the residue is zero when $s_0 = k/3$, $k\equiv 2,4 \pmod{6}$, and is non-zero when $s_0 = k/3$, $k\equiv 1,5 \pmod{6}$. The information in this example was already obtained in \cite{bridges2024asymptotic}. 
	\end{example}
	
	\begin{example}\label{G2_skip_pole_ex}
		For $\xi_{G_2}(s)$, we have $d=4$, $5s_0 - 1 = -n \in \mathbb{Z}_{\leq 0}$. In this case $g(x) = 9f(x/3)$, so $K_g(s;z) = 3^{z-2s} K_f(s;z)$. We will show in Corollary \ref{FG_eval} that $$K_f(s;1-s) = \frac{\sqrt{\pi } 2^{2-11 s} 3^{3 s-1} \Gamma (s)}{\Gamma \left(\frac{5 s}{2}\right) \Gamma \left(\frac{3 s}{2}+\frac{1}{2}\right)} \times \pFq{2}{1}{\frac{3 s}{2},\frac{1}{2} (5 s-1)}{3 s}{\frac{3}{4}}.$$  
		Therefore $$\begin{aligned}\Res{s=s_0} \xi_{G_2}(s) &= \frac{(-1)^{n}}{n!} \frac{\Gamma(1-s_0)}{\Gamma(s_0)} \frac{\zeta(6s_0-1)}{5} [K_f(s_0;1-s_0) + K_g(s_0;1-s_0)] \\ &=  \frac{(-1)^{n}}{n!} \frac{\Gamma(1-s_0)}{\Gamma(s_0)} \frac{\zeta(6s_0-1)}{5} (1+3^{1-3s_0})K_f(s_0;1-s_0) \\
			&= \frac{(-1)^{n}}{n!} \frac{\sqrt{\pi } 2^{2-11 s_0} 3^{3 s_0-1}\Gamma(1-s_0)}{\Gamma \left(\frac{5 s_0}{2}\right) \Gamma \left(\frac{3 s_0}{2}+\frac{1}{2}\right)} \frac{\zeta(6s_0-1)}{5} (1+3^{1-3s_0})\; \pFq{2}{1}{\frac{3 s_0}{2},-\frac{n}{2}}{3 s_0}{\frac{3}{4}} .\end{aligned}$$
		When $n$ is odd, $5s_0/2$ is a non-positive integer. The residue vanishes because of the factor $\Gamma(5s_0/2)$ in the denominator, so $\xi_{G_2}(s)$ is analytic at these points, i.e. $s_0 = k/5$ with $k\leq 1$ and $k\equiv 2,4,6,8 \pmod{10}$. \\
		When $n$ is even, the $_2F_1$ terminates and therefore evaluates to a rational number; the residue is zero if and only if this rational number vanishes. This never seems to occur, so we formulate it as a conjecture here, although proving it might not be difficult. 
		\begin{conjecture}
			Let $k\leq 1$ be an integer and set $s_0=k/5$, where $k\equiv 1,3,7,9 \pmod{10}$. Then the rational number $$\pFq{2}{1}{\frac{3 s_0}{2},\frac{1}{2} (5 s_0-1)}{3 s_0}{\frac{3}{4}}$$ is non-zero.
			Equivalently, the points $s=k/5$, where $k\equiv 1,3,7,9 \pmod{10}$, are genuine poles of $\xi_{G_2}(s)$.
		\end{conjecture}
		
	\end{example}

	\section{A class of convolution integrals and the calculation of $\xi_f'(0)$}\label{convolution_int_section}
	In this section, we evaluate a class of integrals that are indispensable in calculating the derivative of the Witten zeta function at the origin. We shall develop a method applicable to both rank-$2$ and rank-$3$ Witten zeta functions. Throughout this section, let $\varepsilon$ be a small positive number. \par
	
	For $\alpha>0$, define $$\begin{aligned}\mathcal{A}(s,\alpha) &:= \int_0^\infty \frac{x^{s-1}}{(e^{\alpha x}-1)(e^x-1)} dx, \quad\Re(s)>2, \\
	\mathcal{A}(k,s,\alpha) &:= \frac{d^k}{d\alpha^k} \mathcal{A}(s,\alpha).\end{aligned}$$
	
	\subsection{$\mathcal{A}(k,s,\alpha)$ in terms of Hurwitz zeta functions.}
	In this subsection, we assume that $\alpha = p/q$ is a rational number with $(p,q)=1$ and set $l := pq$. We describe a procedure that expresses $\mathcal{A}(k,s,p/q)$ in terms of the Hurwitz zeta functions $\zeta(s,r) = \sum_{n\geq 0} (n+r)^{-s}$ for various rational numbers $r$. \par
	From the definition $\mathcal{A}(s,\alpha) = \int_0^\infty \frac{x^{s-1}}{(e^x-1)(e^{\alpha x}-1)} dx$, we have
	$$\mathcal{A}(k,s,\alpha) = \left. \frac{d^k}{d\alpha^k }\right\vert _{\alpha = p/q} \mathcal{A}(s,\alpha) = q^{s+k} \int_0^\infty \frac{x^{s+k-1}}{e^{qx}-1} \left( \left. \frac{d^k}{du^k }\right\vert _{u = px} \frac{1}{e^u-1} \right) dx.$$
	The expression inside the integrand $$\frac{1}{e^{qx}-1} \left( \left. \frac{d^k}{du^k }\right\vert _{u = px} \frac{1}{e^u-1} \right) $$
	is a rational function of $e^x$. Consider its partial-fraction decomposition
	$$\frac{1}{e^{qx}-1} \left( \left. \frac{d^k}{du^k }\right\vert _{u = px} \frac{1}{e^u-1} \right) = \sum_{\mu: \mu^l = 1} \sum_{i=0}^{k} \frac{c_{\mu,i}}{(e^x-\mu)^{i+1}},\qquad c_{\mu,i}\in \mathbb{Q}(\mu),$$ where we are summing over all $l$-th roots of unity. Recall the formula
	$$\int_0^\infty \frac{x^{s-1}}{e^x-a} dx = \Gamma(s) \frac{\Li_s(a)}{a}, \qquad\Re(s)>1, |a|<1,$$
	which implies 
	$$\int_0^\infty \frac{x^{s-1}}{(e^x-a)^{i+1}} dx = \Gamma(s) \frac{1}{i!}\frac{d^i}{da^i } \left( \frac{\Li_s(a)}{a} \right).$$
	Also recall the derivative of polylogarithm: $\frac{d}{da} \Li_s(a) = \Li_{s-1}(a)/a$. Therefore
	$$\frac{1}{\Gamma(s)} \int_0^\infty \frac{x^{s-1}}{(e^x-\mu)^{i+1}}\,dx = f_{\mu,i}(\Li_{s}(\mu),\Li_{s-1}(\mu),\ldots,\Li_{s-i}(\mu)),$$
	where $f_{\mu,i}(X_0,\ldots,X_i)$ is a homogeneous linear polynomial with coefficients in the field $\mathbb{Q}(\mu)$. Therefore
	$$\begin{aligned}\left. \frac{d^k}{d\alpha^k }\right\vert _{\alpha = p/q} \mathcal{A}(s,\alpha) &= q^{s+k} \int_0^\infty \frac{x^{s+k-1}}{e^{qx}-1} \left( \left. \frac{d^k}{du^k }\right\vert _{u = px} \frac{1}{e^u-1} \right) dx\\ &= q^{s+k}\sum_{\mu: \mu^l = 1} \sum_{i=0}^{k} c_{\mu,i} \int_0^\infty \frac{x^{s+k-1}}{(e^x-\mu)^{i+1}} dx \\ &= q^{s+k} \Gamma(s+k) \sum_{\mu: \mu^l = 1} \sum_{i=0}^{k} c_{\mu,i}   f_{\mu,i}(\Li_{s+k}(\mu),\Li_{s+k-1}(\mu),\ldots,\Li_{s+k-i}(\mu)),\end{aligned}.$$
	Denote the inner sum by $f_\mu(\Li_{s+k}(\mu),\Li_{s+k-1}(\mu),\ldots,\Li_{s}(\mu))$, where $f_\mu$ is a homogeneous linear polynomial with coefficients in $\mathbb{Q}(\mu)$. Then
	$$\left. \frac{d^k}{d\alpha^k }\right\vert _{\alpha = p/q} \mathcal{A}(s,\alpha) = q^{s+k} \Gamma(s+k) 
	\sum_{\mu: \mu^l = 1} f_\mu(\Li_{s+k}(\mu),\Li_{s+k-1}(\mu),\ldots,\Li_{s}(\mu)).$$
	
	Finally, we can convert $\Li_s(\mu)$ into Hurwitz zeta functions as follows:
	$$\Li_s(\mu) = \sum_{n\geq 1} \frac{\mu^n}{n^s} = \sum_{i=0}^\infty \sum_{j=1}^l \frac{\mu^j}{(j+il)^s} = l^{-s} \sum_{j=1}^l \mu^j \zeta(s,\frac{j}{l}).$$
	Therefore, we see that $\mathcal{A}(k,s,p/q)$ can be expressed in terms of Hurwitz zeta functions. 
	\begin{example}
		$$\begin{aligned}
			\mathcal{A}(0,s,1) &= \Gamma(s) (\zeta (s-1)-\zeta (s)),  \\
			\mathcal{A}(1,s,1) &= \frac{\Gamma(s+1)}{2} (\zeta (s)-\zeta (s-1)), \\
			\mathcal{A}(2,s,1) &= \frac{\Gamma(s+2)}{6} (2 \zeta (s-1)-3 \zeta (s)+\zeta (s+1)), \\
			\mathcal{A}(0,s,\frac{1}{2}) &= 2^{s-2} \Gamma(s) \left(2 \zeta (s-1)-\left(2^{1-s}+2\right) \zeta (s)\right), \\
			\mathcal{A}(1,s,\frac{1}{2}) &= 2^{s-2} \Gamma (s+1) \left(-2 \zeta (s-1)+4 \zeta (s)+\left(2^{-s}-2\right) \zeta (s+1)\right), \\
			\mathcal{A}(0,s,3) &= \frac{\Gamma(s)}{3} \left(3^s \zeta (s-1)+\left(-2\times 3^s-1\right) \zeta (s)+\zeta \left(s,\frac{1}{3}\right)\right).
		\end{aligned}$$
	\end{example}
	
	\subsection{A class of convolution integral}
	Recall our notation $(z)_k = z(z+1)\ldots (z+k-1)$ for the rising factorial. Let $n\geq 0$ be a non-negative integer. The following integrals occupy a central place in our framework:
	$$\mathcal{I}(k,n,\alpha) := \fint_{(-n-1-\varepsilon)} \Gamma(z)\zeta(z) \alpha^{-z} (z)_k \Gamma(-n-z)\zeta(-n-z) dz,\qquad \alpha>0,\quad k\in \mathbb{Z}_{\geq 0}.$$
	Our goal is to evaluate them when $\alpha$ is positive and rational.
	
	\begin{lemma}
		When $\alpha > 0$ and $1<c<\Re(s)-1$, we have
		$$\mathcal{A}(s,\alpha) = \fint_{(c)} \Gamma(z)\zeta(z) \alpha^{-z}\Gamma(s-z) \zeta(s-z)  dz.$$
	\end{lemma}
	\begin{proof}
		This is a simple consequence of Mellin convolution: if 
		$$F_i(s) := \int_0^\infty x^{s-1} f_i(x) dx,\qquad \alpha_i < \Re(s) < \beta_i,$$
		then $$\int_0^\infty x^{s-1} f_1(x) f_2(x) dx = \fint_{(c)} F_1(z) F_2(s-z) dz, \qquad \alpha_2+c < \Re(s) < \beta_2+c, \alpha_1<c<\beta_1.$$
		We specialize this to $f_1(x) = \frac{1}{e^x-1}, f_2(x) = \frac{1}{e^{\alpha x}-1}$, with corresponding $F_1(s) = \Gamma(s)\zeta(s)$ and $F_2(s) = \alpha^{-s}\Gamma(s)\zeta(s)$. Here $\alpha_i = 1,\beta_i = \infty$.
	\end{proof}
	
	For fixed $\alpha >0$, we can understand the analytic continuation of $\mathcal{A}(s,\alpha)$ using the method of singularity diagrams, just as for $\xi_f(s)$. More precisely, the singularities of the integrand $\Gamma(z)\zeta(z) \alpha^{-z}\Gamma(s-z) \zeta(s-z)$ are
	$$z= s+j,\quad z=-j,\quad j\in \left\{ -1,0,1,2,\ldots \right\}.$$
	
	\begin{figure}[h]
		\centering
		\begin{tikzpicture}[x=0.68cm,y=0.68cm]
			\begin{scope}
				\clip (-5.2,-5.45) rectangle (5.35,5.55);
				\foreach \n in {0,...,10}
					\draw[gray!55,thin] (-5.1,{-5.1+\n}) -- (5.1,{5.1+\n});
				\foreach \n in {1,...,5}
					\draw[gray!55,thin] (-5.1,-\n) -- (5.1,-\n);
				\draw[black,thick,dashed] (-4.45,-5.45) -- (5.1,4.1);
				\draw[black,thick,dashed] (-5.1,1) -- (5.1,1);
			\end{scope}
			\draw[->] (-5.15,0) -- (5.3,0) node[right] {$\Re(s)$};
			\draw[->] (0,-5.4) -- (0,5.5) node[above] {$c$};
			\foreach \x in {-4,-2,2,4}
				\draw (\x,0.08) -- (\x,-0.08) node[below=2pt] {\small $\x$};
			\foreach \y in {-4,-2,2,4}
				\draw (0.08,\y) -- (-0.08,\y) node[left=2pt] {\small $\y$};
		\end{tikzpicture}
		\caption{\small Singularity diagram of $\Gamma(z)\zeta(z) \alpha^{-z}\Gamma(s-z) \zeta(s-z)$.}
		\label{pi_factor_sing_diagram}
	\end{figure} 
	
	Starting from the unbounded region between the two dashed lines in the first quadrant, we go down vertically, picking up residues at $z=1,0,-1,\ldots,-n-1$, giving
	$$\begin{aligned}\mathcal{A}(s,\alpha)  &= \sum_{i=-n-1}^1 \Res{z=i} \left[ \Gamma(z)\zeta(z) \alpha^{-z}\Gamma(s-z) \zeta(s-z)\right] + \fint_{(-n-1-\varepsilon)} \Gamma(z)\zeta(z) \alpha^{-z}\Gamma(s-z) \zeta(s-z)dz    \\
		&= \left(\alpha^{-1}\Gamma(s-1)\zeta(s-1) + \sum_{i=0}^{n+1} \frac{(-\alpha)^i}{i!} \Gamma(s+i)\zeta(s+i) \zeta(-i) \right) + \fint_{(-n-1-\varepsilon)} \Gamma(z)\zeta(z) \alpha^{-z}\Gamma(s-z) \zeta(s-z)dz\end{aligned}.$$
	
	Applying $(-\alpha)^k \frac{d^k}{d\alpha^k}$ to both sides, we have
	\begin{multline*}(-\alpha)^k \frac{d^k}{d\alpha^k} \mathcal{A}(s,\alpha) - \left(k! \alpha^{-1}\Gamma(s-1)\zeta(s-1) + \sum_{i=0}^{n+1} \frac{(-1)^{i} (-i)_k}{i!} \alpha^{i} \Gamma(s+i)\zeta(s+i) \zeta(-i) \right) \\ = \fint_{(-n-1-\varepsilon)} \Gamma(z)\zeta(z) (z)_k \alpha^{-z}\Gamma(s-z) \zeta(s-z) dz.\end{multline*}
	Therefore,
	\begin{align*}\mathcal{I}(k,n,\alpha) &=  \fint_{(-n-1-\varepsilon)} \Gamma(z)\zeta(z) \alpha^{-z} (z)_k \Gamma(-n-z)\zeta(-n-z) dz \\
		&= \lim_{s\to -n} \left[ (-\alpha)^k \mathcal{A}(k,s,\alpha) - \left(k! \alpha^{-1}\Gamma(s-1)\zeta(s-1) + \sum_{i=0}^{n+1} \frac{(-1)^{i} (-i)_k}{i!} \alpha^{i} \Gamma(s+i)\zeta(s+i) \zeta(-i) \right) \right]. \end{align*}
	When $\alpha$ is a rational number, we have seen how to express $\mathcal{A}(k,s,\alpha)$ in terms of Hurwitz zeta functions. Evaluating this limit then gives our desired value of $\mathcal{I}(k,n,\alpha)$.
	
	\begin{example}
		Consider the case $(k,n,\alpha) = (0,0,1)$. We have $$\mathcal{A}(0,s,1) = (\zeta (s-1)-\zeta (s)) \Gamma(s).$$
		Hence, \begin{align*}\mathcal{I}(0,0,1) &= \lim_{s\to 0} \left[  \mathcal{A}(0,s,1) - \left(\Gamma(s-1)\zeta(s-1) + \sum_{i=0}^{1} \frac{(-1)^{i}}{i!}\Gamma(s+i)\zeta(s+i) \zeta(-i) \right) \right] \\
			&= \lim_{s\to 0} \left[ (\zeta (s-1)-\zeta (s)) \Gamma(s) - \left(\Gamma(s-1)\zeta(s-1) + \frac{1}{12} \zeta (s+1) \Gamma (s+1)-\frac{\zeta (s) \Gamma (s)}{2} \right) \right]
		\end{align*}
		This limit can be calculated in terms of various constants that occur in the Laurent expansions of $\Gamma(s)$ and $\zeta(s)$ when $s \in \{1,0,-1\}$. Some of these constants have more elementary forms; for example,
		$$\begin{aligned}\zeta(s) &= \frac{1}{s-1} + \gamma + O(s-1) \\
			\Gamma(s) &= 1 - \gamma (s-1) + O((s-1)^2) \\
			\zeta(s) &= -\frac{1}{2} -\frac{1}{2}\log (2\pi) s+ O(s^2),\end{aligned}$$
		Here $\gamma$ is Euler's constant.
		Some calculation then gives
		$$\mathcal{I}(0,0,1) = \fint_{(-1-\varepsilon)} \Gamma(z)\zeta(z) \Gamma(-z)\zeta(-z) dz = 2 \zeta'(-1)-\frac{\gamma }{12}-\frac{1}{12}+\frac{\log (2\pi )}{4}.$$
	\end{example}
	
	\begin{example}
		Some more examples: 
		$$\begin{aligned}
			\mathcal{I}(0,1,1) &= \fint_{(-2-\varepsilon)} \Gamma(z)\zeta(z) \Gamma(-1-z)\zeta(-1-z) dz = \frac{1}{24} \left(-36 \zeta '(-2)+12 \zeta'(-1)-1+\log (2 \pi )\right), \\
			\mathcal{I}(0,0,1/2) &= \fint_{(-1-\varepsilon)} \Gamma(z)\zeta(z) 2^z \Gamma(-z)\zeta(-z) dz = \frac{1}{24} \left(60 \zeta'(-1)-\gamma -4+ 5\log 2+6 \log (2\pi) \right),\\
			\mathcal{I}(1,1,2) &= \fint_{(-2-\varepsilon)} \Gamma(z)\zeta(z) 2^{-z} z\Gamma(-1-z)\zeta(-1-z) dz = -\zeta'(-1)+\frac{\log 2}{24}-\frac{\log (2\pi) }{12},\\
			\mathcal{I}(0,0,1/3) &= \fint_{(-1-\varepsilon)} \Gamma(z)\zeta(z) 3^z \Gamma(-z)\zeta(-z) dz =  \frac{10}{3} \zeta'(-1)-\frac{1}{4}-\frac{\gamma}{36}+\frac{\log (2\pi)}{12}+\frac{11 \log 3}{36}+\frac{1}{3} \log\left( \Gamma \left(\frac{1}{3}\right)\right).
		\end{aligned}$$
		In the last example, $\log\left(\Gamma(\frac13)\right)$ occurs because of the identity \cite[p~271]{whittaker1920course},
		$$\left. \frac{d}{ds}\right\vert _{s = 0} \zeta(s,a) = -\frac{1}{2}\log(2\pi) + \log \Gamma(a).$$
	\end{example}
	
	For calculation relevant to rank two and three Witten zeta functions, we shall only require the case $k=0,1$. 
	
	\subsection{Calculating $\xi_f'(0)$}
	
	Finally we are ready to evaluate $\xi_f'(0)$. Write the polynomial $f(x)$ as \begin{equation}\label{f(x)form}f(x) = c_0(1+\alpha_1x)\ldots (1+\alpha_dx) = c_0 + c_1 x + \ldots + c_d x^d.\end{equation} Then $K_f(0;z) = \sum_k \alpha_k^{-z}$ (Proposition \ref{Kf(0,z)}). From equation (\ref{continuation_M}), the following representation is valid in a neighbourhood of $s=0$:
	\begin{multline*}\xi_f(s) = \frac{\Gamma (1-s)\Gamma(d s+s-1)}{\Gamma(s)}\zeta ((d+2)s-1) K_f(s;1-s) + \sum_{i=0}^1 \zeta(s-i)\zeta((d+1)s+i) [f(x)^{-s}][x^i] \\ + \frac{1}{\Gamma(s)}\fint_{(-1-\varepsilon)} \Gamma(z)\Gamma(sd-z) K_f(s;z) \zeta(s+z)\zeta((d+1)s-z) dz.\end{multline*}
	
	Let $A_f$ be the sum of the coefficients of $s$ in the Taylor expansions of the first two terms. Then $$\begin{aligned}\xi_f'(0) &= A_f + \fint_{(-1-\varepsilon)} \Gamma(z)\Gamma(-z) K_f(0;z) \zeta(z)\zeta(-z) dz \\
		&= A_f + \sum_{k=1}^d \fint_{(-1-\varepsilon)} \Gamma(z)\Gamma(-z) \alpha_k^{-z} \zeta(z)\zeta(-z) dz \\ &= A_f + \sum_{k=1}^d \mathcal{I}(0,0,\alpha_k). \end{aligned}$$
	When all $\alpha_k \in \mathbb{Q}$, we know how to calculate $\mathcal{I}(0,0,\alpha_k)$ from the previous subsection. To calculate $A_f$, all that remains is to evaluate $\left. \frac{d}{ds}\right\vert _{s= 0} K_f(s;1-s)$. This is given below.
	
	\begin{proposition}\label{derivativeat0}Writing $f(x)$ as in equation (\ref{f(x)form}), we have $$K_f(s;1-s) = \frac{c_{d-1}}{c_d} + \left(\frac{c_{d-1}}{c_d} ((d-1)\gamma - \log c_d) + (d+1)\sum_{k=1}^d \frac{\log \alpha_k}{\alpha_k}\right)s + O(s^2), \qquad s\to 0,$$
		where $\gamma$ is Euler's constant. 
	\end{proposition}
	
	\begin{proof}
		Recall Proposition \ref{Kf(0,z)}. The constant term is $K_f(0;1) = \sum_k \alpha_k^{-1} = c_{d-1}/c_d$. The coefficient of $s$ is $$\left.\frac{d}{ds}\right\vert_{s=0} K_f(s;1-s) = \left.\frac{d}{ds}\right\vert_{s=0} K_f(s;1) - \left.\frac{d}{ds}\right\vert_{s=1} K_f(0;s),$$
		by the chain rule. The second term above is easy to evaluate by Proposition \ref{Kf(0,z)}: $\left.\frac{d}{ds}\right\vert_{s=1} K_f(0;s) = -\sum_k (\log \alpha_k)/\alpha_k$. For the first term, denote $g(x)=x^d f(1/x)$. Proposition \ref{f_g_relation} gives $K_f(s;1) = K_g(s;sd-1)$; thus
		$$\begin{aligned}\left.\frac{d}{ds}\right\vert_{s=0} K_f(s;1) &= \left.\frac{d}{ds}\right\vert_{s=0} K_g(s;sd-1) \\ &= \left.\frac{d}{ds}\right\vert_{s=0} K_g(s;-1) + d \left.\frac{d}{ds}\right\vert_{s=-1} K_g(0;s) \\
			&= \left.\frac{d}{ds}\right\vert_{s=0} \left( -\frac{\Gamma(s)}{\Gamma(sd+1)} [g(x)^{-s}][x^1] \right) + d \left.\frac{d}{ds}\right\vert_{s=-1} \left( \sum_k \alpha_k^s\right), \end{aligned}$$
		where we have used Propositions \ref{Kf(s,-n)} and \ref{Kf(0,z)}. The observation $[g(x)^{-s}][x^1] = -c_d^{-s-1}c_{d-1}$ then completes the proof. 
	\end{proof}
	
	We have now completed our description of how to calculate $\xi_f'(0)$, provided that all roots of $f(x)$ are rational numbers. 
	
	\begin{example}Let $f(x) = (1+x)(1+3x)$, $$\xi_f(s) = \sum_{n,m\geq 1}\frac{1}{n^s m^s (n+m)^s (n+3m)^s},$$ an explicit calculation gives
		$$\xi_f'(0) = -\frac{4}{9}\zeta'(-1) - \frac{25}{108}\log 3 + \frac{4}{3}\log(2\pi) + \frac{1}{3}\log\Gamma(\frac{1}{3}).$$
	\end{example}
	In view of the occurrence of $\zeta'(-1)$ and $\log\Gamma(1/3)$ in this generic example, the simplicity of the following derivatives might be surprising: 
	
	\begin{theorem}
		$$\begin{aligned}\xi_{A_2}'(0) &= \log (2\pi),\\
			\xi_{B_2}'(0) &=  \frac{3}{2}\log (2\pi)-\frac{1}{4}\log 2, \\
			\xi_{G_2}'(0) &= \frac{5}{2}\log (2\pi) - \frac{1}{2}\log 2 - \frac{1}{2}\log 3.
		\end{aligned}$$
	\end{theorem}
	
	The case of $\xi_{A_2}$ was also proved in \cite{bailey2018computation,borwein2018derivatives,onodera2014functional}. All three derivatives were also obtained by Rutard \cite{rutard2026values}. The simplicity of the derivative at $s=0$ continues for the rank-$3$ Witten zeta functions, as we will see below.

	\section{Integration kernel (rank 3)}
	We now embark on the odyssey of investigating rank-$3$ Witten zeta functions. In this section, we generalize the results of Section 2 to higher dimensions. Throughout, we fix a two-variable polynomial $f(x_1,x_2)$, write $d_1$ for its degree in $x_1$, $d_2$ for its degree in $x_2$, and $d_{12}$ for its total degree. We also assume that $f(x_1,x_2)$ factors into a product of linear polynomials: \begin{equation}\label{def_two_var_f} f(x_1,x_2) = c \prod_{k=1}^{d_{12}} (1 + \alpha_k x_1 + \beta_k x_2) \in \mathbb{R}[x_1,x_2]\end{equation} such that
	\begin{itemize}[leftmargin=*]
		\item $c>0, \alpha_k\geq 0, \beta_k\geq 0$,
		\item not all $\alpha_k$ are 0,
		\item not all $\beta_k$ are 0,
		\item the degrees $(d_1,d_2,d_{12})$ satisfy in addition the inequalities $3d_1 > d_{12}$ and $3d_2 > d_{12}$.\footnote{See Corollary \ref{cor:residue_abs} for the motivation of this rather strange assumption.}\end{itemize}
	
	\begin{lemma}\label{f_1234_lemma}
		For $f$ satisfying the condition in (\ref{def_two_var_f}), define the polynomials $f_1,f_2,f_3,f_4$ by
		\begin{alignat*}{2}
			&f(\frac{1}{x_1},x_2) = x_1^{-d_1} f_1(x_1,x_2),\qquad &&f(x_1,\frac{1}{x_2}) = x_2^{-d_2} f_2(x_1,x_2), \\
			&f(\frac{1}{x_1},\frac{1}{x_1x_2}) = x_1^{-d_{12}} x_2^{-d_2} f_3(x_1,x_2),\qquad &&f(\frac{1}{x_1x_2},\frac{1}{x_2}) = x_1^{-d_1} x_2^{-d_{12}} f_4(x_1,x_2).
		\end{alignat*}
		Then $f_i(x_1,x_2)$ are positive on $[0,\infty)\times [0,\infty)$. Moreover, $f_1(0,x_2)$ has degree $d_{12}-d_1$ and $f_2(x_1,0)$ has degree $d_{12}-d_2$.
	\end{lemma}
	\begin{proof}
		This is completely elementary.
	\end{proof}
	
	Consider the Mellin transform of $f(x_1,x_2)$, $$F_f(s;z_1,z_2) := \int_0^\infty \int_0^\infty f(x_1,x_2)^{-s} x_1^{z_1} x_2^{z_2} \frac{dx_1}{x_1} \frac{dx_2}{x_2}.$$ We will see, in the proof below, that $F_f(s;z_1,z_2)$ converges absolutely on
	\begin{equation}\label{2D_int_conv_range}\left\{(s,z_1,z_2)\in \mathbb{C}^3 : \Re(z_1)>0,\Re(z_2)>0,\Re(d_1s-z_1)>0,\Re(d_2s-z_2)>0,\Re(d_{12}s-z_1-z_2)>0\right\}.\end{equation}
	
	\begin{proposition}\label{K_is_entire_function_2}
		$F_f(s;z_1,z_2)$ has meromorphic continuation to $\mathbb{C}^3$ and 
		$$K_f(s;z_1,z_2) := \frac{\Gamma(s)}{\Gamma(z_1)\Gamma(z_2)\Gamma(d_1s-z_1)\Gamma(d_2s-z_2)\Gamma(d_{12}s-z_1-z_2)} F_f(s;z_1,z_2)$$
		is an entire function on $\mathbb{C}^3$. \end{proposition}
	\begin{proof}
		In the definition of $F_f(s;z_1,z_2)$, we separate the domain of integration $[0,\infty)^2$ into four parts, giving (we abbreviate $\frac{dx_1 dx_2}{x_1x_2} = d\mu$)
		\begin{multline*}F_f(s,z_1,z_2) = \int_{[0,1]^2} f(x_1,x_2)^{-s} x_1^{z_1} x_2^{z_2} d\mu +  \int_{[0,1]^2} f(\frac{1}{x_1},x_2)^{-s} x_1^{-z_1} x_2^{z_2} d\mu \\ + \int_{[0,1]^2} f(x_1,\frac{1}{x_2})^{-s} x_1^{z_1} x_2^{-z_2} d\mu + \int_{[0,1]^2} f(\frac{1}{x_1},\frac{1}{x_2})^{-s} x_1^{-z_1} x_2^{-z_2} d\mu.\end{multline*}
		The last integral can be further separated as:
		$$\int_{[0,1]^2} f(\frac{1}{x_1},\frac{1}{x_1x_2})^{-s} x_1^{-z_1-z_2} x_2^{-z_2} d\mu + \int_{[0,1]^2} f(\frac{1}{x_1 x_2},\frac{1}{x_2})^{-s} x_1^{-z_1} x_2^{-z_1-z_2} d\mu.$$
		Using polynomials $f_1,f_2,f_3,f_4$ defined in the lemma above, we rewrite $F_f(s;z_1,z_2)$ as a sum of five terms \begin{multline*}F_f(s;z_1,z_2) = \int_{[0,1]^2} f(x_1,x_2)^{-s} x_1^{z_1} x_2^{z_2} d\mu +  \int_{[0,1]^2} f_1(x_1,x_2)^{-s} x_1^{d_1s-z_1} x_2^{z_2} d\mu + \\ \int_{[0,1]^2} f_2(x_1,x_2)^{-s} x_1^{z_1} x_2^{d_2s-z_2} d\mu
			+ \int_{[0,1]^2}f_3(x_1,x_2)^{-s} x_1^{d_{12}s-z_1-z_2} x_2^{d_2s-z_2} d\mu  \\ +\int_{[0,1]^2}f_4(x_1,x_2)^{-s} x_1^{d_1s-z_1} x_2^{d_{12}s-z_1-z_2} d\mu.
		\end{multline*}
		
		The polynomials $f,f_1,f_2,f_3,f_4$ are all positive on $[0,1]^2$. The integral defining $F_f(s;z_1,z_2)$ is analytic when all the above integrals converge absolutely, yielding the region in \eqref{2D_int_conv_range}. According to Lemma \ref{integral_over_gamma_entire}, the first integral divided by $\Gamma(z_1)\Gamma(z_2)$ extends to an entire function, the second integral divided by $\Gamma(d_1s-z_1)\Gamma(z_2)$ extends to an entire function, and the third integral divided by $\Gamma(d_2s-z_2)\Gamma(z_1)$ extends to an entire function. Likewise, the fourth integral divided by $\Gamma(d_{12}s-z_1-z_2)\Gamma(d_2s-z_2)$ and the fifth integral divided by $\Gamma(d_{12}s-z_1-z_2)\Gamma(d_1s-z_1)$ extend to entire functions. Thus, $$\frac{F_f(s;z_1,z_2)}{\Gamma(z_1)\Gamma(z_2)\Gamma(d_1s-z_1)\Gamma(d_2s-z_2)\Gamma(d_{12}s-z_1-z_2)}$$ extends to an entire function on $\mathbb{C}^3$. \par
		Next, we need to show that the above expression multiplied by $\Gamma(s)$ is still entire, i.e. $F_f(-n;z_1,z_2) = 0$ for a non-negative integer $n$. This part is largely parallel to the proof of Theorem \ref{K_is_entire_function}. Recall the contour $C(\infty)$ used in its proof. We have, \textit{mutatis mutandis}, \begin{equation}\label{eq_F_contourint}F_f(s;z_1,z_2) = (e^{2\pi i z_1}-1)^{-1} (e^{2\pi i z_2}-1)^{-1} \int_{C(\infty)^2} f(x_1,x_2)^{-s} x_1^{z_1-1} x_2^{z_2-1}\,dx_1\,dx_2,\end{equation}
		as long as the RHS is absolutely convergent. To show $F_f(-n;z_1,z_2) = 0$, we can assume that $\Re(z_1)$ and $\Re(z_2)$ are sufficiently negative. Then
		$$F_f(-n;z_1,z_2) = (e^{2\pi i z_1}-1)^{-1} (e^{2\pi i z_2}-1)^{-1} \int_{C(\infty)^2} f(x_1,x_2)^n x_1^{z_1-1} x_2^{z_2-1}\,dx_1\,dx_2.$$
		That the expression above is zero follows by deforming the contours to $C(R)$ and $C'(R)$ as in the proof of Theorem \ref{K_is_entire_function}.
	\end{proof}
	
	\begin{proposition}\label{rank_3_non-negative-int_red2}
		Let $n_1, n_2$ be non-negative integers. Then
		$$K_f(s;-n_1,-n_2) = \frac{(-1)^{n_1+n_2} (n_1)! (n_2)! \; \Gamma(s)}{\Gamma(d_1s+n_1)\Gamma(d_2s+n_2)\Gamma(d_{12}s+n_1+n_2)} [f(x_1,x_2)^{-s}][x_1^{n_1} x_2^{n_2}],$$
		where the last term means the Taylor coefficient of $x_1^{n_1} x_2^{n_2}$ in the expansion around $(x_1,x_2)=(0,0)$.
	\end{proposition}
	\begin{proof}
		By analytic continuation, it suffices to prove this for $\Re(s)$ sufficiently large. From equation (\ref{eq_F_contourint}), we can use the representation
		$$K_f(s;z_1,z_2) = \frac{\Gamma(s)(e^{2\pi i z_1}-1)^{-1} (e^{2\pi i z_2}-1)^{-1}}{\Gamma(z_1)\Gamma(z_2)\Gamma(d_1s-z_1)\Gamma(d_2s-z_2)\Gamma(d_{12}s-z_1-z_2)} \int_{C(\infty)^2} f(x_1,x_2)^{-s} x_1^{z_1-1} x_2^{z_2-1}\,dx_1\,dx_2.$$
		Letting $z_i \to -n_i$ gives
		\begin{multline*}K_f(s;-n_1,-n_2) = \frac{(-1)^{n_1+n_2}(n_1)! (n_2)!\Gamma(s)}{\Gamma(d_1s+n_1)\Gamma(d_2s+n_2)\Gamma(d_{12}s+n_1+n_2)}\\ \times (\frac{1}{2\pi i})^2 \int_{C(\infty)^2} f(x_1,x_2)^{-s} x_1^{-n_1-1} x_2^{-n_2-1}\,dx_1\,dx_2.\end{multline*}
		The last integral is single-valued across the positive real axes, so its value is $[f(x_1,x_2)^{-s}][x_1^{n_1} x_2^{n_2}]$ by the residue theorem.
	\end{proof}
	
	For our given $f(x_1,x_2)$ satisfying the condition in equation \ref{def_two_var_f}, we can fix $x_1$ and consider it as a polynomial in $x_2$, obtaining the one-variable kernel studied before:
	$$K^{x_2}_f(s;z) := \frac{\Gamma(s)}{\Gamma(z)\Gamma(d_2s-z)} \int_0^\infty f(x_1,x_2)^{-s} x_2^{z-1} dx_2.$$
	Similarly, we can fix $x_2$ and treat $f(x_1,x_2)$ as a polynomial in $x_1$:
	$$K^{x_1}_f(s;z) := \frac{\Gamma(s)}{\Gamma(z)\Gamma(d_1s-z)} \int_0^\infty f(x_1,x_2)^{-s} x_1^{z-1} dx_1.$$
	
	From Theorem \ref{K_is_entire_function}, we know $K_{f(x_1,x_2)}^{x_2}(s;z)$ is entire in $s$ and $z$. It is also easy to see it is analytic\footnote{Because when $x_1$ is sufficiently near $0$, $f(x_1,x_2)$ as a polynomial in $x_2$, still satisfies equation (\ref{cond_poly_f_1var}), this is a consequence of condition (\ref{def_two_var_f}).} around $x_1 = 0$. Therefore for a non-negative integer $n$, it makes sense to write $[K_f^{x_2}(s;z)][x_1^n]$, which is its $x_1^n$-th coefficient.

	\begin{proposition}\label{rank_3_non-negative-int_red}
		Let $n$ be a non-negative integer. Then
		$$K_f(s;-n,z) = \frac{(-1)^n n!}{\Gamma(d_1s+n)\Gamma(d_{12}s+n-z)} [K_f^{x_2}(s;z)][x_1^n],$$
		$$K_f(s;z,-n) = \frac{(-1)^n n!}{\Gamma(d_2s+n)\Gamma(d_{12}s+n-z)} [K_f^{x_1}(s;z)][x_2^n].$$
	\end{proposition}
	\begin{proof}
		Let $\gamma(s;z_1,z_2) = \frac{\Gamma(s)}{\Gamma(z_1)\Gamma(z_2)\Gamma(d_1s-z_1)\Gamma(d_2s-z_2)\Gamma(d_{12}s-z_1-z_2)}$. By analytic continuation, it suffices to consider $\Re(s)$ sufficiently negative, so we have the representation $$K_f(s;-n+\varepsilon,z) = \gamma(s;-n+\varepsilon,z) (e^{2\pi i \varepsilon}-1)^{-1} \int_{x_1\in C(\infty)} \int_{x_2 > 0} f(x_1,x_2)^{-s} x_1^{-n+\varepsilon} x_2^z \frac{dx_1}{x_1}\frac{dx_2}{x_2}.$$
		Both $\gamma(s;-n+\varepsilon,z) (e^{2\pi i \varepsilon}-1)^{-1}$ and the integral have limit as $\varepsilon\to 0$. Calculating the limit, we see
		$$\begin{aligned}K_f(s;-n,z) &= \frac{(-1)^n n! \Gamma(s)}{\Gamma(z)\Gamma(d_1s+n) \Gamma(d_2s-z) \Gamma(d_{12}s+n-z)} \left(\int_0^\infty [f(x_1,x_2)^{-s}][x_1^n] x_2^z \frac{dx_2}{x_2} \right) \\
			&= \frac{(-1)^n n! \Gamma(s)}{\Gamma(z)\Gamma(d_1s+n) \Gamma(d_2s-z) \Gamma(d_{12}s+n-z)} \left[\int_0^\infty f(x_1,x_2)^{-s}x_2^z \frac{dx_2}{x_2}\right][x_1^n]\\
			&=  \frac{(-1)^n n!}{\Gamma(d_1s+n) \Gamma(d_{12}s+n-z)} \left[K^{x_2}_f(s;z)\right][x_1^n].\end{aligned}$$
		This proves the first equation of the proposition. The second equation is proved analogously. 
	\end{proof}
	
	Consider the polynomials $f_1,f_2,f_3,f_4$ associated with $f$ as in Lemma \ref{f_1234_lemma}. We define $$F_{f_i}(s;z_1,z_2) := \int_0^\infty \int_0^\infty f_i(x_1,x_2)^{-s} x_1^{z_1-1} x_2^{z_2-1}\,dx_1\,dx_2.$$
	Then 
	\begin{multline*}F_f(s;z_1,z_2) = F_{f_1}(s;d_1s-z_1,z_2) = F_{f_2}(s;z_1,d_2s-z_2) \\ = F_{f_3}(s;d_{12}s-z_1-z_2,d_2s-z_2)  = F_{f_4}(s;d_1s-z_1,d_{12}s-z_1-z_2).\end{multline*}
	These identities follow by making simple substitutions in the integral. For the next proposition, recall our notation for the rising factorial $(z)_n = z(z+1)\ldots (z+n-1)$. 
	
	\begin{proposition}\label{Kf_rank3_secondreduction}
		Let $n$ be a non-negative integer. \\ (a) Let $u_i(s)$ be polynomials defined by
		$$[f_1(x_1,x_2)^{-s}][x_1^n] = f_1(0,x_2)^{-s} \left(\sum_{i=0}^n u_i(s) x_2^i\right).$$
		Then $$K_f(s;d_1s+n,z) = \frac{(-1)^n n!}{\Gamma(d_1s+n)\Gamma(d_2s-z)} \sum_{i=0}^n u_i(s) (z)_i ((d_{12}-d_1)s-z-n)_{n-i} K_{f_1(0,x_2)}(s;z+i).$$
		(b) Let $v_i(s)$ be polynomials defined by
		$$[f_2(x_1,x_2)^{-s}][x_2^n] = f_2(x_1,0)^{-s} \left(\sum_{i=0}^n v_i(s) x_1^i\right).$$
		Then $$K_f(s;z,d_2s+n) = \frac{(-1)^n n!}{\Gamma(d_2s+n)\Gamma(d_1s-z)} \sum_{i=0}^n v_i(s) (z)_i ((d_{12}-d_2)s-z-n)_{n-i} K_{f_2(x_1,0)}(s;z+i).$$
	\end{proposition}
	\begin{proof}
		We prove only (a); part (b) follows by a similar argument. First, we show that $[f_1(x_1,x_2)^{-s}][x_1^n]$ indeed has the desired form. Write $f(x_1,x_2) = c \prod_{k=1}^{d_{12}} (1+\alpha_k x_1 + \beta_k x_2)$, 
		then $$\begin{aligned} f_1(x_1,x_2) = x_1^{d_1} f(\frac{1}{x_1},x_2) &= c \prod_{k: \alpha_k\neq 0} (x_1+\alpha_k +\beta_k x_1x_2) \prod_{k: \alpha_k = 0} (1+\beta_k x_2)
			\\ &= c \left(\prod_{k:\alpha_k \neq 0} \alpha_k \right) \left(\prod_{k: \alpha_k = 0} (1+\beta_k x_2)\right) \prod_{k: \alpha_k\neq 0} (1+x_1\frac{1+\beta_k x_2}{\alpha_k}) \\
			&= f_1(0,x_2) \prod_{k: \alpha_k\neq 0} (1+x_1\frac{1+\beta_k x_2}{\alpha_k}).
		\end{aligned}$$
		Raise both sides to the power $-s$ and then use the series expansion of $(1+x_1)^{-s}$. We see that $$[f_1(x_1,x_2)^{-s}][x_1^n] = f_1(0,x_2)^{-s} \left(\sum_{i=0}^n u_i(s) x_2^i\right)$$
		for some polynomials $u_i(s)$. Let $\gamma(s;z_1,z_2)$ be as in the previous proof. Then
		$$K_f(s;\varepsilon+d_1s+n,z) = \gamma(s,\varepsilon+d_1s+n,z) F_f(s;\varepsilon+d_1s+n,z) = \gamma(s;\varepsilon+d_1s+n,z) F_{f_1}(s;-\varepsilon-n,z).$$
		A similar argument as in previous proof shows \begin{multline*}K_f(s;d_1s+n,z) = \lim_{\varepsilon\to 0}\left[\gamma(s,\varepsilon+d_1s+n,z) (e^{-2\pi i \varepsilon}-1)^{-1} (2\pi i)\right] \int_0^\infty [f_1(x_1,x_2)^{-s}][x_1^n] x_2^z\frac{dx_2}{x_2} \\
			=\frac{(-1)^n n! \Gamma (s)}{\Gamma (z) \Gamma \left(n+d_1s\right) \Gamma \left(d_2s-z\right) \Gamma \left(-n-z-s d_1+s d_{12}\right)} \int_0^\infty [f_1(x_1,x_2)^{-s}][x_1^n] x_2^z\frac{dx_2}{x_2}
			.\end{multline*}
		We cannot take the $[x_1^n]$ out of the integral since the previous footnote does not hold\footnote{For example, if $f(x_1,x_2) = (1+x_1)(1+x_2)(1+x_1+x_2)$, then $f_1(x_1,x_2) = (1+x_1)(1+x_2)(1+x_1+x_1x_2)$.} for $f_1$. However, we can still expand $[f_1(x_1,x_2)^{-s}][x_1^n]$ explicitly, giving
		$$\begin{aligned}\int_0^\infty [f_1(x_1,x_2)^{-s}][x_1^n] x_2^z\frac{dx_2}{x_2} &= \sum_{i=0}^n u_i(s) \int_0^\infty f_1(0,x_2)^{-s} x_2^{z+i} \frac{dx_2}{x_2} \\
			&= \sum_{i=0}^n u_i(s) F_{f_1(0,x_2)}(s;z+i) \\
			&= \sum_{i=0}^n u_i(s) \frac{\Gamma(s+i)\Gamma((d_{12}-d_1)s-z-i)}{\Gamma(s)} K_{f_1(0,x_2)}(s;z+i).\end{aligned}$$
		In the last line, we use the fact that $f_1(0,x_2)$ has degree $d_{12}-d_1$ in $x_2$ (Lemma \ref{f_1234_lemma}). Substituting this back to the previous equation completes the proof.
	\end{proof}
	
	The next proposition is an extension of Proposition \ref{Kf(0,z)}:
	\begin{proposition}\label{Kf(0,z)rank3}
		Let $f(x_1,x_2) = c\prod_{k=1}^{d_{12}} (1+\alpha_k x_1 + \beta_k x_2)$ with $f$ satisfying condition \ref{def_two_var_f}. Then $$K_f(0;z_1,z_2) = \frac{1}{\Gamma(-z_1)\Gamma(-z_2)} \left(\sum_{\alpha_k \beta_k\neq 0} \alpha_k^{-z_1} \beta_k^{-z_2} \right),$$
		Here the sum runs over $k$ such that both $\alpha_k$ and $\beta_k$ are non-zero.
	\end{proposition}
	\begin{proof}
		We assume that $\Re(z_1)$ and $\Re(z_2)$ are sufficiently negative. Thus, we have the representation
		\begin{multline}\label{aux_6}K_f(s;z_1,z_2) = \frac{\Gamma(s) \times (e^{2\pi i z_1}-1)^{-1} (e^{2\pi i z_2}-1)^{-1}}{\Gamma(z_1)\Gamma(z_2)\Gamma(d_1s-z_1)\Gamma(d_2s-z_2)\Gamma(d_{12}s-z_1-z_2)}  \int_{C(\infty)^2} f(x_1,x_2)^{-s} x_1^{z_1-1} x_2^{z_2-1}\,dx_1\,dx_2.\end{multline}
		Letting $s\to 0$ gives
		$$K_f(0;z_1,z_2) = \frac{-(e^{2\pi i z_1}-1)^{-1} (e^{2\pi i z_2}-1)^{-1} }{\Gamma(z_1)\Gamma(z_2)\Gamma(-z_1)\Gamma(-z_2)\Gamma(-z_1-z_2)} \times \int_{C(\infty)^2} \log f(x_1,x_2) x_1^{z_1-1} x_2^{z_2-1}\,dx_1\,dx_2.$$
		Write $\log f(x_1,x_2) = \log c + \sum_{k=1}^{d_{12}} \log(1+\alpha_k x_1 + \beta_k x_2)$. Since $\int_{C(\infty)} x^{w-1}\,dx = 0$, only terms for which both $\alpha_k$ and $\beta_k$ are non-zero contribute to the integral. The substitutions $x_1\mapsto x_1/\alpha_k$ and $x_2\mapsto x_2/\beta_k$ give
		$$K_f(0;z_1,z_2) = \frac{-(e^{2\pi i z_1}-1)^{-1} (e^{2\pi i z_2}-1)^{-1} }{\Gamma(z_1)\Gamma(z_2)\Gamma(-z_1)\Gamma(-z_2)\Gamma(-z_1-z_2)} \left( \int_{C(\infty)^2} \log(1+x_1+x_2) x_1^{z_1} x_2^{z_2} \frac{dx_1}{x_1}\frac{dx_2}{x_2} \right)\left(\sum_{\alpha_k \beta_k\neq 0} \alpha_k^{-z_1} \beta_k^{-z_2} \right).$$
		We denote the expression above by \begin{equation}\label{aux_5}K_f(0;z_1,z_2) = h(z_1,z_2) \left(\sum_{\alpha_k \beta_k\neq 0} \alpha_k^{-z_1} \beta_k^{-z_2} \right).\end{equation} It remains to find $h(z_1,z_2)$, which is independent of $f$. To do this, specialize to $f= 1+x_1+x_2$. Then 
		$$F_{1+x_1+x_2}(s;z_1,z_2) = \int_0^\infty \int_0^\infty (1+x_1+x_2)^{-s} x_1^{z_1-1} x_2^{z_2-1}\,dx_1\,dx_2 = \frac{\Gamma (z_1) \Gamma (z_2) \Gamma (s-z_1-z_2)}{\Gamma (s)},$$
		here the integral is evaluated by two applications of beta integrals. Therefore $$K_{1+x_1+x_2}(s;z_1,z_2) = \frac{1}{\Gamma (s-z_1) \Gamma (s-z_2)}.$$
		Comparing this to equation (\ref{aux_5}) shows $h(z_1,z_2) = (\Gamma(-z_1)\Gamma(-z_2))^{-1}$. 
	\end{proof}
	
	The following proposition generalizes the above result when the first argument is a non-negative integer.
	\begin{proposition}\label{Kf(-n,z)rank3}
		Let $n$ be a non-negative integer, then
		\begin{multline*}K_f(-n;z_1,z_2) = \frac{(-1)^n}{n! \Gamma(-nd_1-z_1)\Gamma(-nd_2-z_2)} \\ \times \sum_{0\leq i_1+i_2\leq nd_{12}} \left(\sum_{\alpha_k \beta_k \neq 0} \alpha_k^{-z_1-i_1} \beta_k^{-z_2-i_2} \right)[f(x_1,x_2)^n][x_1^{i_1} x_2^{i_2}](z_1)_{i_1} (z_2)_{i_2} (-z_1-z_2-n d_{12})_{nd_{12}-i_1-i_2}\end{multline*}
		where the outer sum ranges over non-negative integers $i_1,i_2$ such that $i_1+i_2\leq nd_{12}$. In particular, $K_f(-n;m,z) = K_f(-n;z,m) = 0$ for non-negative integers $n,m$ and $z\in \mathbb{C}$.
	\end{proposition}
	\begin{proof}
		The assertion on vanishing follows immediately from the formula since the gamma factor is $0$ and the summation is a polynomial in $z_1,z_2$. \par
		The proof largely follows that of the proposition above. Letting $s\to -n$ in equation (\ref{aux_6}) gives
		\begin{multline*}K_f(-n;z_1,z_2) = \frac{-(-1)^n(e^{2\pi i z_1}-1)^{-1} (e^{2\pi i z_2}-1)^{-1} }{n!\Gamma(z_1)\Gamma(z_2)\Gamma(-nd_1-z_1)\Gamma(-nd_2-z_2)\Gamma(-nd_{12}-z_1-z_2)}  \\ \times \int_{C(\infty)^2} \log f(x_1,x_2) f(x_1,x_2)^n x_1^{z_1-1} x_2^{z_2-1}\,dx_1\,dx_2,\end{multline*}
		Since $\int_{C(\infty)} x^{s-1} dx = 0$, only terms for which both $\alpha_k$ and $\beta_k$ are non-zero contribute to the integral. 
		$$\begin{aligned}
			&\quad \int_{C(\infty)^2} \log f(x_1,x_2) f(x_1,x_2)^n x_1^{z_1-1} x_2^{z_2-1}\,dx_1\,dx_2   \\ &= \sum_{0\leq i_1+i_2\leq nd_{12}} \sum_{\alpha_k \beta_k \neq 0} [f(x_1,x_2)^n][x_1^{i_1} x_2^{i_2}] \int_{C(\infty)^2} \log(1+\alpha_k x_1 + \beta_k x_2) x_1^{i_1+z_1-1} x_2^{i_2+z_2-1}\,dx_1\,dx_2 \\
			&=\sum_{0\leq i_1+i_2\leq nd_{12}} \sum_{\alpha_k \beta_k \neq 0} [f(x_1,x_2)^n][x_1^{i_1} x_2^{i_2}] \alpha_k^{-z_1-i_1} \beta_k^{-z_2-i_2}  \int_{C(\infty)^2} \log(1+ x_1 + x_2) x_1^{i_1+z_1-1} x_2^{i_2+z_2-1}\,dx_1\,dx_2
		\end{aligned}$$
		From the proof of the above proposition, we know $$\int_{C(\infty)^2} \log(1+x_1+x_2) x_1^{z_1-1} x_2^{z_2-1}\,dx_1\,dx_2 = -(e^{2\pi i z_1}-1)(e^{2\pi i z_2}-1) \Gamma(z_1)\Gamma(z_2)\Gamma(-z_1-z_2).$$
		Substituting this back into the previous expression proves the claim.
	\end{proof}
	
	For our application to Witten zeta functions, we mainly focus on three examples:
	\begin{align}\label{f_ABC_rank3}
		\begin{split}
			f_A(x_1,x_2)&:= (1+x_1)(1+x_2)(1+x_1+x_2),\\
			f_B(x_1,x_2)&:= (1+x_1)(1+x_2)(1+x_1+x_2)(2+x_2)(2x_1+2+x_2)(x_1+2+x_2),\\
			f_C(x_1,x_2)&:= (1+x_1)(1+x_2)(1+x_1+x_2)(1+2x_2)(x_1+1+2x_2)(x_1+2+2x_2)
		\end{split}
	\end{align}
	The degrees $(d_1,d_2,d_{12})$ for $f_A$ are $(2,2,3)$; for $f_B$ and $f_C$, they are $(4,5,6)$. We denote the corresponding $F_f$ by $F_A, F_B, F_C$, respectively, and similarly for $K_f$. \par 
	Note that $f_B(x,2y) = 4 f_C(x,y)$ and $F_B(s;z_1,z_2) = 2^{-2s+z_2} F_C(s;z_1,z_2)$. Statements about $F_B$ that we give below can be easily translated to $F_C$, so we focus only on $F_A$ and $F_B$. Also note that $f_A$ is symmetric, so $F_A(s;z_1,z_2) = F_A(s;z_2,z_1)$ and $K_A(s;z_1,z_2) = K_A(s;z_2,z_1)$.
	
	As applications of the above results, we assemble some facts about $K_A, K_B$ that will be used later. 
	\begin{example}\label{derivative_zero_example}
		Let $z\in \mathbb{C}$. Then
		$$\left. \frac{d}{ds}\right\vert _{s=0} K_A(s;1-s,z) = -\frac{5+z}{\Gamma(-z)},$$ and $$\left. \frac{d}{ds}\right\vert _{s=0} K_B(s;1-s,z) = -\frac{2 \left(2 z+9\times 2^z+7\right)}{\Gamma (-z)},\qquad \left. \frac{d}{ds}\right\vert _{s=0} K_B(s;z,1-s) = -\frac{5 \left(z+2^{z+3}+4\right)}{2 \Gamma (-z)}.$$
	\end{example}
	\begin{proof}
		We shall only prove the last two evaluations for $K_B$, the other evaluation for $K_A$ is proved similarly. By chain rule, we have $$\left. \frac{d}{ds}\right\vert _{s=0} K_B(s;1-s,z) = \left. \frac{d}{ds}\right\vert _{s=0} K_B(s;4s+1,z) - 5 \left. \frac{d}{ds}\right\vert _{s=1} K_B(0;s,z)$$
		We calculate these two derivatives. For the former derivative, by Proposition \ref{Kf_rank3_secondreduction}, $$K_B(s;4s+1,z) = \frac{5s}{2\Gamma(4s+1)\Gamma(5s-z)} \left[(4s-2-2z)K_g(s;z) + z K_g(s;1+z)\right],\quad \text{ where }g(x) = 2(1+x)(2+x).$$
		Thus
		$$\left. \frac{d}{ds}\right\vert _{s=0} K_B(s;4s+1,z) = \frac{5}{2\Gamma(-z)} \left[(-2-2z)K_g(0;z) + z K_g(0;1+z)\right].$$ 
		Proposition \ref{Kf(0,z)} says $K_g(0;z) = 1+2^z$, thus $\left. \frac{d}{ds}\right\vert _{s=0} K_B(s;4s+1,z) = -\frac{5 \left(z+2^{z+1}+2\right)}{2 \Gamma (-z)}$. For the latter derivative $\left. \frac{d}{ds}\right\vert _{s=1} K_B(0;s,z)$, one simply uses Proposition \ref{Kf(0,z)rank3} which says $K_B(0;s,z) = \frac{2^{s+z}+2^z+1}{\Gamma (-s) \Gamma (-z)}$. Hence we have the claimed evaluation of $\left. \frac{d}{ds}\right\vert _{s=0} K_B(s;1-s,z)$. For the second claim about $K_B$, start with
		$$\left. \frac{d}{ds}\right\vert _{s=0} K_B(s;z,1-s) = \left. \frac{d}{ds}\right\vert _{s=0} K_B(s;z,5s+1) - 6 \left. \frac{d}{ds}\right\vert _{s=1} K_B(0;z,s),$$
		These two derivatives can be calculated similarly.
	\end{proof}
	
	The following proposition allows us to express $K_A(s;z_1,z_2)$ in terms of gamma functions when one of $z_i$ is a non-positive integer. This property does not generalize to $K_f(s;z_1,z_2)$ for general $f$.
	\begin{proposition}\label{KA3_non-negative_int}
		Let $n\geq 0$ be a non-negative integer,
		$$K_A(s;-n,z) = K_A(s;z,-n) = \frac{(-1)^n n! \Gamma(s)}{\Gamma(2s)\Gamma(2s+n)\Gamma(3s+n-z)} \sum_{k=0}^n \binom{-2s+z}{n-k} \frac{(s)_k (z)_k}{(1)_k (2s)_k} (-1)^k.$$
	\end{proposition}
	\begin{proof}
		The first equality is true because the last two arguments of $K_A$ are symmetric. From Proposition \ref{rank_3_non-negative-int_red}, $$K_A(s;-n,z) = \frac{\Gamma(s) (-1)^n n!}{\Gamma(2s+n)\Gamma(2s-z) \Gamma(z)\Gamma(3s+n-z)} \left( \int_0^\infty f_A(x_1,x_2)^{-s} x_2^{z-1} dx_2\right) [x_1^n],$$
		the integral evaluates to $$\frac{\Gamma(z)\Gamma(2s-z)}{\Gamma(2s)} (1+x_1)^{-2s+z} \pFq{2}{1}{s,z}{2s}{-x_1},$$ using series expansion of $_2F_1$ at origin proves the claim. 
	\end{proof}
	
	An extension of Lemma \ref{rapid_decrease_Ff} also holds for $F_f(s;z_1,z_2)$: it has exponential decay when one of $z_i$ has a large imaginary part. We omit the statement and proof here.
	
	\section{Creative telescoping and some definite integrals}
	We insert here a self-contained section with radically different flavor. Our goal is to evaluate $$\begin{aligned}F_G(s;1-s) &:= \int_0^\infty \left(x(1+x)(1+2x)(1+3x)(2+3x) \right)^{-s} dx, \\
		F_A(s;1-s,1-s)&= \int_0^\infty\int_0^\infty \left(xy(1+x)(1+y)(1+x+y)\right)^{-s} dxdy, \\
		F_B(s;1-s,1-s)&= \int_0^\infty \int_0^\infty \left(x (x+1) y (y+1) (y+2) (x+y+1) (2 x+y+2) (x+y+2) \right)^{-s} dxdy.
	\end{aligned}$$ 
	Surprisingly, these functions can all be evaluated in terms of gamma or hypergeometric functions. Having such forms is useful in finding the residues of the corresponding Witten zeta functions. \par
	
	\begin{proposition}\label{F_Adiag_eval}
		$$F_A(s;1-s,1-s) = \frac{2^{2 s-2} \Gamma (1-s)^2 \Gamma \left(\frac{3 s}{2}-\frac{1}{2}\right) \Gamma \left(\frac{5 s}{2}-1\right)}{\Gamma \left(\frac{s}{2}+\frac{1}{2}\right) \Gamma \left(\frac{3 s}{2}\right)}.$$
	\end{proposition}
	\begin{proof}
		It is possible to prove a more general formula: $$F_A(s;z_1,z_2) = \frac{\Gamma(z_1)\Gamma(z_2)\Gamma(2s-z_1)\Gamma(2s-z_2)}{\Gamma(2s)^2} \pFq{3}{2}{z_1,z_2,s}{2s,2s}{1},$$
		Indeed, 
		$$\begin{aligned}
			F_A(s;z_1,z_2) &= \int_0^\infty \int_0^\infty x^{z_1-1} y^{z_2-1} ((1+x)(1+y)(1+x+y))^{-s} dxdy \\
			&= \frac{\Gamma(z_2)\Gamma(2s-z_2)}{\Gamma(2s)} \int_0^\infty x^{z_1-1} (1+x)^{-2s} \pFq{2}{1}{z_2,s}{2s}{\frac{x}{1+x}} dx \\
			&= \frac{\Gamma(z_2)\Gamma(2s-z_2)}{\Gamma(2s)} \sum_{n\geq 0} \frac{(z_2)_n (s)_n}{(1)_n (2s)_n} \int_0^\infty x^{z_1+n-1} (1+x)^{-2s-n} dx \\
			&= \frac{\Gamma(z_2)\Gamma(2s-z_2)}{\Gamma(2s)} \frac{\Gamma(z_1)\Gamma(2s-z_1)}{\Gamma(2s)} \pFq{3}{2}{z_1,z_2,s}{2s,2s}{1}.
		\end{aligned}$$
		Invoking Dixon's $_3F_2$ summation formula completes the proof:
		\begin{equation*}\pFq{3}{2}{a,b,c}{1+a-b,1+a-c}{1} = \frac{\Gamma\left(1 + \frac{a}{2}\right) \Gamma\left(1 + a - b\right) \Gamma\left(1 + a - c\right) \Gamma\left(1 + \frac{a}{2} - b - c\right)}{\Gamma\left(1 + a\right) \Gamma\left(1 + \frac{a}{2} - b\right) \Gamma\left(1 + \frac{a}{2} - c\right) \Gamma\left(1 + a - b - c\right)}. \qedhere\end{equation*}
	\end{proof}
	The evaluations of $F_G(s;1-s)$ and $F_B(s;1-s,1-s)$ are more challenging. We need the following non-trivial definite integral:
	\begin{lemma}Let $1/5<s<1$ and $y>-1$. Then
		$$\int_0^\infty (x(1+x)(1+x+y)(2x+2+y)(x+y+2))^{-s} dx = 2^{-1-s} \frac{\Gamma(1-s)\Gamma(\frac{5s-1}{2})}{\Gamma(\frac{3s+1}{2})} (1+y)^{1-5s} \pFq{2}{1}{\frac{3s}{2},\frac{5s-1}{2}}{3s}{\frac{y(2+y)}{(1+y)^2}}.$$
	\end{lemma}
	\begin{proof}
		We give a verification proof by a technique known as \emph{creative telescoping}. Fix $s$ between $1/5$ and $1$, let $H(x,y) := (x(1+x)(1+x+y)(2x+2+y)(x+y+2))^{-s}$, and set $$g(x,y) := \frac{s x (x+1) \left(6 x^2 y+6 x^2+10 x y^2+26 x y+18 x+4 y^3+17 y^2+24 y+12\right)}{(x+y+1) (x+y+2) (2 x+y+2)}.$$
		One easily verifies\footnote{Note that after dividing by $H(x,y)$, both sides become rational functions in $x$ and $y$, so the verification reduces to the familiar operations on rational functions.}
		$$\left(y(1+y)(2+y) \frac{\partial^2}{\partial y^2} + 2s(1+2y)(3+2y) \frac{\partial}{\partial y} + 3s(5s-1)(1+y) \right) H(x,y) = \frac{\partial}{\partial x} (g(x,y) H(x,y)).$$
		Set $I(y) := \int_0^\infty H(x,y) dx$. Then $I(y)$ satisfies
		$$y(1+y)(2+y) I''(y) + 2s(1+2y)(3+2y) I'(y)+ 3s(5s-1)(1+y) I(y)= \int_0^\infty\frac{\partial}{\partial x} (g(x,y) H(x,y))  dx.$$
		By the fundamental theorem of calculus, the RHS equals $$\lim_{x\to \infty} g(x,y) H(x,y) - \lim_{x\to 0} g(x,y) H(x,y) = 0.$$ Thus $I(y)$ solves the second-order ODE $$y(1+y)(2+y) I''(y) + 2s(1+2y)(3+2y) I'(y)+ 3s(5s-1)(1+y) I(y) = 0.$$ This ODE has a regular singular point at $y=0$, where the two exponents are $0$ and $1-3s$. Since $I(y) = \int_0^\infty H(x,y) dx$ is analytic around $y=0$, it must be a multiple of the solution with exponent $0$. One checks that $$(1+y)^{1-5s}\pFq{2}{1}{\frac{3s}{2},\frac{5s-1}{2}}{3s}{\frac{y(2+y)}{(1+y)^2}}$$ solves the above ODE and has exponent $0$ at $y=0$. Therefore $I(y)$ equals a constant multiple of the above expression; this constant is $I(0)$. 
		Finally, $$I(0) = 2^{-s} \int_0^\infty x^{-s} (1+x)^{-3s} (x+2)^{-s} dx = 2^{-1-s} \frac{\Gamma(1-s)\Gamma(\frac{5s-1}{2})}{\Gamma(\frac{3s+1}{2})}$$ by Remark \ref{remark_d_12}. This concludes the proof.
	\end{proof}
	
	\begin{remark}
		In the above proof, the certificate of creative telescoping\footnote{This refers to the rational function $g(x,y)$ used in the proof.} is obtained with the help of the Mathematica package \textsf{HolonomicFunctions} developed by C. Koutschan \cite{Koutschan09}: one simply executes the following command:
		\begin{verbatim}
CreativeTelescoping[(x(1+x)(1+x+y)(2x+2+y)(x+y+2))^(-s),Der[x],{Der[y]}]//Factor\end{verbatim}
	\end{remark}
	
	\begin{corollary}\label{FG_eval}
		$$F_G(s;1-s) = \frac{\Gamma(1-s)\Gamma(\frac{5s-1}{2})}{\Gamma(\frac{3s+1}{2})} 2^{-6s} 3^{3s-1} \pFq{2}{1}{\frac{3s}{2},\frac{5s-1}{2}}{3s}{\frac34}.$$
		In particular, $F_G(\frac13;\frac23) = \dfrac{\Gamma(\frac{1}{3})^3}{\pi  3^{1/2} 2^{5/3}}.$
	\end{corollary}
	\begin{proof}
		Set $y=1$ in the above lemma and make substitution $x\mapsto 1/x$. The evaluation at $s=1/3$ follows from the following special values of $_2F_1$ (see \cite{ebisu2017special} or \cite{goursat1881equation}): 
		\begin{equation*}\pFq{2}{1}{a,\frac12}{3a}{\frac34} = \frac{2}{\sqrt{3}}\frac{\Gamma \left(a+\frac{1}{3}\right) \Gamma \left(a+\frac{2}{3}\right)}{\Gamma \left(a+\frac{1}{2}\right)^2}.\qedhere \end{equation*}
	\end{proof}
	
	\begin{corollary}\label{FB_diag_eval}
		When $\Re(s)>0$, we have
		$$F_B(s;1-s,1-s) = \frac{2^{2 s-3} \Gamma (1-s)^2 \Gamma(\frac{5 s-1}{2}) \Gamma(3 s-\frac{1}{2}) \Gamma (4 s-1)}{\sqrt{\pi } \Gamma (3 s) \Gamma(\frac{9 s-1}{2})} \pFq{3}{2}{\frac{1+s}{2},\frac{3s}{2},4s-1}{3s,\frac{9s-1}{2}}{1}.$$
		In particular, 
		$$F_B\left(\frac{1}{3};\frac23,\frac23\right) = \frac{3}{32 \pi ^2}\Gamma \left(\frac{1}{3}\right)^6.$$
	\end{corollary}
	\begin{proof}
		The evaluation at $s=1/3$ follows from Whipple's identity:
		$$\pFq{3}{2}{a,1-a,c}{e,1+2c-e}{1} = \frac{2^{1-2c}\pi \Gamma(e)\Gamma(1+2c-e)}{\Gamma(\frac{a+e}{2})\Gamma(\frac{a+1-e}{2}+c) \Gamma(\frac{1-a+e}{2}) \Gamma(1+c-\frac{a+e}{2})}.$$
		It remains to prove the formula for $F_B(s;1-s,1-s)$. Starting from the definition of $F_B(s;1-s,1-s)$, the above lemma gives the integral over $x$,
		$$\begin{aligned}F_B(s;1-s,1-s) &= 2^{-1-s} \frac{\Gamma(1-s)\Gamma(\frac{5s-1}{2})}{\Gamma(\frac{3s+1}{2})} \int_0^\infty (1+y)^{1-6s} y^{-s}(y+2)^{-s} \pFq{2}{1}{\frac{3s}{2},\frac{5s-1}{2}}{3s}{\frac{y(2+y)}{(1+y)^2}}dy \\
			&= 2^{-1-s} \frac{\Gamma(1-s)\Gamma(\frac{5s-1}{2})}{\Gamma(\frac{3s+1}{2})} \sum_{n\geq 0} \frac{(\frac{3s}{2})_n (\frac{5s-1}{2})_n}{(1)_n (3s)_n} \int_0^\infty (1+y)^{1-6s-2n} y^{-s+n} (y+2)^{-s+n} dy \\
			&= 2^{-1-s} \frac{\Gamma(1-s)\Gamma(\frac{5s-1}{2})}{\Gamma(\frac{3s+1}{2})} \sum_{n\geq 0} \frac{(\frac{3s}{2})_n (\frac{5s-1}{2})_n}{(1)_n (3s)_n} \frac{\Gamma(1+n-s)\Gamma(4s-1)}{\Gamma(n+3s)}.
		\end{aligned}$$
		Hence
		$$F_B(s;1-s,1-s) = \frac{2^{-s-2} \Gamma (1-s)^2 \Gamma \left(\frac{5 s-1}{2}\right) \Gamma (4 s-1)}{\Gamma (3 s) \Gamma \left(\frac{3 s}{2}+\frac{1}{2}\right)} \pFq{3}{2}{1-s,\frac{3 s}{2},\frac{5 s-1}{2}}{3 s,3 s}{1}.$$
		This $_3F_2$ converges absolutely\footnote{More generally, ${_3F_2}(a,b,c;d,e;1)$ converges absolutely when $\Re(d+e-a-b-c)>0$.} only for $\Re(s)>1/6$. To get the ${_3F_2}$ in the statement, which converges on a wider region $\Re(s)>0$, one uses the transformation
		\begin{equation*}\pFq{3}{2}{a,b,c}{d,e}{1} =  \frac{\Gamma(d)\Gamma(d+e-a-b-c)}{\Gamma(d-a)\Gamma(d+e-b-c)} \pFq{3}{2}{a,e-b,e-c}{e,d+e-b-c}{1}.\qedhere\end{equation*}
		
	\end{proof}
	
	\section{Witten zeta functions of rank 3}
	This section begins our investigation of rank-$3$ Witten zeta functions. They are:
	\begin{align*}\xi_{A_3}(s) &= \sum_{m_1,m_2,m_3\geq 1} \left( m_1m_2m_3 (m_1+m_2)(m_2+m_3)(m_1+m_2+m_3)\right)^{-s}, \\
		\xi_{B_3}(s) &= \sum_{m_1,m_2,m_3\geq 1}  \big[ m_1m_2m_3 (m_1+m_2)(m_2+m_3)(m_1+m_2+m_3) \\& \qquad \qquad \qquad \times(2m_2+m_3) (2m_1+2m_2+m_3)(m_1+2m_2+m_3)\big]^{-s}, \\
		\xi_{C_3}(s) &= \sum_{m_1,m_2,m_3\geq 1} \big[ m_1m_2m_3 (m_1+m_2)(m_2+m_3)(m_1+m_2+m_3) \\ & \qquad \qquad \qquad \times (m_2+2m_3)(m_1+m_2+2m_3)(m_1+2m_2+2m_3)\big]^{-s}.
	\end{align*}
	For any $f(x_1,x_2)$ satisfying the condition in \ref{def_two_var_f} from the last section, we define
	$$\xi_f(s) := \sum_{m_1,m_2,n\geq 1} \left[ n^{d_{12}+1} m_1 m_2 f\left(\frac{m_1}{n},\frac{m_2}{n}\right) \right]^{-s}.$$
	Recall the definition of $F_f(s;z_1,z_2)$ and $K_f(s;z_1,z_2)$ as in Section 7. By the inverse Mellin transform, 
	$$f(x_1,x_2)^{-s} = \fint_{(c_1)} \fint_{(c_2)} F_f(s;z_1,z_2) x_1^{-z_1} x_2^{-z_2} dz_1 dz_2, \quad x_1,x_2>0,\quad  c_i \text{ sufficiently large.}$$
	In complete analogy with the rank-$2$ cases, we have
	\begin{align*}\xi_f(s) &= \sum_{m_1,m_2,n\geq 1} \left[ n^{d_{12}+1} m_1 m_2 \right]^{-s} \fint_{(c_1)} \fint_{(c_2)} F_f(s;z_1,z_2) (\frac{m_1}{n})^{-z_1} (\frac{m_2}{n})^{-z_2} dz_1 dz_2 \\
		&=	\fint_{(c_1)}\fint_{(c_2)} \zeta((d_{12}+1)s-z_1-z_2) \zeta(s+z_1)\zeta(s+z_2)F_f(s;z_1,z_2) dz_1 dz_2,\end{align*}
	for $\Re(c_i)$ sufficiently large. When $f$ is one of the polynomials $f_A, f_B, f_C$ defined in (\ref{f_ABC_rank3}), $\xi_f(s)$ is the corresponding Witten zeta function.
	~\\[0.02in]
	
	Consider an integral of the form $$\fint_{(c_1)} \fint_{(c_2)} p(z_1,z_2) dz_1 dz_2.$$ We want to move the lines of integration from $\Re(z_1) = c_1$ to $\Re(z_1) = c_1'$ and from $\Re(z_2) = c_2$ to $\Re(z_2) = c_2'$, with $c_i' < c_i$. We assume that the integrand has sufficiently fast decay in the vertical directions. For generic $z_2$, let $P_1$ be the set of poles picked up when shifting the contour from $\Re(z_1) = c_1$ to $\Re(z_1) = c_1'$. Similarly, let $P_2$ be the set of poles picked up when shifting the contour from $\Re(z_2) = c_2$ to $\Re(z_2) = c_2'$. We also assume that $P_1$ is independent of $z_2$ and $P_2$ is independent of $z_1$. Under these assumptions, the above integral equals
	\begin{multline}\label{double-contour-shift}\sum_{s_1 \in P_1, s_2\in P_2} \Res{z_1 = s_1} \Res{z_2 = s_2} p(z_1,z_2) + \sum_{s_2 \in P_2}\fint_{(c_1')} \left(\Res{z_2 = s_2} p(z_1,z_2) \right)dz_1 \\
		+ \sum_{s_1 \in P_1} \fint_{(c_2')}\left(\Res{z_1 = s_1} p(z_1,z_2) \right) dz_2+ \fint_{(c_1')} \fint_{(c_2')} p(z_1,z_2) dz_1 dz_2.\end{multline}
	~\\[0.02in]
	Let us apply this observation to $$\begin{aligned} p=p(z_1,z_2) &:= \zeta(s+z_1)\zeta(s+z_2)\zeta((d_{12}+1)s-z_1-z_2) F_f(s;z_1,z_2) \\
		&= \frac{1}{\Gamma(s)} \Gamma \left(z_1\right) \Gamma \left(z_2\right)  \Gamma \left(d_2 s-z_2\right) \Gamma \left(s d_{12}-z_1-z_2\right)  \\ &\quad \times \Gamma \left(d_1 s-z_1\right)\zeta \left(s+z_1\right) \zeta \left(s+z_2\right) \zeta \left(\left(d_{12}+1\right) s-z_1-z_2\right) K_f(s;z_1,z_2),
	\end{aligned}$$ 
	Because $K_f(s;z_1,z_2)$ is entire, singularities of the integrand come solely from gamma and zeta factors. Similar to the $2$-dimensional case before, one can draw a (now $3$-dimensional) singularity diagram. When shifting the contour, the poles that we picked up are $P_1 = P_2 = \{1-s,0,-1,\ldots,-M\}$, where $M$ is a fixed non-negative integer. A total of 9 terms will appear in expansion (\ref{double-contour-shift}), namely \begin{multline}\xi_f(s) = 
		\Res{z_1=1-s}\Res{z_2=1-s} p  + \sum_{i=0}^M \Res{z_1=1-s} \: \Res{z_2=-i} p + \sum_{i=0}^M \Res{z_2=1-s} \: \Res{z_1=-i} p  +
		\sum_{i,j=0}^M \Res{z_1=-i} \: \Res{z_2=-j} p  \\
		+ \fint_{(-M-\varepsilon)} \left(\Res{z_1=1-s} p \right) dz_2 + \fint_{(-M-\varepsilon)} \left(\Res{z_2=1-s}p \right) dz_1 \\
		+ \sum_{i=0}^M  \fint_{(-M-\varepsilon)} \left(\Res{z_1=-i} p \right)dz_2 
		+ \sum_{i=0}^M  \fint_{(-M-\varepsilon)} \left(\Res{z_2=-i} p \right) dz_1 + \fint_{(-M-\varepsilon)} \fint_{(-M-\varepsilon)} p dz_1 dz_2.\end{multline}
	We shall denote these 9 terms by $T_f^1(s),\ldots,T_f^9(s)$. More precisely, we write
	\begin{equation}\label{9_term_expr}\xi_f(s) = \sum_{k=1}^9 T^k_f(s),\end{equation}
	with {\small \allowdisplaybreaks  \begin{align*}
			T^1_f(s) &:= \Res{z_1=1-s} \: \Res{z_2=1-s} p = F_f(s;1-s,1-s) \zeta((d_{12}+3)s-2) \\
			T^2_f(s) &:=  \sum_{i=0}^M \Res{z_1=1-s} \: \Res{z_2=-i} p \\ &= \frac{\Gamma (1-s) \Gamma((d_1+1) s-1) }{\Gamma(s)}\sum_{i=0}^M \frac{(-1)^i}{i!} \Gamma((1+d_{12})s+i-1)\Gamma(i+d_2s) \zeta (s-i) \zeta ((d_{12}+2)s+i-1) K_f(s;1-s,-i) \\
			&= \frac{\Gamma (1-s) \Gamma((d_1+1) s-1) }{\Gamma(s)}\sum_{i=0}^M \zeta (s-i) \zeta ((d_{12}+2)s+i-1) [K_f^{x_1}(s;1-s)][x_2^i] \\
			T^3_f(s) &:= \sum_{i=0}^M \Res{z_2=1-s}  \: \Res{z_1=-i} p \\
			&=\frac{\Gamma (1-s) \Gamma((d_2+1) s-1) }{\Gamma(s)}\sum_{i=0}^M \zeta (s-i) \zeta ((d_{12}+2)s+i-1) [K_f^{x_2}(s;1-s)][x_1^i] \\
			T^4_f(s) &:= \sum_{0\leq i,j\leq M}  \Res{z_1=-i} \: \Res{z_2=-j} p \\ &=  \sum_{0\leq i,j\leq M} \zeta(s-i)\zeta(s-j) \zeta(i+j+(d_{12}+1)s) [f(x_1,x_2)^{-s}][x_1^i x_2^j]\\
			T^5_f(s) &:= \fint_{(-M-\varepsilon)} \left(\Res{z_1=1-s} p \right) dz_2 = \frac{\Gamma (1-s)\Gamma((d_1+1) s-1)}{\Gamma(s)} \\  & \qquad \times \fint_{(-M-\varepsilon)} \Gamma(z_2)  \zeta(s+z_2) \zeta((d_{12}+2) s-z_2-1) \Gamma(s d_2-z_2) \Gamma ((d_{12}+1) s-z_2-1) K_f(s;1-s,z_2) dz_2  \\
			T^6_f(s) &:= \fint_{(-M-\varepsilon)} \left(\Res{z_2=1-s}p \right) dz_1 = \frac{\Gamma (1-s)\Gamma((d_2+1) s-1)}{\Gamma(s)} \\ &  \qquad \times \fint_{(-M-\varepsilon)} \Gamma(z_1)  \zeta(s+z_1) \zeta((d_{12}+2) s-z_1-1) \Gamma(s d_1-z_1) \Gamma ((d_{12}+1) s-z_1-1) K_f(s;z_1,1-s) dz_1 \\
			T^7_f(s) &:= \sum_{i=0}^M  \fint_{(-M-\varepsilon)} \left(\Res{z_1=-i} p \right)dz_2  \\
			&= \sum_{i=0}^M \frac{(-1)^i\zeta(s-i) \Gamma(i+d_1s)}{i! \Gamma(s)} \fint_{(-M-\varepsilon)} \Gamma(d_2s-z_2)\Gamma(z_2) \Gamma(i+sd_{12}-z_2)  \zeta((d_{12}+1)s-z_2+i) \zeta(s+z_2) K_f(s;-i,z_2) dz_2 \\
			&= \sum_{i=0}^M \frac{\zeta(s-i)}{\Gamma(s)} \fint_{(-M-\varepsilon)} \Gamma(d_2s-z_2)\Gamma(z_2) \zeta((d_{12}+1)s-z_2+i) \zeta(s+z_2) [K_f^{x_2}(s;z_2)][x_1^i] dz_2 \\
			T^8_f(s) &:= \sum_{i=0}^M  \fint_{(-M-\varepsilon)} \left(\Res{z_2=-i} p \right) dz_1 \\
			&= \sum_{i=0}^M \frac{\zeta(s-i)}{\Gamma(s)} \fint_{(-M-\varepsilon)} \Gamma(d_1s-z_1)\Gamma(z_1) \zeta((d_{12}+1)s-z_1+i) \zeta(s+z_1) [K_f^{x_1}(s;z_1)][x_2^i] dz_1 \\
			T^9_f(s) &:= \fint_{(-M-\varepsilon)} \fint_{(-M-\varepsilon)} p dz_1 dz_2 \\
			&= \frac{1}{\Gamma(s)} \fint_{(-M-\varepsilon)} \fint_{(-M-\varepsilon)} \Gamma \left(z_1\right) \Gamma \left(z_2\right)  \Gamma\left(d_1 s-z_1\right) \Gamma \left(d_2 s-z_2\right) \Gamma \left(s d_{12}-z_1-z_2\right)  \zeta \left(s+z_1\right) \\  & \qquad \times \zeta \left(s+z_2\right) \zeta \left(\left(d_{12}+1\right) s-z_1-z_2\right) K_f(s;z_1,z_2) dz_1 dz_2
	\end{align*} }
	where we have used Proposition \ref{rank_3_non-negative-int_red} to simplify $T_f^2(s),T_f^3(s),T_f^7(s),T_f^8(s)$ and Proposition \ref{rank_3_non-negative-int_red2} to simplify $T_f^4(s)$. This $9$-term formula for $\xi_f(s)$ is the rank-$3$ generalization of the $3$-term formula (\ref{continuation_M}). 
	
	For fixed $N\in \mathbb{R}$, all (single or double) integrals occurring in the above formula are analytic on $\Re(s)>N$ when $M$ is chosen sufficiently large. This provides the analytic continuation of $\xi_f(s)$. Also note that $T_f^7(s), T_f^8(s), T_f^9(s)$ can be made analytic at any given $s\in \mathbb{C}$ with $\Re(s)<1$ by choosing $M$ sufficiently large.\par
	As in the rank-$2$ case, the term $T_f^1(s) = F_f(s;1-s,1-s)\zeta((d_{12}+3)s-2)$ dictates the pole at the abscissa of convergence.
	\begin{corollary}\label{cor:residue_abs}
		The abscissa of convergence of $\xi_f(s)$ is $\frac{3}{d_{12}+3}$. It has a simple pole at this point with residue
		$$\frac{1}{d_{12}+3} F_f\left(\frac{3}{d_{12}+3};\frac{d_{12}}{d_{12}+3},\frac{d_{12}}{d_{12}+3}\right) = \frac{1}{d_{12}+3}\int_0^\infty\int_0^\infty (x_1x_2f(x_1,x_2))^{-3/(d_{12}+3)} dx_1dx_2.$$
	\end{corollary}
	\begin{proof}
	Here we recall the inequalities $3d_1 > d_{12}$ and $3d_2 > d_{12}$, which are part of the assumptions defining $f(x_1,x_2)$ in \eqref{def_two_var_f}. One checks easily that these two inequalities are equivalent to the assertion that the point $\left(\frac{3}{d_{12}+3},\frac{d_{12}}{d_{12}+3},\frac{d_{12}}{d_{12}+3}\right)$ lies in the absolute-convergence region defined in \eqref{2D_int_conv_range}. Therefore $F_f(s;1-s,1-s)$ is analytic at this point, and the displayed quantity equals the residue of $T_f^1(s)$ there.\par
	Moreover, from the expression of $T_f^2(s), \ldots T_f^6(s)$, we see the locations of their largest poles are
	$$\max\left(\frac{1}{d_1+1},\frac{2}{d_{12}+2}\right),\quad \max\left(\frac{1}{d_2+1},\frac{2}{d_{12}+2}\right),\quad \frac{1}{d_{12}+1},\quad \frac{1}{d_1+1},\quad \frac{1}{d_2+1},$$
	respectively. Since $3d_1 > d_{12}$ and $3d_2 > d_{12}$, $\frac{3}{d_{12}+3}$ is strictly larger than all of these numbers, so $\frac{3}{d_{12}+3}$ is indeed the abscissa of convergence of $\xi_f(s)$.
	\end{proof}

	\subsection{Poles and residues of $\xi_{A_3}(s)$}
	For $f(x_1,x_2) = f_A(x_1,x_2) = (1+x_1)(1+x_2)(1+x_1+x_2)$, we have $(d_1,d_2,d_{12}) = (2,2,3)$. Because $f_A(x_1,x_2)$ is symmetric, $F_A(s;z_1,z_2)$ and $K_A(s;z_1,z_2)$ are symmetric with respect to their last two arguments. We write $T_f^1,\ldots,T_f^9$ as $T_A^1,\ldots,T_A^9$. By this symmetry, we have $T_A^2 = T_A^3, T_A^5 = T_A^6$, and $T_A^7 = T_A^8$. Therefore 
	\begin{equation}\label{9_term_expr_A3}\xi_{A_3}(s) = T_A^1(s) + 2T_A^2(s) + T_A^4(s) + 2T_A^5(s) + 2T_A^7(s) + T_A^9(s), \end{equation}
	where
	\begin{align*}T_A^1(s) &= F_A(s;1-s,1-s) \zeta(6s-2) = \frac{2^{2 s-2} \Gamma (1-s)^2 \Gamma \left(\frac{3 s}{2}-\frac{1}{2}\right) \Gamma \left(\frac{5 s}{2}-1\right)}{\Gamma \left(\frac{s}{2}+\frac{1}{2}\right) \Gamma \left(\frac{3 s}{2}\right)} \zeta(6s-2)  \qquad \text{ (by Proposition \ref{F_Adiag_eval})} \\
		T_A^2(s) &= \frac{\Gamma (1-s) \Gamma(3 s-1) }{\Gamma(s)}\sum_{i=0}^M \zeta (s-i) \zeta (5s+i-1) [K_A^{x_1}(s;1-s)][x_2^i]\\
		T_A^4(s) &= \sum_{0\leq i,j\leq M} \zeta(s-i)\zeta(s-j) \zeta(i+j+4s) [f_A(x_1,x_2)^{-s}][x_1^i x_2^j] \\
		T_A^5(s) &= \frac{\Gamma (1-s)\Gamma(3 s-1)}{\Gamma(s)}  \fint_{(-M-\varepsilon)} \Gamma(z_2)   \Gamma (4s-z_2-1) \Gamma(2s-z_2) \zeta(s+z_2) \zeta(5s-z_2-1)K_A(s;1-s,z_2) dz_2 \\
		T_A^7(s)  &= \sum_{i=0}^M \frac{\zeta(s-i)}{\Gamma(s)} \fint_{(-M-\varepsilon)} \Gamma(2s-z_2)\Gamma(z_2) \zeta(4s-z_2+i) \zeta(s+z_2) [K_A^{x_2}(s;z_2)][x_1^i] dz_2 \\
		T_A^9(s) &= \frac{1}{\Gamma(s)} \fint_{(-M-\varepsilon)} \fint_{(-M-\varepsilon)} \Gamma \left(z_1\right) \Gamma \left(z_2\right)  \Gamma \left(2s-z_2\right) \Gamma \left(3s -z_1-z_2\right) \\ &\quad  \times \Gamma \left(2s-z_1\right)\zeta \left(s+z_1\right) \zeta \left(s+z_2\right) \zeta \left(4s-z_1-z_2\right) K_A(s;z_1,z_2) dz_1 dz_2
	\end{align*}
	For large $M$, the integrals in $T_A^5, T_A^7, T_A^9$ are analytic on $\Re(s)>0$. 
	
	\begin{theorem}\label{A3_residues}
		(a) The Witten zeta function $\xi_{A_3}(s)$ has four poles with $\Re(s)>0$, all of which are simple. Their respective residues are$$\begin{aligned}
			\Res{s=1/2} \xi_{A_3}(s) &= \frac{\Gamma(\frac14)^4}{24\pi}, \\  \Res{s=2/5} \xi_{A_3}(s) &= \frac{\left(\sqrt{5}+5\right) \Gamma \left(\frac{1}{5}\right) \Gamma \left(\frac{3}{5}\right)}{10 \Gamma \left(\frac{4}{5}\right)} \zeta \left(\frac{2}{5}\right), \\
			\Res{s=1/3} \xi_{A_3}(s) &= \frac{2}{3}\xi_{A_2}\left(\frac{1}{3}\right), \\
			\Res{s=1/4} \xi_{A_3}(s) &= \frac{1}{4}\zeta \left(\frac14\right)^2.
		\end{aligned}$$
		(b) All poles of $\xi_{A_3}(s)$ on the complex plane are simple.
	\end{theorem}
	\begin{proof}
		(a) In the region $\Re(s)>0$, \begin{itemize}
			\item $T_A^1(s)$ has poles at $s=1/4,2/5$ and $1/2$.
			\item $T_A^2(s)$ has poles at $s=1/3, 1/5, 2/5$, where $s=1/3$ comes from $\Gamma(3s-1)$, $s=1/5, 2/5$ comes from $\zeta(5s+i-1)$ for $i=0$ and $1$.
			\item $T_A^4(s)$ has a pole at $s=1/4$ coming from $\zeta(4s)$.
			\item $T_A^5(s)$ has a pole at $s=1/3$ coming from $\Gamma(3s-1)$.
			\item $T_A^7(s)$ and $T_A^9(s)$ are analytic.
		\end{itemize}
		In all cases, the poles are at most simple and
		$$\begin{aligned}
			\Res{s=1/2} \xi_{A_3}(s) &= \Res{s=1/2} T_A^1(s), \\
			\Res{s=2/5} \xi_{A_3}(s) &= \Res{s=2/5} T_A^1(s) +  2\Res{s=2/5}T_A^2(s), \\
			\Res{s=1/5} \xi_{A_3}(s) &= \Res{s=1/5} 2T_A^2(s), \\
			\Res{s=1/3} \xi_{A_3}(s) &= \Res{s=1/3} T_A^1(s) + 2\Res{s=1/3} T_A^2(s) + 2\Res{s=1/3}T_A^5(s),\\
			\Res{s=1/4} \xi_{A_3}(s) &= \Res{s=1/4} T_A^4(s).
		\end{aligned}$$
		Next we show that the residue at $s=1/5$ is actually zero, so this point is not a pole of $\xi_{A_3}(s)$. Indeed, $$\Res{s=1/5} \xi_{A_3}(s) = \Res{s=1/5} 2T_A^2(s) = 2\Res{s=1/5} \frac{\Gamma(1-s)\Gamma(3s-1)}{\Gamma(s)} \zeta(s-1)\zeta(5s) [K_A^{x_1}(s;1-s)][x_2^1].$$
		From Proposition \ref{KA3_non-negative_int}, one calculates $$[K_A^{x_1}(s;1-s)][x_2^1] = \frac{(1-5s)\Gamma(s)}{2\Gamma(2s)}.$$ This vanishes at $s=1/5$ and thus cancels the pole from $\zeta(5s)$. \par
		As $T_A^1(s), T_A^2(s), T_A^4(s)$ can be expressed in terms of gamma and zeta functions, it is easy to calculate the residues at $s=1/2, 2/5, 1/4$. We concentrate on the residue at $s=1/3$, which is slightly non-obvious. 
		We have \begin{align*}\Res{s=1/3} \xi_{A_3}(s) &= \Res{s=1/3} T_A^1(s) + 2\Res{s=1/3} T_A^2(s) + 2\Res{s=1/3}T_A^5(s) \\ & = -\frac{\Gamma \left(-\frac{1}{6}\right) \Gamma \left(\frac{2}{3}\right)}{6 \sqrt[3]{2} \sqrt{\pi }} +\frac{2\Gamma(\frac23)}{3\Gamma(\frac13)} \sum_{i=0}^M \frac{(-1)^i}{i!} \Gamma(\frac13+i)\Gamma(\frac23+i) \zeta (\frac13-i) \zeta (\frac23+i) K_A(\frac13;\frac23,-i)\\ &\quad +\frac{2\Gamma(\frac23)}{3\Gamma(\frac13)} \fint_{(-M-\varepsilon)} \Gamma(z_2)  \Gamma(\frac13-z_2) \Gamma (\frac23-z_2) \zeta(\frac13+z_2) \zeta(\frac23-z_2) K_A(\frac13;\frac23,z_2) dz_2 \\
			&=-\frac{\Gamma \left(-\frac{1}{6}\right) \Gamma \left(\frac{2}{3}\right)}{6 \sqrt[3]{2} \sqrt{\pi }}  + \frac{2\Gamma(\frac23)}{3\Gamma(\frac13)} \fint_{(0+\varepsilon)} \Gamma(z_2)  \Gamma(\frac13-z_2) \Gamma (\frac23-z_2) \zeta(\frac13+z_2) \zeta(\frac23-z_2) K_A(\frac13;\frac23,z_2) dz_2 
			\\ & = -\frac{\Gamma \left(-\frac{1}{6}\right) \Gamma \left(\frac{2}{3}\right)}{6 \sqrt[3]{2} \sqrt{\pi }}  + \frac{2}{3\Gamma(\frac13)} \fint_{(0+\varepsilon)} \Gamma(z_2)  \Gamma(\frac13-z_2)\zeta(\frac13+z_2) \zeta(\frac23-z_2) dz_2,
		\end{align*}
		where we used the equality $K_A(1/3;2/3,s) = \left(\Gamma (\frac{2}{3}) \Gamma(\frac{2}{3}-s)\right)^{-1}$, which follows from Proposition \ref{Kf_rank3_secondreduction}. On the other hand, formula (\ref{continuation_M}) gives us the following expression for $\xi_{A_2}(s)$:
		\begin{equation}\label{A2_repre}\xi_{A_2}(s) = \frac{\Gamma(1-s)\Gamma(2s-1)}{\Gamma(s)} \zeta(3s-1) + \frac{1}{\Gamma(s)} \fint_{(0+\varepsilon)} \Gamma(z)\Gamma(s-z)\zeta(s+z)\zeta(2s-z) dz.\end{equation}
		Comparing both formulas implies $\Res{s=1/3} \xi_{A_3}(s) = \frac{2}{3}\xi_{A_2}(1/3)$. \par
		(b) It is a general fact that Witten zeta functions are analytic at each non-positive integer (see \cite{au2024vanishing}). Moreover, $\xi_{A_3}(s)$ is analytic on $\Re(s)>1/2$. Hence it suffices to prove each term in the equation (\ref{9_term_expr_A3}) has at most a simple pole at \emph{non-integer} $s$ with $\Re(s)\leq 0$. Note that $T_A^9(s)$ and $T_A^7(s)$ are analytic, so we can safely ignore these two terms. \par
		$T_A^1(s)$ has only simple poles by inspection; $T_A^4(s)$ has at most a simple pole at each $s\in \mathbb{Z}/4$, coming from $\zeta(i+j+4s)$; and $T_A^3(s)$ has at most a simple pole at each $s\in \mathbb{Z}/3$, coming from $\Gamma(3s-1)$. For $T_A^2(s)$,
		$$T_A^2(s) = \frac{2\Gamma (1-s) \Gamma(3s-1)}{\Gamma(s)} \sum_{i=0}^M \zeta(s-i)\zeta(5s+i-1) [K_A^{x_1}(s;1-s)][x_2^i],$$
		a double pole can occur only when both $3s-1$ and $5s+i-1$ are integers. Since $3$ and $5$ are coprime, this forces $s\in \mathbb{Z}$, a case we have already excluded. Thus $T_A^2(s)$ also has at most a simple pole at each non-integer $s$. Consequently, $\xi_{A_3}(s)$ has only simple poles. 
	\end{proof}
	
	\begin{remark}
		The multi-variable $A_3$-zeta function $$\xi_{A_3}(s_1,\ldots,s_6) := \sum_{m_k\geq 1} \frac{1}{m_1^{s_1}m_2^{s_2}m_3^{s_3} (m_1+m_2)^{s_4} (m_2+m_3)^{s_5} (m_1+m_2+m_3)^{s_6}}$$
		has been studied in \cite{matsumoto2006witten}, where it is shown that $\xi_{A_3}(s_1,\ldots,s_6)$ extends to a meromorphic function on $\mathbb{C}^6$ with hyperplane singularities, 
		\begin{equation}\label{A3_multi_sing}\begin{cases}
				s_1+s_4+s_6 = 1-l,\quad l\in \mathbb{Z}_{\geq 0} \\
				s_3+s_5+s_6 = 1-l, \quad l\in \mathbb{Z}_{\geq 0}\\
				s_2+s_4+s_5+s_6 = 1-l, \quad l\in \mathbb{Z}_{\geq 0}\\
				s_1+s_2+s_4+s_5+s_6 = 2-l, \quad l\in \mathbb{Z}_{\geq 0}\\
				s_1+s_3+s_4+s_5+s_6 = 2-l, \quad l\in \mathbb{Z}_{\geq 0}\\
				s_2+s_3+s_4+s_5+s_6 = 2-l, \quad l\in \mathbb{Z}_{\geq 0}\\
				s_1+s_2+s_3+s_4+s_5+s_6 = 3 
		\end{cases}\end{equation}
		and each of these are true singularities. Restricting $\xi_{A_3}(s_1,\ldots,s_6)$ to the diagonal gives $\xi_{A_3}(s)$. In view of (\ref{A3_multi_sing}), $\xi_{A_3}(s)$ being analytic at $s=1/5$ is perhaps a bit surprising. One also expects that many candidate poles of $\xi_{A_3}(s)$ obtained by restricting (\ref{A3_multi_sing}) to the diagonal are actually not true poles. \par 
		We already saw this phenomenon for $G_2$. Below we will see $\xi_{B_3}(s), \xi_{C_3}(s)$ also skip poles in comparison to their multi-variable counterpart: both $\xi_{B_3}(s), \xi_{C_3}(s)$ are analytic at $s=1/8$.
	\end{remark}
	
	
	\subsection{Poles and residues of $\xi_{B_3}(s)$ and $\xi_{C_3}(s)$}
	Consider $$\begin{aligned}f_B(x_1,x_2) &= (1 + x_1) (1 + x_2) (1 + x_1 + x_2) (2 + x_2) (2 x_1 + 2 + x_2) (x_1 + x_2 + 2),\\
		f_C(x_1,x_2)&= (1+x_1)(1+x_2)(1+x_1+x_2)(1+2x_2)(x_1+1+2x_2)(x_1+2+2x_2).\end{aligned}$$ In both cases, we have $(d_1,d_2,d_{12}) = (4,5,6)$. Using our notation in equation (\ref{9_term_expr}), we write \begin{equation}\label{9_term_expr_B3}\xi_{B_3}(s) = \sum_{k=1}^9 T_B^k(s), \qquad \xi_{C_3}(s) = \sum_{k=1}^9 T_C^k(s).\end{equation} First we focus on $B_3$, $T_B^k(s)$ are explicitly
	{\allowdisplaybreaks \begin{align*}T_B^1(s) &= F_B(s;1-s,1-s) \zeta(9s-2) \\  &= \frac{2^{2 s-3} \Gamma (1-s)^2 \Gamma(\frac{5 s-1}{2}) \Gamma(3 s-\frac{1}{2}) \Gamma (4 s-1)}{\sqrt{\pi } \Gamma (3 s) \Gamma(\frac{9 s-1}{2})} \pFq{3}{2}{\frac{1+s}{2},\frac{3s}{2},4s-1}{3s,\frac{9s-1}{2}}{1} \zeta(9s-2)  \qquad \text{ (by Corollary \ref{FB_diag_eval})} \\
			T_B^2(s) &= \frac{\Gamma (1-s) \Gamma(5s-1) }{\Gamma(s)}\sum_{i=0}^M \zeta (s-i) \zeta (8s+i-1) [K_B^{x_1}(s;1-s)][x_2^i]\\
			T_B^3(s) &= \frac{\Gamma (1-s) \Gamma(6s-1) }{\Gamma(s)}\sum_{i=0}^M \zeta (s-i) \zeta (8s+i-1) [K_B^{x_2}(s;1-s)][x_1^i]\\
			T_B^4(s) &= \sum_{0\leq i,j\leq M} \zeta(s-i)\zeta(s-j) \zeta(i+j+7s) [f_B(x_1,x_2)^{-s}][x_1^i x_2^j] \\
			T_B^5(s) &= \frac{\Gamma (1-s)\Gamma(5s-1)}{\Gamma(s)}  \fint_{(-M-\varepsilon)} \Gamma(z_2)   \Gamma (8s-z_2-1) \Gamma(5s-z_2) \zeta(s+z_2) \zeta(7s-z_2-1)K_B(s;1-s,z_2) dz_2 \\
			T_B^6(s)  &= \frac{\Gamma (1-s)\Gamma(6s-1)}{\Gamma(s)}  \fint_{(-M-\varepsilon)} \Gamma(z_1)   \Gamma (7s-z_1-1)\Gamma(5s-z_1)  \zeta(s+z_1) \zeta(8s-z_1-1)K_B(s;z_1,1-s) dz_1  \\
			T_B^7(s) &= \sum_{i=0}^M \frac{\zeta(s-i)}{\Gamma(s)} \fint_{(-M-\varepsilon)} \Gamma(5s-z_2)\Gamma(z_2) \zeta(7s-z_2+i) \zeta(s+z_2) [K_B^{x_2}(s;z_2)][x_1^i] dz_2 \\
			T_B^8(s) &= \sum_{i=0}^M \frac{\zeta(s-i)}{\Gamma(s)} \fint_{(-M-\varepsilon)} \Gamma(4s-z_1)\Gamma(z_1) \zeta(7s-z_1+i) \zeta(s+z_1) [K_B^{x_1}(s;z_1)][x_2^i] dz_1 \\
			T_B^9(s) &= \frac{1}{\Gamma(s)} \fint_{(-M-\varepsilon)} \fint_{(-M-\varepsilon)} \Gamma \left(z_1\right) \Gamma \left(z_2\right)  \Gamma \left(5 s-z_2\right) \Gamma \left(6s -z_1-z_2\right) \\ &\quad \times \Gamma \left(4 s-z_1\right)\zeta \left(s+z_1\right) \zeta \left(s+z_2\right) \zeta \left(7s-z_1-z_2\right) K_B(s;z_1,z_2) dz_1 dz_2.
	\end{align*} }
	For large $M$, the integrals in $T_B^5(s), T_B^6(s), T_B^7(s), T_B^8(s), T_B^9(s)$ are analytic on $\Re(s)>0$. 
	
	\begin{theorem}\label{B3_residues}
		(a) The Witten zeta function $\xi_{B_3}(s)$ has five poles with $\Re(s)>0$, all of which are simple. Their respective residues are {\allowdisplaybreaks \begin{align*}
				\Res{s=1/3} \xi_{B_3}(s) &= \frac{\Gamma(\frac13)^6}{96\pi^2}, \\  \Res{s=1/4} \xi_{B_3}(s) &= \frac{\sqrt[4]{2}+\csc(\pi/8)}{32\sqrt{\pi}} \Gamma\left(\frac18\right)\Gamma\left(\frac38\right) \zeta\left(\frac{1}{4}\right), \\
				\Res{s=1/5} \xi_{B_3}(s) &= \frac{1}{5\times 2^{1/5}}\xi_{B_2}\left(\frac{1}{5}\right), \\
				\Res{s=1/6} \xi_{B_3}(s) &= \frac{1}{6}\xi_{A_2}\left(\frac{1}{6}\right), \\
				\Res{s=1/7} \xi_{B_3}(s) &= \frac{1}{7\times 2^{3/7}}\zeta\left(\frac{1}{7}\right)^2.
		\end{align*} }
		(b) On the complex plane, $\xi_{B_3}(s)$ has possibly double poles at $s=-1/2,-3/2,\ldots$, all other poles are simple.
	\end{theorem}
	\begin{proof}
		(a) In the region $\Re(s)>0$, \begin{itemize}[leftmargin=*]
			\item $T_B^1(s)$ has poles at $s=1/3,1/4, 1/5$ and $1/6$,
			\item $T_B^2(s)$ has poles at $s=1/5, 1/4$ and $1/8$, where $s=1/5$ comes from $\Gamma(5s-1)$, $s=1/4, 1/8$ comes from $\zeta(8s+i-1)$ for $i=0$ and $1$,
			\item $T_B^3(s)$ has poles at $s=1/6, 1/4$ and $1/8$, where $s=1/6$ comes from $\Gamma(6s-1)$, $s=1/4, 1/8$ comes from $\zeta(8s+i-1)$ for $i=0$ and $1$,
			\item $T_B^4(s)$ has a pole at $s=1/7$ coming from $\zeta(7s)$,
			\item $T_B^5(s)$ has a pole at $s=1/5$ coming from $\Gamma(5s-1)$,
			\item $T_B^6(s)$ has a pole at $s=1/6$ coming from $\Gamma(6s-1)$,
			\item $T_B^7(s), T_B^8(s)$ and $T_B^9(s)$ are analytic.
		\end{itemize}
		In all cases, the poles are at most simple and
		\begin{align}\label{aux_11}\begin{split}
				\Res{s=1/3} \xi_{B_3}(s) &= \Res{s=1/3} T_B^1(s), \\
				\Res{s=1/4} \xi_{B_3}(s) &= \Res{s=1/4} T_B^1(s) +  \Res{s=1/4}T_B^2(s) + \Res{s=1/4}T_B^3(s), \\
				\Res{s=1/8} \xi_{B_3}(s) &= \Res{s=1/8} T_B^2(s) + \Res{s=1/8}T_B^3(s), \\
				\Res{s=1/5} \xi_{B_3}(s) &= \Res{s=1/5} T_B^1(s) + \Res{s=1/5}T_B^2(s) + \Res{s=1/5}T_B^5(s), \\
				\Res{s=1/6} \xi_{B_3}(s) &= \Res{s=1/6} T_B^1(s) + \Res{s=1/6}T_B^3(s) + \Res{s=1/6}T_B^6(s), \\
				\Res{s=1/7} \xi_{B_3}(s) &= \Res{s=1/7} T_B^4(s).\end{split}
		\end{align}
		The residue at the abscissa of convergence $s=1/3$ is evaluated in Corollary \ref{FB_diag_eval}. The residue at $s=1/7$ is also easy. For $s=1/4$, we calculate
		\begin{align}\label{aux_10}
			\Res{s=1/4} \xi_{B_3}(s) &= \frac{\zeta \left(\frac{1}{4}\right) \Gamma \left(\frac{1}{8}\right) \Gamma \left(\frac{1}{4}\right) \Gamma \left(\frac{3}{4}\right)}{16 \sqrt{2 \pi } \Gamma \left(\frac{5}{8}\right)} + \frac{1}{8} \zeta \left(\frac{1}{4}\right) \Gamma \left(\frac{3}{4}\right) K_{f_B(x_1,0)}(\frac14;\frac34) + \frac{\sqrt{\pi } \zeta \left(\frac{1}{4}\right) \Gamma \left(\frac{3}{4}\right)}{8 \Gamma \left(\frac{1}{4}\right)} K_{f_B(0,x_2)}(\frac14;\frac34) \nonumber \\
			&= \frac{ \sqrt{\pi} \zeta \left(\frac{1}{4}\right) \Gamma \left(\frac{1}{8}\right)}{16\Gamma \left(\frac{5}{8}\right)} + \frac{\zeta \left(\frac{1}{4}\right)}{8}\int_0^\infty (2 x_1(1 + x_1)^2 (2 + x_1) (2 + 2 x_1))^{-1/4} dx_1  \nonumber \\ &\quad + \frac{\zeta \left(\frac{1}{4}\right)}{8} \int_0^\infty  (x_2(1 + x_2)^2 (2 + x_2)^3)^{-1/4} dx_2  \nonumber\\
			&= \frac{ \sqrt{\pi} \zeta \left(\frac{1}{4}\right) \Gamma \left(\frac{1}{8}\right)}{16\Gamma \left(\frac{5}{8}\right)} + \frac{\zeta \left(\frac{1}{4}\right)}{8} I_1 + \frac{\zeta \left(\frac{1}{4}\right)}{8} I_2,
		\end{align}
		where $I_1$ and $I_2$ denote the last two definite integrals. They can be calculated in terms of the gamma function: $I_1 = \frac{\pi ^{3/2}}{\sqrt[4]{2} \Gamma \left(\frac{5}{8}\right) \Gamma \left(\frac{7}{8}\right)}, I_2 = \sqrt{\sqrt{2}-1} I_1$. Substituting them back gives the value of $\Res{s=1/4}\xi_{B_3}(s)$.\par
		Next we claim $T_B^2(s), T_B^3(s)$ are actually both analytic at $s=1/8$, so this point is not a pole of $\xi_{B_3}(s)$. Indeed,
		$$\begin{aligned}\Res{s=1/8} T_B^2(s) &= \Res{s=1/8}\frac{\Gamma(1-s)\Gamma(5s-1)}{\Gamma(s)} \zeta(s-1)\zeta(8s) [K_B^{x_1}(s;1-s)][x_2^1]\\
			&= \Res{s=1/8} \zeta(s-1)\zeta(8s) [F_B^{x_1}(s;1-s)][x_2^1];\\
			\Res{s=1/8} T_B^3(s) &= \Res{s=1/8}\frac{\Gamma(1-s)\Gamma(6s-1)}{\Gamma(s)} \zeta(s-1)\zeta(8s) [K_B^{x_2}(s;1-s)][x_1^1] \\
			&= \Res{s=1/8} \zeta(s-1)\zeta(8s) [F_B^{x_2}(s;1-s)][x_1^1].\end{aligned}$$
		Extracting the $x_1^1$ and $x_2^1$ coefficients gives
		$$\begin{aligned}[F_B^{x_1}(s;1-s)][x_2^1] &= -2^{-2 s-1} s  \int_0^\infty \left[(x_1+1)^2 \left(3 x_1^2+14 x_1+14\right) x_1^{-s} \left((x_1+1)^3 (x_1+2)\right)^{-s-1} \right] dx_1, \\ [F_B^{x_2}(s;1-s)][x_1^1] &= -s \int_0^\infty \left[(x_2+1) (x_2+2)^2 \left(x_2^2+7 x_2+7\right) x_2^{-s} \left((x_2+1)^2 (x_2+2)^3\right)^{-s-1} \right] dx_2.
		\end{aligned}$$
		By Remark \ref{remark_d_12}, they can be expressed in terms of the $_2F_1$ function. Their vanishing (after analytic continuation) at $s=1/8$ can then be checked easily; we omit the details.\footnote{If one trusts Mathematica, the vanishing can be checked by the command \begin{verbatim} Integrate[x^-s(1+x)^2((1+x)^3(2+x))^(-1-s)(14+14x+3x^2),{x,0,Infinity},GenerateConditions -> False] 
				/. s->1/8 // FullSimplify\end{verbatim} 
		}
		
		At $s=1/6$, $$\begin{aligned}\Res{s=1/6} \xi_{B_3}(s) &= \Res{s=1/6} T_B^1(s) + \Res{s=1/6}T_B^3(s) + \Res{s=1/6}T_B^6(s) \\
			&= \frac{\zeta \left(-\frac{1}{2}\right) \Gamma \left(-\frac{1}{3}\right) \Gamma \left(\frac{5}{6}\right)}{12\ 2^{2/3} \sqrt{\pi }} + \frac{1}{6\Gamma(\frac16)} \fint_{(0+\varepsilon)}  \Gamma (\frac{1}{6}-z_1) \Gamma (z_1) \zeta(\frac{1}{3}-z_1) \zeta (z_1+\frac{1}{6}) dz_1\end{aligned},$$
		where we have used $K_B(1/6;z_1,5/6) = \left( \Gamma(\frac{5}{6}) \Gamma (\frac{2}{3}-z_1) \right)^{-1}$ (Proposition \ref{Kf_rank3_secondreduction}) and an argument similar to the $A_3$ case above. From equation (\ref{A2_repre}), we see that this equals $\frac{1}{6}\xi_{A_2}\left(\frac16\right)$. \par
		Finally, at $s=1/5$, an argument similar to the $A_3$-case above shows 
		$$\Res{s=1/5} \xi_{B_3}(s) = \Res{s=1/5} T_B^1(s) +  \frac{\Gamma (\frac45)}{5\Gamma(\frac15)}  \fint_{(0+\varepsilon)} \Gamma(z_2)   \Gamma (\frac35-z_2) \Gamma(1-z_2) \zeta(\frac15+z_2) \zeta(\frac25-z_2)K_B(\frac{1}{5};\frac{4}{5},z_2) dz_2.$$
		Proposition \ref{Kf_rank3_secondreduction} says $$K_B(1/5;4/5,z) = \frac{1}{\Gamma(\frac 45)\Gamma(1-z)} K_{2(1+y)(2+y)}^y(\frac15;z) = \frac{2^{-1/5}}{\Gamma(\frac 45)\Gamma(1-z)} K_{(1+y)(2+y)}^y(\frac15;z).$$
		Hence $$\begin{aligned}\Res{s=1/5} \xi_{B_3}(s) &=-\frac{\sqrt{\pi } \zeta \left(-\frac{1}{5}\right) \Gamma \left(\frac{4}{5}\right)}{2^{4/5} \Gamma \left(\frac{3}{10}\right)} + \frac{2^{-1/5}}{5\Gamma(\frac15)} \fint_{(0+\varepsilon)}\Gamma(\frac{2}{5}-z) \Gamma (z) \zeta(\frac{3}{5}-z) \zeta(z+\frac{1}{5}) K^y_{(1+y)(2+y)}(\frac{1}{5},z) dz \\
			&= -\frac{\sqrt{\pi } \zeta \left(-\frac{1}{5}\right) \Gamma \left(\frac{4}{5}\right)}{2^{4/5} \Gamma \left(\frac{3}{10}\right)} + \frac{2^{-1/5}}{5\Gamma(\frac15)} \fint_{(0+\varepsilon)}\Gamma(\frac{2}{5}-z) \Gamma (z) \zeta(\frac{1}{5}+z) \zeta(\frac35-z) K^x_{(1+x)(1+2x)}(\frac{1}{5},z) dz
			,\end{aligned}$$
		where we made the substitution $z\mapsto 2/5-z$ inside the integral. 
		On the other hand, formula (\ref{continuation_M}) applied to $f(x) = (1+x)(1+2x)$ gives
		\begin{equation}\label{aux_12}\xi_{B_2}(s) = \frac{\Gamma(1-s)\Gamma(3s-1)}{\Gamma(s)} \zeta(4s-1)K_f(s;1-s) + \frac{1}{\Gamma(s)} \fint_{(0+\varepsilon)} \Gamma(z)\Gamma(2s-z)\zeta(s+z)\zeta(3s-z) K_f(s;z) dz,\end{equation}
		comparing both formulas says the residue is $\frac{2^{-1/5}}{5} \xi_{B_2}(\frac{1}{5})$. \par

		(b) Similar to the $A_3$-case discussed above, we investigate the order of poles for each term in equation (\ref{9_term_expr_B3}) at \emph{non-integer} $s$ with $\Re(s)\leq 0$. \par
		The terms $T_B^7(s), T_B^8(s), T_B^9(s)$ can be ignored because they are analytic at any given $s$ for sufficiently large $M$. The terms $T_B^5(s), T_B^6(s)$, which contain the factors $\Gamma(5s-1)$ and $\Gamma(6s-1)$, can contribute at most simple poles. The term $T_B^4(s)$ contributes simple poles when $s\in \mathbb{Z}/7$. For $$T_B^2(s) = \frac{\Gamma (1-s) \Gamma(5s-1) }{\Gamma(s)}\sum_{i=0}^M \zeta (s-i) \zeta (8s+i-1) [K_B^{x_1}(s;1-s)][x_2^i],$$ a double pole can occur only when both $5s-1$ and $8s+i-1$ are integers. Because $5$ and $8$ are coprime, this forces $s\in \mathbb{Z}$, a case we have already excluded. For $$T_B^3(s) = \frac{\Gamma (1-s) \Gamma(6s-1) }{\Gamma(s)}\sum_{i=0}^M \zeta (s-i) \zeta (8s+i-1) [K_B^{x_2}(s;1-s)][x_1^i],$$ the argument does not work because $6$ and $8$ are not coprime; this term could therefore have a double pole at $s=-1/2,-3/2,\ldots$. For $T_B^1(s) =F_B(s;1-s,1-s) \zeta(9s-2)$, we rewrite it in terms of $K_B$:
		$$T_B^1(s) = \frac{\Gamma (1-s)^2 \Gamma (5 s-1) \Gamma (6 s-1) \Gamma (8 s-2) }{\Gamma (s)} K_B(s;1-s,1-s)\zeta(9s-2).$$
		Since $K_B(s;1-s,1-s)$ is entire, we see the expression above has at most double poles at $s=-1/2,-3/2,\ldots$ and has at most simple poles at other non-integer $s$.
	\end{proof}
	
	\begin{example}
		It is not difficult to calculate the leading coefficient of $\xi_{B_3}(s)$ at $s=-1/2,-3/2,\ldots$. We illustrate this for $s=-1/2$. From the proof above, the terms contributing to a double pole are $T_B^1(s)$ and $T_B^3(s)$. For the latter term, only the summand with $i=6$ contributes, so
		$$\begin{aligned}\xi_{B_3}(s) &= T_B^1(s) + T_B^3(s) +O(\frac{1}{s+\frac12}) \\
			&= F_B(s;1-s,1-s) \zeta(9s-2)+ \frac{\zeta (s-6) \zeta (8 s+5) \Gamma (1-s) \Gamma (6 s-1)}{\Gamma (s)} [K_{f_B(x_1,x_2)}^{x_2}(s;1-s)][x_1^6] + O(\frac{1}{s+\frac12}) \\
			&= -\frac{\pi  K_B(\frac{-1}{2};\frac{3}{2},\frac{3}{2}) \zeta \left(-\frac{13}{2}\right)}{43545600 \left(s+\frac{1}{2}\right)^2} -\frac{\zeta \left(-\frac{13}{2}\right)[K_{f_B(x_1,x_2)}^{x_2}(-1/2;3/2)][x_1^6]}{4608(s+\frac12)^2}+ O(\frac{1}{s+\frac12}).
		\end{aligned}$$
		For the first term on the RHS, Proposition \ref{Kf_rank3_secondreduction} gives $$K_B(-1/2;z,3/2) = -\frac{3(16 z^4+144 z^3+500 z^2+792 z+945)}{8 \sqrt{\pi } \Gamma (-z-2)}\implies K_B(-1/2;3/2,3/2) = -\frac{1204875}{128 \pi }.$$ For the second term on the RHS, Proposition \ref{Kf(s,-n)} gives $$[K_{f_B(x_1,x_2)}^{x_2}(-1/2;3/2)][x_1^6] = \left[12 \sqrt{x_1+1} (x_1^4+2 x_1^3+3 x_1^2+2 x_1+1)\right][x_1^6] = -\frac{255}{256}.$$
		Combining these two pieces of information, we obtain$$\xi_{B_3}(s) = \frac{85 \zeta \left(-\frac{13}{2}\right)}{196608 \left(s+\frac{1}{2}\right)^2} + O(\frac{1}{s+\frac12}),$$
		so $\xi_{B_3}(s)$ indeed has a double pole at $s=-1/2$.
	\end{example}
	
	\begin{theorem}\label{C3_residues}
		(a) The Witten zeta function $\xi_{C_3}(s)$ has five poles with $\Re(s)>0$, all of which are simple. Their respective residues are$$\begin{aligned}
			\Res{s=1/3} \xi_{C_3}(s) &= \frac{\Gamma(\frac13)^6}{96\pi^2}, \\  \Res{s=1/4} \xi_{C_3}(s) &= \frac{2+\csc \left(\frac{\pi }{8}\right)}{32 \sqrt[4]{2} \sqrt{\pi }} \Gamma \left(\frac{1}{8}\right) \Gamma \left(\frac{3}{8}\right) \zeta\left(\frac{1}{4}\right), \\
			\Res{s=1/5} \xi_{C_3}(s) &= \frac{1}{5}\xi_{B_2}\left(\frac{1}{5}\right), \\
			\Res{s=1/6} \xi_{C_3}(s) &= \frac{1}{6\sqrt{2}}\xi_{A_2}\left(\frac{1}{6}\right), \\
			\Res{s=1/7} \xi_{C_3}(s) &= \frac{1}{7\times 2^{1/7}}\zeta\left(\frac{1}{7}\right)^2.
		\end{aligned}$$
		(b) On the complex plane, $\xi_{C_3}(s)$ has possibly double poles at $s=-1/2,-3/2,\ldots$, all other poles are simple.
	\end{theorem}
	\begin{proof}
		Note that the positions of the poles of $T_f^k(s)$ depend only on the value of $(d_1,d_2,d_{12})$ for $f(x_1,x_2)$. To obtain the residues of $\xi_{C_3}(s)$, one simply replaces $B_3$ by $C_3$ throughout equation (\ref{aux_11}). \par
		We already noted that $F_C(s;z_1,z_2) = 2^{2s-z_2} F_B(s;z_1,z_2)$ and $K_C(s;z_1,z_2) = 2^{2s-z_2} K_B(s;z_1,z_2)$. Hence we have $T_C^1(s) = 2^{3s-1} T_B^1(s), T_C^6(s) = 2^{3s-1} T_B^6(s)$, therefore
		$$\Res{s=1/3} \xi_{C_3}(s) = \Res{s=1/3} 2^{3s-1} \xi_{B_3}(s), \quad \Res{s=1/6} \xi_{C_3}(s) = \Res{s=1/6} 2^{3s-1} \xi_{B_3}(s),$$
		which gives the residues at the points $1/3$ and $1/6$. The residue at $s=1/7$ is also easy to calculate. 
		Let us denote the summands of $T_B^2(s)$ and $T_B^3(s)$ as $$T_B^2(s) := \sum_{i=0}^M T_B^{2,i}(s),\quad T_B^3(s) := \sum_{i=0}^M T_B^{3,i}(s),$$
		Then $$T_C^2(s) = \sum_{i=0}^M 2^{2s+i} T_B^{2,i}(s), \qquad T_C^3(s) = \sum_{i=0}^M 2^{3s-1} T_B^{3,i}(s).$$
		For the residue at $s=1/4$, only the term with $i=0$ contributes, so
		$$\begin{aligned}\Res{s=1/4} \xi_{C_3}(s) &= \Res{s=1/4} 2^{3s-1} T_B^1(s) + \Res{s=1/4} 2^{2s} T_B^2(s) + \Res{s=1/4} 2^{3s-1} T_B^3(s) \\
			&= 2^{-1/4} \Res{s=1/4} T_B^1(s) + 2^{1/2}\Res{s=1/4} T_B^2(s) + 2^{-1/4}\Res{s=1/4} T_B^3(s) \\
			&= 2^{-1/4}\frac{\zeta \left(\frac{1}{4}\right) \Gamma \left(\frac{1}{8}\right)}{16\Gamma \left(\frac{5}{8}\right)} + 2^{1/2} I_1 + 2^{-1/4} I_2,\end{aligned}$$
		where $I_1$ and $I_2$ are defined in equation (\ref{aux_10}). This gives the value of $\Res{s=1/4} \xi_{C_3}(s)$. 
		
		For the residue at $s=1/8$, only the term with $i=1$ contributes, so
		$$\Res{s=1/8} \xi_{C_3}(s) = \Res{s=1/8} 2^{2s+1} T_B^2(s) + \Res{s=1/8} 2^{3s-1} T_B^3(s).$$ According to the proof above, both residues vanish, hence $\xi_{C_3}(s)$ is analytic at $s=1/8$. \par
		It remains to calculate the residue at $s=1/5$. We proceed largely in parallel with the $B_3$ case. First, $$\Res{s=1/5} \xi_{C_3}(s) = \Res{s=1/5} T_C^1(s) +  \frac{\Gamma (\frac45)}{5\Gamma(\frac15)}  \fint_{(0+\varepsilon)} \Gamma(z_2)   \Gamma (\frac35-z_2) \Gamma(1-z_2) \zeta(\frac15+z_2) \zeta(\frac25-z_2)K_C(\frac{1}{5};\frac{4}{5},z_2) dz_2,$$
		and
		$$K_C(1/5;4/5,z) = 2^{2/5-z} K_B(1/5;4/5,z) = \frac{2^{1/5-z}}{\Gamma(\frac 45)\Gamma(1-z)} K_{(1+y)(2+y)}^y(\frac15;z).$$
		Hence $$\Res{s=1/5} \xi_{C_3}(s) = 2^{-2/5} \left(-\frac{\sqrt{\pi } \zeta \left(-\frac{1}{5}\right) \Gamma \left(\frac{4}{5}\right)}{2^{4/5} \Gamma \left(\frac{3}{10}\right)}\right) + \frac{1}{5\Gamma(\frac15)} \fint_{(0+\varepsilon)}\Gamma(\frac{2}{5}-z) \Gamma (z) \zeta(\frac{1}{5}+z) \zeta(\frac35-z) 2^{1/5-z} K^y_{(1+y)(2+y)}(\frac{1}{5},z) dz$$
		Note that $$K^y_{(1+y)(2+y)}(\frac{1}{5},z) = 2^z K^y_{(1+2y)(2+2y)}(\frac{1}{5},z) = 2^{z-1/5} K^y_{(1+y)(1+2y)}(\frac{1}{5},z) ,$$
		comparison with equation (\ref{aux_12}) shows that the residue is $\frac{1}{5}\xi_{B_2}(\frac15)$.\par
		(b) The proof is the same as in the $B_3$ case, \textit{mutatis mutandis}.
	\end{proof}

	\section{$\xi(0)$ and $\xi'(0)$ for rank 3 Witten zeta functions}
	The goal of this section is to establish the following evaluations.
	\begin{theorem}
		We have
		$$\xi_{A_3}(0) = -\frac{1}{4},\qquad \xi_{B_3}(0) = \xi_{C_3}(0) = -\frac{5}{16},$$ as well as
		$$\xi'_{A_3}(0) = -\frac{3}{2}\log(2\pi),\qquad \xi'_{B_3}(0) = -\frac{45}{16}\log(2\pi)+\frac{9}{16}\log 2,\qquad \xi'_{C_3}(0) = -\frac{45}{16}\log(2\pi)+\frac{5}{8}\log 2.$$
	\end{theorem}
	
	In equation (\ref{9_term_expr}), we take $M=2$. \footnote{One can also take any positive integer $M\geq 2$, but the following calculations will be longer.} Note that $T_f^k(0) = 0$ for $k=5,6,7,8,9$: here the vanishing for $k=7,8,9$ comes from the $\Gamma(s)$ in their denominators, while that for $k=5,6$ comes from Proposition \ref{Kf(-n,z)rank3}. Define $A_f, A_f'\in \mathbb{R}$ by $$\sum_{k=1}^4 T_f^k(s) = A_f + sA_f' + O(s^2).$$
	Then $\xi_f(0) = A_f$ and \begin{multline*}\xi_f'(0) = A'_f + \sum_{k=5}^9 T_f'^k(0) \\ = A'_f - \frac{1}{d_1+1} \fint_{(-2-\varepsilon)} \Gamma(z_2) \zeta(z_2) \zeta(-z_2-1) \Gamma(-z_2) \Gamma(-z_2-1) \left(\left. \frac{d}{ds}\right\vert_{s=0}K_f(s;1-s,z_2)\right) dz_2 \\
		- \frac{1}{d_2+1} \fint_{(-2-\varepsilon)} \Gamma(z_1) \zeta(z_1) \zeta(-z_1-1) \Gamma(-z_1) \Gamma(-z_1-1) \left(\left. \frac{d}{ds}\right\vert_{s=0}K_f(s;z_1,1-s)\right) dz_1 \\
		+ \sum_{i=0}^1 \zeta(-i) \fint_{(-2-\varepsilon)} \Gamma(-z_2)\Gamma(z_2)\zeta(i-z_2)\zeta(z_2) [K^{x_2}_f(0;z_2)][x_1^i] dz_2\\
		+ \sum_{i=0}^1 \zeta(-i) \fint_{(-2-\varepsilon)} \Gamma(-z_1)\Gamma(z_1)\zeta(i-z_1)\zeta(z_1) [K^{x_1}_f(0;z_1)][x_2^i] dz_1 \\
		+ \fint_{(-2-\varepsilon)} \fint_{(-2-\varepsilon)} \Gamma (z_1) \Gamma (z_2)  \Gamma (-z_2) \Gamma (-z_1-z_2)\Gamma (-z_1)\zeta (z_1) \zeta (z_2) \zeta (-z_1-z_2) K_f(0;z_1,z_2) dz_1 dz_2.
	\end{multline*}
	Note that $[K^{x_2}_f(0;z_2)][x_1^i]$ and $[K^{x_1}_f(0;z_1)][x_2^i]$ can be calculated by Proposition \ref{Kf(0,z)}; $K_f(0;z_1,z_2)$ can be calculated by Proposition \ref{Kf(0,z)rank3}; we will also use Example \ref{derivative_zero_example} for $\left. \frac{d}{ds}\right\vert_{s=0}K_f(s;z_1,1-s)$ and $\left. \frac{d}{ds}\right\vert_{s=0}K_f(s;1-s,z_2)$. For positive $\alpha_1,\alpha_2$, we shall denote $$\mathcal{J}(\alpha_1,\alpha_2) := \fint_{(-2-\varepsilon)}\fint_{(-2-\varepsilon)} \Gamma(z_1)\Gamma(z_2)\Gamma(-z_1-z_2)\zeta(z_1)\zeta(z_2)\zeta(-z_1-z_2) \alpha_1^{-z_1} \alpha_2^{-z_2}dz_1dz_2,$$
	
	\begin{example}[The case for $A_3$]
		Thanks to Proposition \ref{KA3_non-negative_int}, we can explicitly evaluate $K_A(s;1-s,-i) = K_A(s;-i,1-s)$ in terms of the gamma and Riemann zeta function. After some calculations,  \begin{multline*}\sum_{k=1}^4 T_A^k(s) = \frac{2^{2 s-2} \Gamma (1-s)^2 \Gamma \left(\frac{3 s}{2}-\frac{1}{2}\right) \Gamma \left(\frac{5 s}{2}-1\right)}{\Gamma \left(\frac{s}{2}+\frac{1}{2}\right) \Gamma \left(\frac{3 s}{2}\right)} \zeta(6s-2) \\ + \frac{2\Gamma (1-s) \Gamma(3s-1)}{\Gamma(s)}\sum_{i=0}^2 \frac{(-1)^i}{i!} \Gamma(4s+i-1)\Gamma(i+2s) \zeta (s-i) \zeta (5s+i-1) K_A(s;1-s,-i) \\ + \sum_{0\leq i,j\leq 2} \zeta(s-i)\zeta(s-j) \zeta(i+j+4s) [f(x_1,x_2)^{-s}][x_1^i x_2^j],\end{multline*}
		when expanded up to $O(s^2)$ is 
		$$-\frac{1}{4}+\left(2\zeta '(-2)+\frac{25 \zeta'(-1)}{6}+\frac{\pi ^2}{864}-\frac{\gamma }{6}-\frac{1}{9}-\frac{71}{72} \log (2 \pi )\right)s+O\left(s^2\right).$$
		Thus $A_f = -1/4 = \xi_{A_3}(0)$ and $A'_f$ is the number inside parenthesis. Now we evaluate other terms in the equation (\ref{derivative_terms}). Note that $[K_A^{x_1}(0;z)][x_2^i] = [1+(1+x_2)^z][x_2^i]$ and by Example \ref{derivative_zero_example} and Proposition \ref{Kf(0,z)rank3}, 
		$$\left. \frac{d}{ds}\right\vert_{s=0}K_A(s;1-s,z) = -\frac{5+z}{\Gamma(-z)},\qquad K_A(0;z_1,z_2) = \frac{1}{\Gamma(-z_1)\Gamma(-z_2)}.$$
		Plugging them into equation (\ref{derivative_terms}) gives 
		\begin{multline*}\xi_{A_3}'(0) = A'_f - \frac{2}{3} \fint_{(-2-\varepsilon)} \Gamma(z_2) \zeta(z_2) \zeta(-z_2-1) \Gamma(-z_2-1) (-5-z_2) dz_2 \\
			+ 2\sum_{i=0}^1 \zeta(-i) \fint_{(-2-\varepsilon)} \Gamma(-z_2)\Gamma(z_2)\zeta(i-z_2)\zeta(z_2) [1+(1+x_1)^{z_2}][x_1^i] dz_2
			+\mathcal{J}(1,1) \\
			= A'_f - \frac{2}{3} \fint_{(-2-\varepsilon)} \Gamma(z_2) \zeta(z_2) \zeta(-z_2-1) \Gamma(-z_2-1) (-5-z_2) dz_2 - 2\fint_{(-2-\varepsilon)} \Gamma(z_2) \zeta(z_2) \zeta(-z_2) \Gamma(-z_2) dz_2 \\ - \frac{1}{6} \fint_{(-2-\varepsilon)} \Gamma(z_2) \zeta(z_2) \zeta(-z_2) \Gamma(-z_2) z_2 dz_2 + \mathcal{J}(1,1)
		\end{multline*}
		Therefore all single integrals in equation (\ref{derivative_terms}) reduce to integrals of the form
		$$\mathcal{I}(k,n,\alpha) = \fint_{(n-1-\varepsilon)} \Gamma(z)\zeta(z) \alpha^{-z} (z)_k \Gamma(n-z)\zeta(n-z) dz, \qquad n\in\{0,-1\},\quad \alpha = 1,\quad k\in \{0,1\}.$$
		which have been investigated in Section \ref{convolution_int_section}. Using the procedure described there, we can evaluate them, giving 
		$$\xi'_{A_3}(0) = A'_f + \left( -5 \zeta'(-2)-\frac{8}{3} \zeta'(-1)-\frac{1}{18}-\frac{\pi ^2}{432}+\frac{\gamma }{12}-\frac{11}{36} \log (2 \pi ) \right) + \mathcal{J}(1,1).$$
		We will find the value of $\mathcal{J}(1,1)$ at the end of this section, thereby finishing the evaluation of $\xi'_{A_3}(0)$. Note the miraculous cancellation of various terms, leaving only $\xi'_{A_3}(0) = -\frac{3}{2}\log(2\pi)$. 
	\end{example}
	
	\begin{example}[The cases for $B_3$ and $C_3$]
		It can be shown $T_B^1(s) = F_B(s;1-s,1-s) \zeta(9s-2) = \frac{225}{32}\zeta'(-2) s + O(s^2)$, we can of course also explicitly evaluate $T_B^4(s)$ up to $O(s^2)$. For $T_B^2(s), T_B^3(s)$, we need to know the expansion of $K_B^{x_1}(s;1-s)$ and $K_B^{x_2}(s;1-s)$ up to $O(s^2)$, this is given by Proposition \ref{derivativeat0}:
		$$\begin{aligned}
			K_B^{x_2}(s;1-s) &= 4 \left(x_1+2\right)+s \left(2 \left(x_1+2\right) \left(-3 \log \left(x_1+2\right)+8 \gamma -6 \log (2)\right)-2 \left(11 x_1+13\right) \log \left(x_1+1\right)\right)+O\left(s^2\right),\\
			K_B^{x_1}(s;1-s) &= \frac{5}{2} \left(x_2+2\right)+\frac{5s}{2} \left(\left(x_2+2\right) \left(3 \gamma -4 \log \left(x_2+2\right)\right)-\left(3 x_2+4\right) \log \left(x_2+1\right)\right)+O\left(s^2\right).\end{aligned}$$
		Therefore \begin{multline*}\sum_{k=1}^4 T_B^k(s) = T_B^1(s) + T_B^4(s) + \frac{\Gamma (1-s) \Gamma(5s-1) }{\Gamma(s)}\sum_{i=0}^2 \zeta (s-i) \zeta (8s+i-1) [K_B^{x_1}(s;1-s)][x_2^i] \\
			+\frac{\Gamma (1-s) \Gamma(6s-1) }{\Gamma(s)}\sum_{i=0}^2 \zeta (s-i) \zeta (8s+i-1) [K_B^{x_2}(s;1-s)][x_1^i] \\
			=\frac{-5}{16}+s\left(\frac{11}{32} \zeta '(-2)+\frac{125}{12} \zeta'(-1)+\frac{7 \pi ^2}{3456}-\frac{7 \gamma }{24}-\frac{61}{144}+\frac{9 \log (2)}{8}-\frac{119}{72} \log (2 \pi) \right) + O(s^2)\end{multline*}
		Hence $\xi_{B_3}(0) = -5/16$. Now we evaluate other terms in the equation (\ref{derivative_terms}). From Example \ref{derivative_zero_example}, Propositions \ref{Kf(0,z)} and \ref{Kf(0,z)rank3}, we have
		$$\begin{aligned}\left. \frac{d}{ds}\right\vert _{s=0} K_B(s;1-s,z) &= -\frac{2 \left(2 z+9\times 2^z+7\right)}{\Gamma (-z)},\qquad \left. \frac{d}{ds}\right\vert _{s=0} K_B(s;z,1-s) = -\frac{5 \left(z+2^{z+3}+4\right)}{2 \Gamma (-z)},\\
			K^{x_2}_B(0;z) &= 2^z (x_1+1)^z+(x_1+1)^z+(x_1+2)^z+2^z+1,\\ K^{x_1}_B(0;z) &= (x_2+1)^z+2^{-z} (x_2+2)^z+(x_2+2)^z+1,\\
			K_B(0;z_1,z_2) &= \frac{2^{z_1+z_2}+2^{z_2}+1}{\Gamma(-z_1)\Gamma(-z_2)}.
		\end{aligned}$$
		Plugging them into equation (\ref{derivative_terms}) gives  \begin{multline}\label{derivative_terms}\xi_{B_3}'(0) = A'_f + \frac{2}{5} \fint_{(-2-\varepsilon)} \Gamma(z_2) \zeta(z_2) \zeta(-z_2-1) \Gamma(-z_2-1) (2 z_2+9\times 2^{z_2}+7) dz_2 \\
			+ \frac{5}{6} \fint_{(-2-\varepsilon)} \Gamma(z_1) \zeta(z_1) \zeta(-z_1-1) \Gamma(-z_1-1) (z_1+2^{z_1+3}+4) dz_1 \\
			+ \sum_{i=0}^1 \zeta(-i) \fint_{(-2-\varepsilon)} \Gamma(-z_2)\Gamma(z_2)\zeta(i-z_2)\zeta(z_2) [K^{x_2}_f(0;z_2)][x_1^i] dz_2\\
			+ \sum_{i=0}^1 \zeta(-i) \fint_{(-2-\varepsilon)} \Gamma(-z_1)\Gamma(z_1)\zeta(i-z_1)\zeta(z_1) [K^{x_1}_f(0;z_1)][x_2^i] dz_1 \\
			+ \fint_{(-2-\varepsilon)} \fint_{(-2-\varepsilon)} \Gamma (z_1) \Gamma (z_2) \Gamma (-z_1-z_2)\zeta (z_1) \zeta (z_2) \zeta (-z_1-z_2)(2^{z_1+z_2}+2^{z_2}+1)dz_1 dz_2.
		\end{multline}
		All single integrals above reduce to integrals of the form
		$$\mathcal{I}(k,n,\alpha) = \fint_{(n-1-\varepsilon)} \Gamma(z)\zeta(z) \alpha^{-z} (z)_k \Gamma(n-z)\zeta(n-z) dz, \qquad n\in\{0,-1\},\quad \alpha \in\{1,2,\frac{1}{2}\},\quad k\in \{0,1\},$$
		which can be calculated by the method outlined in Section \ref{convolution_int_section}, giving
		\begin{multline*}\xi'_{B_3}(0) = A'_f + \left(-17 \zeta '(-2)-\frac{71}{12}\zeta'(-1)-\frac{5}{144}-\frac{7 \pi ^2}{1728}+\frac{\gamma }{48}-\frac{25}{24} \log (2)-\frac{13}{72} \log (2 \pi )-\frac{\log (\pi )}{3} \right) \\ + \fint_{(-2-\varepsilon)}\fint_{(-2-\varepsilon)} \Gamma(z_1)\Gamma(z_2)\Gamma(-z_1-z_2)\zeta(z_1)\zeta(z_2)\zeta(-z_1-z_2) (1+2^{z_1} + 2^{z_1+z_2}) dz_1dz_2.
		\end{multline*}
		We will evaluate the last integral, which is $\mathcal{J}(\frac{1}{2},\frac{1}{2}) + \mathcal{J}(\frac12,1)+\mathcal{J}(1,1)$, at the end of this section. So we obtain the evaluation of $\xi'_{B_3}(0)$, note we have again miraculous cancellation of various terms. \par
		For $\xi_{C_3}(s)$, the steps are completely analogous, and we omit them here. \end{example}
	
	\subsection{Double integrals arising from Mellin convolution}
	Now we describe how to calculate double integrals
	$$\mathcal{J}(\alpha_1,\alpha_2) = \fint_{(-2-\varepsilon)}\fint_{(-2-\varepsilon)} \Gamma(z_1)\Gamma(z_2)\Gamma(-z_1-z_2)\zeta(z_1)\zeta(z_2)\zeta(-z_1-z_2) \alpha_1^{-z_1} \alpha_2^{-z_2}dz_1dz_2,$$
	where $\alpha_1,\alpha_2$ are positive rational numbers. This is an extension of the method described in Section \ref{convolution_int_section}. We will only need $\alpha_k \in \{1,2,\frac12\}$. To begin, note that by Mellin convolution,
	$$\int_0^\infty \frac{x^{s-1} dx}{(e^x-1)(e^{\alpha_1 x}-1)(e^{\alpha_2 x}-1)} = \fint_{(c)}\fint_{(c)} \Gamma(z_1)\Gamma(z_2)\Gamma(s-z_1-z_2)\zeta(z_1)\zeta(z_2)\zeta(s-z_1-z_2) \alpha_1^{-z_1} \alpha_2^{-z_2} dz_1dz_2.$$
	for $\Re(s)$ and $c$ sufficiently large. Using partial fractions, the integral on the LHS can be evaluated in terms of Hurwitz zeta functions. For the RHS, we shall shift the contours of integration into a region that allows us to set $s=0$. We recall formula \ref{double-contour-shift}, putting $$p(z_1,z_2) = \Gamma(z_1)\Gamma(z_2)\Gamma(s-z_1-z_2)\zeta(z_1)\zeta(z_2)\zeta(s-z_1-z_2) \alpha_1^{-z_1} \alpha_2^{-z_2}.$$
	Considering its singularity diagram in Figure \ref{pi_factor_sing_diagram}, we can let $P_1 = P_2 = \{1,0,-1,-2\}$, and we have
	\begin{multline*}\fint_{(c)} \fint_{(c)} p(z_1,z_2) dz_1 dz_2 = \sum_{s_1 \in P_1, s_2\in P_2} \Res{z_1 = s_1} \Res{z_2 = s_2} p(z_1,z_2) + \sum_{s_2 \in P_2}\fint_{(-2-\varepsilon)} \Res{z_2 = s_2} p(z_1,z_2) dz_1\\
		+ \sum_{s_1 \in P_1} \fint_{(-2-\varepsilon)}\Res{z_1 = s_1} p(z_1,z_2) dz_2+ \fint_{(-2-\varepsilon)} \fint_{(-2-\varepsilon)} p(z_1,z_2) dz_1 dz_2.\end{multline*}
	The left-hand side can be expressed in terms of Hurwitz zeta and gamma functions, with variable $s$. The same is true of the first term on the right-hand side and, by the procedure described in Section \ref{convolution_int_section}, of the second and third terms. Letting $s\to 0$ therefore gives the value of $$\fint_{(-2-\varepsilon)} \fint_{(-2-\varepsilon)} p(z_1,z_2) dz_1 dz_2,$$ which is $\mathcal{J}(\alpha_1,\alpha_2)$.
	
	\begin{example}
		$$\begin{aligned}
			\mathcal{J}(1,1) &= 3 \zeta '(-2)-\frac{3}{2} \zeta'(-1)+\frac{\pi ^2}{864}+\frac{1}{6}+\frac{\gamma }{12}-\frac{5}{24} \log (2 \pi ),\\
			\mathcal{J}(\frac12,\frac12) &= \frac{13}{4} \zeta '(-2)-\frac{3}{2}\zeta'(-1)+\frac{\pi ^2}{3456}+\frac{1}{8}+\frac{\gamma }{12}+\frac{\log (2)}{12}-\frac{5}{24} \log (2 \pi ), \\
			\mathcal{J}(\frac12,1) = \mathcal{J}(1,\frac12) &= \frac{27}{8} \zeta '(-2)-\frac{3}{2} \zeta'(-1)+\frac{\pi ^2}{1728}+\frac{1}{6}+\frac{5 \gamma }{48}+\frac{3\log (2)}{24}-\frac{11}{48} \log (2 \pi ),\\
			\mathcal{J}(\frac12,2) &= \frac{27}{16}\zeta'(-2)-\frac{9}{4}\zeta'(-1)+\frac{\pi ^2}{864}+\frac{13}{48}+\frac{17 \gamma }{96}-\frac{25}{96} \log (2)-\frac{17 \log (\pi )}{96}-\frac{1}{4} \log\Gamma \left(\frac{1}{4}\right).\\
		\end{aligned}$$
	\end{example}
	
	\subsection{Conjecture on $\xi_\Phi'(0)$ for general root system}
	While computing $\xi'_\Phi(0)$ for rank-$2$ root systems $\Phi$, we observed that these derivatives exhibit an unusually simple form. We have just seen that this simplicity extends to the rank-$3$ case as well. The author conjectures that this pattern is universal for Witten zeta functions.\par
	
	When $\Phi$ has rank $\leq 3$, we have already defined $\xi_\Phi(s)$ in terms of $\zeta_\Phi(s)$. This can be extended to all root systems, with $K_\Phi \in \mathbb{N}$ given in Table \ref{K_Phi_value_table}, 
	$$\xi_\Phi(s) := K_\Phi^{-s} \zeta_\Phi(s).$$

	\begin{table}[h]
		\begin{tabular}{l|l|l|l}
			\multicolumn{1}{c}{$\Phi$} & \multicolumn{1}{c|}{$K_\Phi$}       & \multicolumn{1}{c}{$\Phi$} & \multicolumn{1}{c}{$K_\Phi$}                                                            \\ \hline
			$A_n$                      & $\prod_{k=1}^n k!$                & $G_2$  & $120$                                                                \\
			$B_n, C_n$                 & $\prod_{k=1}^{n} (2k-1)!$        & $F_4$  & $2^{15} \ 3^7 \ 5^4 \ 7^2 \ 11$                                                     \\
			$D_n$                      & $(n-1)! \prod_{k=1}^{n-2} (2k+1)!$ & $E_6$  & $2^{25}\ 3^{10}\ 5^5\ 7^3\ 11^1$                                       \\
			&                                     & $E_7$  & $2^{47}\ 3^{22}\ 5^{10}\ 7^6\ 11^3\ 13^2\ 17^1$                      \\
			&                                     & $E_8$  & $2^{97}\ 3^{47}\ 5^{21}\ 7^{14}\ 11^8\ 13^6\ 17^4\ 19^3\ 23^2\ 29^1$
		\end{tabular}
		\caption{\small Values of $K_\Phi$ for all root systems. For a geometric interpretation of them, see \cite{au2024vanishing} or \cite[Chapter~3]{komori2023theory}.}
		\label{K_Phi_value_table}
	\end{table}
	\vspace{-3mm}
	
	\begin{conjecture}
		We have
		$$\xi_\Phi'(0) \in   \begin{cases}
			\mathbb{Q} \log(2\pi)  \qquad &\text{ for }\Phi = A_n, \\
			\mathbb{Q} \log(2\pi) + \mathbb{Q} \log(2)  \qquad &\text{ for }\Phi = B_n, C_n, D_n,\\
			\mathbb{Q} \log(2\pi) + \mathbb{Q} \log(2) + \mathbb{Q} \log(3) \qquad &\text{ for }\Phi = G_2, F_4, E_6, E_7,\\
			\mathbb{Q} \log(2\pi) + \mathbb{Q} \log(2) + \mathbb{Q} \log(3) + \mathbb{Q} \log(5) \qquad &\text{ for }\Phi = E_8.
		\end{cases}$$
	\end{conjecture}
	For root system $A_n$, we have a more precise formulation.
	\begin{conjecture}For $n\geq 1$, 
		$$\xi_{A_n}(0) = \frac{(-1)^n}{n+1},\qquad \xi_{A_n}'(0) = \frac{(-1)^n n}{2}\log(2\pi).$$
	\end{conjecture}

	\section{Number of representations for rank 2 and 3 simple Lie algebras}
	Let $\mathfrak{g}$ be a finite-dimensional simple Lie algebra over $\mathbb{C}$, and let $r_\mathfrak{g}(n)$ denote the number of isomorphism classes of representations of $\mathfrak{g}$ of dimension $n$. For the corresponding root system $\Phi$, we abuse notation by writing $r_\Phi(n)$ for $r_\mathfrak{g}(n)$. Recall that every representation of $\mathfrak{g}$ is a unique direct sum of irreducible representations; thus 
	$$\sum_{n\geq 0} r_\Phi(n) q^n = \prod_{\rho \text{ irreducible }}(1-q^{\dim \rho})^{-1}.$$
	For the simplest case $A_1$, $\sum_{n\geq 0}r_{A_1}(n)q^n = \prod_{n\geq 1} (1-q^n)^{-1}$ is the well-known partition generating function. For rank-$2$ and rank-$3$ Lie algebras, the generating functions are explicitly
	{\allowdisplaybreaks \begin{align*}\sum_{n\geq 0} r_{A_2}(n) q^n &= \prod_{n,m\geq 1} (1-q^{mn(m+n)/2})^{-1},\\
			\sum_{n\geq 0} r_{B_2}(n) q^n &= \prod_{n,m\geq 1} (1-q^{mn(m+n)(m+2n)/6})^{-1},  \\ 
			\sum_{n\geq 0} r_{G_2}(n) q^n &= \prod_{n,m\geq 1} (1-q^{mn(m+n)(m+2n)(m+3n)(2m+3n)/120})^{-1},  \\
			\sum_{n\geq 0} r_{A_3}(n) q^n &= \prod_{m_k\geq 1} (1-q^{m_1m_2m_3 (m_1+m_2)(m_2+m_3)(m_1+m_2+m_3)/12})^{-1}, \\
			\sum_{n\geq 0} r_{B_3}(n) q^n &= \prod_{m_k\geq 1} (1-q^{m_1m_2m_3 (m_1+m_2)(m_2+m_3)(m_1+m_2+m_3)(2m_2+m_3)(2m_1+2m_2+m_3)(m_1+2m_2+m_3)/720})^{-1}, \\
			\sum_{n\geq 0} r_{C_3}(n) q^n &= \prod_{m_k\geq 1} (1-q^{m_1m_2m_3 (m_1+m_2)(m_2+m_3)(m_1+m_2+m_3)(m_2+2m_3)(m_1+m_2+2m_3)(m_1+2m_2+2m_3)/720})^{-1}.
	\end{align*}}
	
	A combination of circle method and saddle point method, adapted in \cite{bridges2024asymptotic, debruyne2020saddle}, can be used to find asymptotic of $r_\mathfrak{g}(n)$. We briefly recall the main framework and result of \cite{bridges2024asymptotic}. Let $f:\mathbb{N}\to \mathbb{Z}_{\geq 0}$ be an arithmetic function, $$L_f(s) := \sum_{n\geq 1} \frac{f(n)}{n^s},\qquad \sum_{n\geq 0} p_f(n) q^n := \prod_{n\geq 1} \frac{1}{(1-q^n)^{f(n)}}.$$
	The following analytic conditions have to be satisfied:
	\begin{enumerate}[leftmargin=*]
		\item Let $\alpha>0$ be the largest pole of $L_f(s)$, and let $\Lambda = \mathbb{N} \setminus f^{-1}(\{0\})$. There exists $L\in \mathbb{N}$ such that for all primes $p$, we have $| \Lambda \setminus (p\mathbb{N}\cap \Lambda)| \geq L > \frac{\alpha}{2}$.
		\item $L_f(s)$ admits meromorphic continuation to $\mathbb{C}$; all its poles with $\Re(s) > 0$ are real and simple, and $L_f(s)$ is analytic at $s=0$.
		\item $L_f(s)$ is moderately increasing in the imaginary direction.
	\end{enumerate}
	
	Here the last hypothesis (3) means that for every compact subset $K\subset \mathbb{R}$, there exists $M>0$ such that $$|L_f(s)| = O(|\Im(s)|^M), \qquad \Re(s)\in K,\quad |\Im(s)|\gg 1.$$
	This is true for every Witten zeta function, or more generally, for all rank-$2$ and rank-$3$ series $\xi_f(s)$ considered in this article. This follows immediately from Theorem 3 and Hypothesis $H_{0}S$ in Essouabri \cite{essouabri1997singularite}. The second hypothesis\footnote{The article \cite{bridges2024asymptotic} additionally requires all negative poles to be simple; this is not true for $B_3$ and $C_3$ (Theorem \ref{B3_residues}). However, they do not influence the following discussion since we focus only on the asymptotic inside the exponential.} is satisfied for the rank-$2$ and rank-$3$ zeta functions. The first hypothesis is also easy to check (see the end of \cite{rutard2026values}, for example). \par
	
	Let $\mathcal{P}$ be the set of poles of $L_f(s)$ with positive real parts. For $p\in \mathcal{P}$, define \begin{equation}\label{aux_8}r_p := \Gamma(p) \zeta(p+1) \Res{s=p} L_f(s).\end{equation} Let $\alpha := \max \mathcal{P}$ denote the largest pole of $L_f(s)$. Note that $r_\alpha >0$, as $\alpha$ is the abscissa of convergence of a Dirichlet series with non-negative coefficients.
	
	\begin{theorem} (Theorem 1.4 of \cite{bridges2024asymptotic})
		Under hypothesis and notations above, we have $$p_f(n) \sim \frac{C}{n^b} \exp\left( \mathcal{L}(n) \right),\qquad n\to \infty,$$
		where $\mathcal{L}(n)$ is a Puiseux series\footnote{that is, a formal power series in $n$ whose exponents are allowed to be rational numbers} with leading term $$\left(1+\frac{1}{\alpha}\right) (\alpha r_\alpha)^{1/(\alpha+1)} n^{\frac{\alpha}{\alpha+1}},$$
		$C$ and $b$ are constants given by
		$$C = \frac{e^{L_f'(0)}}{\sqrt{2\pi (\alpha+1)}} (\alpha r_\alpha)^{\frac{1/2-L_f(0)}{\alpha+1}}, \qquad b = \frac{1-L_f(0) +\frac{\alpha}{2}}{\alpha+1}.$$
	\end{theorem}
	Although an explicit expression for $\mathcal{L}(n)$ in the general case is almost impossible to obtain, it is easy to describe the procedure for calculating it. We assume that readers know how to invert a Puiseux series with positive leading term. Namely, given a small positive real number $w$ and $z$ related by
	$$w = c_0 z^{r_0} + c_1 z^{r_1} + \ldots + c_{k-1} z^{r_{k-1}}, \qquad z\to 0^{+},\qquad c_0>0,\quad 0<r_0<r_1<\ldots,$$ 
	we can find the series expansion of $z$ in terms of $w$, starting with $$z = c_0^{-1/r_0} w^{1/r_0}+\ldots, \qquad w\to 0^+.$$

	\begin{tcolorbox}[colback=gray!5!white,colframe=gray!75!black,title= Calculating the Puiseux series $\mathcal{L}(n)$]
		Input: Location of positive poles $\mathcal{P}$, and $r_p$ as in equation (\ref{aux_8}). \par
		Output: The Puiseux series $\mathcal{L}(n)$. 
		\tcblower
		\begin{enumerate}[leftmargin=*]
			\item Let $\Phi(z) = \sum_{p\in \mathcal{P}} p r_p z^{-p-1}$, with $z\to 0^+$. Compute the Puiseux-series inversion of $w = 1/\Phi(z)$ up to an $o(w)$ term, where $w\to 0^+$. Denote the resulting series by $z=\Psi(w)$.
			\item Compute $w^{-1} \Psi(w) + \sum_{p\in \mathcal{P}} r_p \Psi(w)^{-p}$ up to an $o(1)$ term.
			\item Substitute $w = 1/n$ into the above series; the result is $\mathcal{L}(n)$.
		\end{enumerate}
	\end{tcolorbox}
	
	\begin{remark}
		(a) The above procedure is implicit in Lemma 3.2 and Theorem 4.1 of \cite{bridges2024asymptotic}, here $\Psi(w)$ is the series expansion of $\varrho_n$ in the paper. \\
		(b) A Mathematica code performing the above step can be found at \url{https://sites.google.com/view/kc-au/2412-17196}. It is used to calculate examples given below.
	\end{remark}
	
	Now we are ready to write down the asymptotics of $$r_\Phi(n) \sim \frac{C}{n^b} \exp\left( \mathcal{L}(n) \right),\qquad n\to \infty,$$ 
	Note that the values of $C$ and $b$ are known, since we found $\zeta_\Phi(0)$ and $\zeta_\Phi'(0)$ in Sections 6 and 10.
	
	\begin{example}
		For $\Phi = A_2$, we have $\mathcal{P} = \{\frac23,\frac12\}$, $\alpha = \max \mathcal{P} = 
		\frac{2}{3}$, $$b=\frac{3}{5},\qquad C = 2^{14/15} \: 3^{3/10} \: 5^{-1/2}\sqrt{\pi }\zeta\left(\frac{5}{3}\right)^{1/10} \Gamma \left(\frac{1}{3}\right)^{1/5},$$
		$$\mathcal{L}(n) =\frac{5 r_{2/3}^{3/5}}{2^{2/5} 3^{3/5}} \colorn{2/5} + \frac{r_{1/2}}{\alpha ^{3/10} r_{2/3}^{3/10}} \colorn{3/10} -\frac{9 r_{1/2}^2}{80 \sqrt[5]{\alpha } r_{2/3}^{6/5}}\colorn{1/5}+\frac{99 r_{1/2}^3}{3200 \sqrt[10]{\alpha } r_{2/3}^{21/10}}\colorn{1/10}-\frac{27 r_{1/2}^4}{2560 r_{2/3}^3}.$$
		Here $r_{2/3}$ and $r_{1/2}$ are
		$$r_p = \Gamma(p) \zeta(p+1) 2^p \Res{s=p} \xi_{A_2}(s),$$
		with values of residues (Theorem \ref{rank2-residues-abiscssa} and Corollary \ref{second_pole_residue}):
		$$\Res{s=2/3} \xi_{A_2}(s) = \frac{\Gamma(\frac{1}{3})^3}{2 \sqrt{3} \pi },\qquad \Res{s=1/2}\xi_{A_2}(s) = \zeta\left(\frac{1}{2}\right).$$
	\end{example}

	\begin{example}
		For $\Phi = B_2$, we have $\mathcal{P} = \{\frac12,\frac13\}$, $\alpha = \max \mathcal{P} = 
		\frac{1}{2}$,
		$$b=\frac{7}{12},\qquad C = 2^{5/4} \: 3^{-1/8} \pi \Gamma(\frac14)^{1/6} \zeta(\frac{3}{2})^{1/12},$$
		$$\mathcal{L}(n) =3 \alpha ^{2/3} r_{1/2}^{2/3} \colorn{1/3} + \frac{r_{1/3}}{\alpha ^{2/9} r_{1/2}^{2/9}}\colorn{2/9}-\frac{2 r_{1/3}^2}{27 \sqrt[9]{\alpha } r_{1/2}^{10/9}}\colorn{1/9}+\frac{4 r_{1/3}^3}{243 r_{1/2}^2}.$$
		Here $r_{1/2}$ and $r_{1/3}$ are
		$$r_p = \Gamma(p) \zeta(p+1) 6^p \Res{s=p} \xi_{B_2}(s),$$
		with values of residues (Theorem \ref{rank2-residues-abiscssa} and Corollary \ref{second_pole_residue}):
		$$\Res{s=1/2} \xi_{B_2}(s) = \frac{\Gamma \left(\frac{1}{4}\right)^2}{8 \sqrt{2 \pi }},\qquad \Res{s=1/3} \xi_{B_2}(s) = \frac{1+2^{-1/3}}{3}\zeta\left(\frac{1}{3}\right).$$
	\end{example}

	\begin{example}
		For $\Phi = G_2$, we have $\mathcal{P} = \{\frac13,\frac15\}$, $\alpha = \max \mathcal{P} = 
		\frac{1}{3}$,
		$$b = \frac{9}{16},\qquad C = 2^{17/12} \: 3^{25/96} \: 5^{5/12} \pi ^{31/16} \Gamma(\frac{1}{3})^{1/4} \zeta(\frac{4}{3})^{1/16},$$
		$$\mathcal{L}(n) =4 \alpha ^{3/4} r_{1/3}^{3/4} \colorn{1/4} + \frac{r_{1/5}}{\alpha ^{3/20} r_{1/3}^{3/20}}\colorn{3/20}-\frac{9 r_{1/5}^2}{200 \sqrt[20]{\alpha } r_{1/3}^{21/20}} \colorn{1/20}.$$
		Here $r_{1/3}$ and $r_{1/5}$ are 
		$$r_p = \Gamma(p) \zeta(p+1) 120^p \Res{s=p} \xi_{G_2}(s),$$
		with values of residues (Theorem \ref{rank2-residues-abiscssa} and Corollary \ref{second_pole_residue}):
		$$\Res{s=1/3} \xi_{G_2}(s) = \frac{\Gamma \left(\frac{1}{3}\right)^3}{2^{8/3} \: 3^{3/2} \pi },\qquad \Res{s=1/5} \xi_{G_2}(s) = \frac{2^{-1/5}+18^{-1/5}}{5}\zeta\left(\frac{1}{5}\right).$$
	\end{example}
	
	\begin{example}
		For $\Phi = A_3$, we have $\mathcal{P} = \{\frac12, \frac25,\frac13,\frac14\}$, $\alpha = \max \mathcal{P} = \frac{1}{2}$,
		$$b = 1,\qquad C = \frac{\sqrt{\zeta \left(\frac{3}{2}\right)} \Gamma \left(\frac{1}{4}\right)^2}{48 \sqrt[4]{3} \pi ^{9/4}},$$
		\begin{multline*}\mathcal{L}(n) = 3 \alpha ^{2/3} r_{1/2}^{2/3}\colorn{1/3}
			+\frac{r_{2/5}}{\alpha ^{4/15} r_{1/2}^{4/15}}\colorn{4/15}
			+\frac{r_{1/3}}{\alpha ^{2/9} r_{1/2}^{2/9}}\colorn{2/9}
			-\frac{8 r_{2/5}^2}{75 \sqrt[5]{\alpha } r_{1/2}^{6/5}} \colorn{1/5} 
			+\frac{r_{1/4}}{\sqrt[6]{\alpha } \sqrt[6]{r_{1/2}}} \colorn{1/6} 
			-\frac{8 r_{1/3} r_{2/5}}{45 \alpha ^{7/45} r_{1/2}^{52/45}} \colorn{7/45} \\
			+\frac{544 r_{2/5}^3}{16875 \alpha ^{2/15} r_{1/2}^{32/15}} \colorn{2/15}
			-\frac{2 r_{1/3}^2}{27 \sqrt[9]{\alpha } r_{1/2}^{10/9}} \colorn{1/9}
			-\frac{2 r_{1/4} r_{2/5}}{15 \sqrt[10]{\alpha } r_{1/2}^{11/10}} \colorn{1/10}
			+\frac{784 r_{1/3} r_{2/5}^2}{10125 \alpha ^{4/45} r_{1/2}^{94/45}} \colorn{4/45} \\
			-\frac{15872 r_{2/5}^4}{1265625 \sqrt[15]{\alpha } r_{1/2}^{46/15}} \colorn{1/15} 
			-\frac{r_{1/4} r_{1/3}}{9 \sqrt[18]{\alpha } r_{1/2}^{19/18}} \colorn{1/18}
			+\frac{376 r_{1/3}^2 r_{2/5}}{6075 \alpha ^{2/45} r_{1/2}^{92/45}} \colorn{2/45}
			+\frac{62 r_{1/4} r_{2/5}^2}{1125 \sqrt[30]{\alpha } r_{1/2}^{61/30}} \colorn{1/30} \\
			-\frac{267904 r_{1/3} r_{2/5}^3}{6834375 \sqrt[45]{\alpha } r_{1/2}^{136/45}} \colorn{1/45}  +\left(\frac{-r_{1/4}^2}{24 r_{1/2}}+\frac{4 r_{1/3}^3}{243 r_{1/2}^2}+\frac{256 r_{2/5}^5}{46875 r_{1/2}^4} \right).\end{multline*}
		The $r_p$'s are
		$$r_p = \zeta(p+1)\Gamma(p) 12^p  \Res{s=p} \xi_{A_3}(s),$$
		with values of residues given in Theorem \ref{A3_residues}.
	\end{example}
	
	\begin{example}
		For $\Phi = B_3$ or $C_3$, we have $\mathcal{P} = \{\frac13,\frac14,\frac15,\frac16,\frac17\}$, $\alpha = \max \mathcal{P} = \frac{1}{3}$,
		$$b = \frac{71}{64},\qquad C = \frac{\zeta \left(\frac{4}{3}\right)^{39/64} \Gamma \left(\frac{1}{3}\right)^{273/64}}{384\cdot 2^\bullet \: 3^{11/32} \: 5^{5/16} \pi ^{145/32}},$$
		where $\bullet = \frac{67}{64}$ or $\frac{63}{64}$ according as $\Phi = B_3$ or $C_3$, respectively, and
		\begin{multline*}\mathcal{L}(n) = 4 \alpha ^{3/4} r_{1/3}^{3/4} \colorn{1/4}
			+\frac{r_{1/4}}{\alpha ^{3/16} r_{1/3}^{3/16}} \colorn{3/16}
			+\frac{r_{1/5}}{\alpha ^{3/20} r_{1/3}^{3/20}} \colorn{3/20}
			+\frac{128 r_{1/6} r_{1/3}-9 r_{1/4}^2}{128 \sqrt[8]{\alpha } r_{1/3}^{9/8}} \colorn{1/8} 
			+ \frac{r_{1/7}}{\alpha ^{3/28} r_{1/3}^{3/28}} \colorn{3/28} \\
			-\frac{9 r_{1/5} r_{1/4}}{80 \alpha ^{7/80} r_{1/3}^{87/80}} \colorn{7/80}
			+\frac{3 \left(51 r_{1/4}^3-256 r_{1/6} r_{1/4} r_{1/3}\right)}{8192 \sqrt[16]{\alpha } r_{1/3}^{33/16}}\colorn{1/16}
			-\frac{9 r_{1/5}^2}{200 \sqrt[20]{\alpha } r_{1/3}^{21/20}} \colorn{1/20} 
			-\frac{9 r_{1/7} r_{1/4}}{112 \alpha ^{5/112} r_{1/3}^{117/112}} \colorn{5/112}\\
			+\frac{3 r_{1/5} \left(369 r_{1/4}^2-640 r_{1/6} r_{1/3}\right)}{25600 \sqrt[40]{\alpha } r_{1/3}^{81/40}} \colorn{1/40}
			-\frac{9 r_{1/7} r_{1/5}}{140 \sqrt[140]{\alpha } r_{1/3}^{141/140}} \colorn{1/140} 
			-\frac{27 r_{1/4}^4-144 r_{1/6} r_{1/3} r_{1/4}^2+128 r_{1/6}^2 r_{1/3}^2}{4096 r_{1/3}^3}.\end{multline*}
		The $r_p$'s are
		$$r_p = \zeta(p+1)\Gamma(p) 720^p  \Res{s=p} \xi_{\Phi}(s),$$
		with values of residues given in Theorems \ref{B3_residues} or \ref{C3_residues}.
	\end{example}
	The $A_2$ example above was obtained by Romik in \cite{romik2017number}, and the $B_2$ example in \cite{bridges2024asymptotic}. For $G_2$, Rutard \cite{rutard2026values} wrote down the same expression for $\mathcal{L}(n)$ with unevaluated residues. The rank-$3$ results are new.
	
	\bibliographystyle{plain} 
	\bibliography{ref.bib} 
	
\end{document}